\documentclass{amsart}
\usepackage[dvips]{graphicx}
\usepackage{amssymb}
\usepackage{latexsym}
\usepackage{xypic}
\usepackage{multirow}

\newtheorem{theorem}{Theorem}[section]
\newtheorem{lemma}[theorem]{Lemma}
\newtheorem{proposition}[theorem]{Proposition}

\newtheorem{corollary}[theorem]{Corollary}
\newtheorem{Thm}{Theorem}
\newtheorem{Cor}{Corollary}

\theoremstyle{definition}
\newtheorem{definition}[theorem]{Definition}
\newtheorem{example}[theorem]{Example}

\theoremstyle{remark}
\newtheorem{remark}[theorem]{Remark}

\newenvironment{proaf}{{\noindent \sc Proof} }{\mbox{ }\hfill$\Box$                \vspace{1.5ex}\par}

\newcommand{\C}{\mathbb{C}}
\newcommand{\R}{\mathbb{R}}
\newcommand{\s}{\mathbb{S}}
\newcommand{\PC}{{\mathbb P}_{\C}}
\newcommand{\vv}{\vskip.1cm}
\newcommand{\PSL}{\rm{PSL}}
\newcommand{\SO}{\rm{SO}}
\newcommand{\SL}{\rm{SL}}
\newcommand{\PU}{\rm{PU}}
\newcommand{\Aff}{\rm{Aff}}
\newcommand{\Kul}{\rm{Kul}}
\newcommand{\Gr}{\rm{Gr}}
\newcommand{\Eq}{\rm{Eq}}

\numberwithin{equation}{section}

\begin{document}
\author{Angel Cano and Jos\'e Seade}

\title{  \sc{On Discrete groups of automorphisms of $\mathbb{P}^2_{\mathbb{C}}$} }

  \thanks{Partially supported by  CONACYT  and PAPIIT-UNAM, Mexico, and  CNPq, Brazil} 





\begin{abstract}
We study  the geometry and dynamics of discrete subgroups $\Gamma$  of
$\PSL(3,\mathbb{C})$ with
an open invariant set $\Omega \subset \PC^2$ where the action is properly discontinuous and the quotient $\Omega/\Gamma$ contains a connected component whicis compact.
We call such groups {\it quasi-cocompact}. 
In this case $\Omega/\Gamma$ is a compact complex projective orbifold and $\Omega$ is a {\it divisible set}.
Our first theorem refines classical work by Kobayashi-Ochiai and others about complex surfaces with a projective structure: 
We prove that  every such group is either
virtually  affine or complex hyperbolic.  
We then classify the  divisible sets  that  appear in this way, the corresponding quasi-cocompact groups and
 the   orbifolds $\Omega/\Gamma$.   We also
  prove  that excluding  
 a few exceptional cases,  the   Kulkarni region of discontinuity  coincides with the equicontinuity region and is the largest
  open invariant set where the action is properly discontinuous.
\end{abstract}

\maketitle

\section*{Introduction}

Classical Kleinian groups were introduced by H. Poincar\'e at the end of the 19th Century and their
study has played a major role in several areas of mathematics.   These are discrete subgroups of
$\PSL(2,\C)$ that act on the Riemann sphere $\s^2 \cong
\mathbb{P}^1_{\mathbb{C}}$ with non-empty region of discontinuity. 

  Such a group $\Gamma$ determines a splitting of the Riemann sphere  as a union of two $\Gamma$-invariant sets, $\s^2 = \Omega \cup \Lambda$,   where $\Omega$ is the 
 discontinuity region  and $\Lambda$ is its limit set.
   The limit set is where the dynamics concentrates and its study was for decades the paradigm of holomorphic dynamics, as enlightened by the Sullivan-McMullen dictionary between Kleinian groups and rational maps. On the other hand, the quotient $\Omega/\Gamma$ is a Riemann surface equipped with a rich geometry. In particular $\Omega/\Gamma$ has the structure of a complex projective orbifold. The study of the Riemann surfaces one gets in this way is a fascinating subject that somehow begins with B. Riemann and P. K\"obe, and passes through the work of many authors as for instance L. Ahlfors, D.  Sullivan and W. Thurston, to name a few.

More recently, there has been great interest in studying
generalizations of these groups to higher dimensional complex projective
spaces. Of particular interest are the results of P. Deligne, W. Goldman, N. Gusevskii , G. D. Mostow, J. Parker, R.
Schwartz and others, about discrete
subgroups of $\PU(n,1) \subset \PSL(n+1,\mathbb{C})$, the group of
holomorphic isometries of the complex hyperbolic space $\mathbb H^n_\C$.

In \cite {sv1} the authors introduce the concept
of a {\it complex Kleinian group}, which means a discrete subgroup
$\Gamma$ of some $\PSL(n+1,\mathbb{C})$ acting on
$\mathbb{P}^n_{\mathbb{C}}$ so that there is a non-empty open invariant set $\Omega$  where the action is properly discontinuous.  And in a couple of
subsequent articles (\cite {sv2, sv3}) they study some interesting
families of such groups acting on odd-dimensional projective
spaces. Here we look at  the case of groups acting on
$\mathbb{P}^2_{\mathbb{C}}$, continuing the work begun in \cite {cano,
pablo1, pablo2}.

In that setting, of a discrete subgroup
$\Gamma$ of  $\PSL(3,\mathbb{C})$ acting properly discontinuously on an open subset $\Omega \subset \PC^2$,
 the quotient $\Omega/\Gamma$ is a  2-dimensional complex projective orbifold, which may or may not be compact. 
 If   there exists such an invariant set $\Omega$ whose quotient $\Omega/\Gamma$ containas a connected component which is compact, then we say that the group $\Gamma$ is   {\it quasi-cocompact} and, following Y. Benoist \cite{Benoist}, we call $\Omega$  a divisible set. The study of quasi-cocompact complex Kleinian groups is the subject of this article.

Notice that $\PSL(3,\C)$ contains as subgroups the affine group $\Aff(\C^2)$ and  the group $\PU(2,1)$ of holomorphic isometries of the complex hyperbolic 2-space. Hence every 
$(\C^2,\Aff(\C^2))$-orbifold and every  $(\mathbb H^2_\C,\PU(2,1))$-orbifold, is also a complex projective orbifold. Thus our work is naturally linked with the classical uniformization problem for compact two-dimensional complex manifolds and orbifolds.

We say that a subgroup $\Gamma\subset \PSL(3,\mathbb{C})$ is {\it complex hyperbolic} if it is conjugate in $\PSL(3,\mathbb{C})$ to a subgroup of $\PU(2,1)$. The group is {\it affine} if it is conjugate to a subgroup of $\Aff(\C^2)$. And 
the group is {\it virtually affine} it has a finite index subgroup which is affine.

\vv
Our first theorem is:

\begin{Thm} \label{t:main1}
Let $\Gamma\subset \PSL(3,\mathbb{C})$ be  an infinite quasi-cocompact group. Then $\Gamma$ is   complex hyperbolic or virtually affine.
\end{Thm}

 In  Section \ref{s:tor} we construct an example of a
 quasi-cocompact
group which is virtually affine but  is not affine nor complex hyperbolic, showing that this theorem fails if we drop from it the word ``virtually".  And in section \ref{s:sch} we provide an example  showing that the theorem  fails
 if we drop the condition of being quasi-cocompact: We construct a ``kissing Schottky group",  which is complex Kleinian but  is not  complex hyperbolic nor virtually affine.
 
\vv

Theorem \ref{t:main1} is obviously inspired by the work of Kobayashi and Ochiai in \cite{KO}, where they prove that every compact complex surface with a projective structure is either affine or complex hyperbolic. In fact, results of  Kobayashi, Ochiai and  Inoue \cite{KO, KOI}, Klingler \cite{klingler, klingler2}, Mok,  and Yeung   \cite {Mok1, Mok2, MY} can be put together (in Theorem \ref{t:klingler} below)  to assert that  for each
compact complex surface which admits a projective structure, the
holonomy is either  complex hyperbolic  or  affine, and   the image of the
developing map  is one of the  following 8 complex manifolds: $\mathbb{P}^2_{\mathbb{C}}$,
 $\mathbb{C}^2$, $\mathbb{C}^2\setminus\{0\}$, 
$\mathbb{C}\times \mathbb{C}^*$, $\mathbb{C}^*\times\mathbb{C}^*$, $\mathbb{H}^2_{\mathbb{C}}$,
$\mathbb{H}\times \mathbb{C}$ and $D\times \mathbb{C}^*$, where $\mathbb{H}$ is a half-plane in $\C$ and $D$ is
a hyperbolic domain of the extended complex plane $\widehat \C \cong \mathbb{P}^1_{\mathbb{C}}$.  This is our starting point for proving 
Theorem \ref{t:main1}. 

\vskip.1cm

An important problem in the theory  of  discrete group actions 
 is  providing a ``nice" definition of  the limit set. The point is that for all $n \ge 2$, there are examples of discrete subgroups of $\PSL(n+1,\C)$  acting on $\PC^n$ so that if we take as limit set the usual one, the set of accumulation points of the orbits, then the action is not properly discontinuous on the complement. This question was  addressed by 
 R. S. Kulkarni  in \cite {kulkarni} in a rather general setting. He introduced 
 an interesting notion of  a  ``limit set", that we call the Kulkarni limit set  and denote it $\Lambda_{\Kul}(\Gamma)$,
 which  has the nice property of granting that the action of  $\Gamma$ on  its complement is properly discontinuous. 
 The complement 
 $\Omega_{\Kul}(\Gamma) := \PC^2 \setminus \Lambda_{\Kul}(\Gamma)$ is the Kulkarni region of
discontinuity of $\Gamma$.

\vv

Though the action of  $\Gamma$ on  $\Omega_{\Kul}(\Gamma) $ is properly discontinuous, unlike the classical case of Kleinian subgroups of $\PSL(2,\C)$, it is not truth  that this set is always the largest open invariant set where the action is properly discontinuous (see \cite {pablo1, CNS} for explicit examples, and see also  Corollary
\ref{c:maxc} below). Neither is truth that $\Omega_{\Kul}(\Gamma) $ necessarily coincides  with the region of equicontinuity $\Eq(\Gamma)$  (see for instance
Corollary \ref{r:equi}).   Yet, the Kulkarni sets $\Omega_{\Kul}(\Gamma)$ and $\Lambda_{\Kul}(\Gamma)$ do pick a lot of information about the group action, and they have very interesting geometric and dynamical properties which are used in the sequel. Moreover,
there seems to be evidence that  ``generically", in dimension two  one has $\Omega_{\Kul}(\Gamma)= \Eq(\Gamma)$ and this is the largest set where the action is properly discontinuous (cf. \cite{pablo1}). In fact  we prove:

\begin{Thm}\label{t:main2}
Let $\Gamma\subset \PSL(3,\mathbb{C})$ be a quasi-cocompact group which is not virtually cyclic, then:

\begin{enumerate}

\item The set  $\Omega_{\Kul}(\Gamma)$ is the largest open
set on which $\Gamma$ acts properly discontinuously.

\item If $\Gamma$ is not a finite extencion of the fundamental group of a primary Kodaira surface nor 
 a finite extension  of the fundamental group of an Inoue surface, then $\Omega_{\Kul}(\Gamma)$ equals the equicontinuity region $\Eq(\Gamma)$.

\end{enumerate}

\end{Thm}

We remark that unlike the 1-dimensional case, in higher dimensions it is usual to have discrete subgroups of $\PSL(n+1,\C)$ with several  regions in $\PC^n$ where the action is properly discontinuous  but there is not  one such region which is the largest.  We also remark that Theorem \ref{c:maxc} below describes 
  the maximal domains of
discontinuity for  cyclic groups;  this is based on \cite{pablo2}. 
In a subsequent paper we shall prove that Theorem \ref{t:main2} extends to subgroups of $\PSL(3,\C)$, not necessarily quasi-cocompact, whose Kulkarni limit set has ``enough" lines in general position.
\vv

Our next theorem classifies the open subsets of $\PC^2$ that can appear as the Kulkarni region of discontinuity of some quasi-cocompact complex Kleinian group:

\begin{Thm}\label{t:main3}
Let $\Gamma\subset \PSL(3,\mathbb{C})$ be  an infinite quasi-cocompact group. Then,
up to projective equivalence,  $\Omega_{\Kul}(\Gamma)$ is either  $\Bbb{C}^2,\,\Bbb{C}^2\setminus \{0\},
\Bbb{C}^*\times \Bbb{C},\,  \Bbb{C}^*\times \Bbb{C}^*,\,  \Bbb{C}^*\times (\Bbb{H}^+\cup \Bbb{H}^-),\, D\times  \Bbb{C}^*$ or $\Bbb{H}^2_{\Bbb{C}}$, where $\Bbb{H}^{\pm}$ denotes the upper and lower half spaces, respectively, and $D$ is a hyperbolic open set in $\Bbb{P}^1_{\Bbb{C}}$.   
\end{Thm}

By  a hyperbolic  set in $\Bbb{P}^1_{\C}$ we mean a subset whose connected components have the hyperbolic plane $\mathbb H^2$ as universal cover.

\vv

To state our next theorem we ought to say first  a few words about compact complex surfaces. Recall first that the Enriques-Kodaira classification splits these into 
 10 types:  rational, ruled (genus  $> 0$), of ``type VII", K3, Enriques, Kodaira, toric, hyperelliptic, properly quasi-elliptic, and of ``general type". {\it Class VII}  means that they have first Betti number $b_1=1$ and Kodaira dimension $-\infty$. 
{\it Kodaira surfaces}  are compact, non-algebraic,  complex surfaces of Kodaira dimension 0 and odd first Betti number, and such a
  surface is called  {\it primary}  if it further has  trivial canonical bundle. A
{\it Hopf surface} means a quotient of $\C^2 \setminus \{0\}$ by a free action of a discrete group. The surface is called {\it primary} if
 its fundamental group is isomorphic to $\mathbb Z$; each of these is  diffeomorphic to $\s^3 \times \s^1$.
All Hopf surfaces are  of class VII.

We also consider {\it elliptic surfaces}, this means a surface that has an elliptic fibration, {\it i.e.}, a proper connected morphism onto an algebraic curve, almost all of whose fibers are elliptic curves.
Another important class of surfaces of  class VII are the {\it Inoue surfaces}. These  come in 
  three families: $S^0$, $S^+$ and $S^-$. They are all compact quotients of ${\Bbb H} \times \mathbb C$ by a solvable discrete group which acts holomorphically on this product space.
 In Theorem (\ref{t:sol40} ) we give a description of the  Inoue surfaces due to C. T. C. Wall \cite{wall} that we use in the sequel.

\vv

Now we may describe
 the type of  projective orbifolds  one gets  as quotients of quasi-cocompact complex Kleinian groups:

\begin{Thm}\label{t:main4}
Let $\Gamma\subset \PSL(3,\mathbb{C})$ be   an infinite  quasi-cocompact group. Then 
 there exists a   finite covering  
$S_\Gamma \twoheadrightarrow \Omega_{\Kul}(\Gamma)/\Gamma$, ramified at the points in $\Omega_{\Kul}(\Gamma)/\Gamma$ with non-trivial isotropy, where the group of the covering is
 of the form  $\Gamma/\Gamma_0$ with $\Gamma_0$ being a finite index, torsion free normal subgroup of $\Gamma$, and the surface  $S_\Gamma$ is of the following type:

\begin{enumerate}
\item If $\Omega_{\Kul}(\Gamma)=\Bbb{C}^2$, then 
$S_\Gamma$ is  biholomorphic to a complex torus $\s^1 \times \s^1  \times \s^1 \times \s^1$ or a
primary  Kodaira surface.

\item If $\Omega_{\Kul}(\Gamma)= \Bbb{C}^2\setminus \{0\}$,  then $S_\Gamma$ is  biholomorphic to a complex torus or a primary Hopf surface.

\item  If $\Omega_{\Kul}(\Gamma) =\Bbb{C}^*\times \Bbb{C}$,  then $S_\Gamma$ is  biholomorphic  to a complex torus.

\item  If $\Omega_{\Kul}(\Gamma) =\Bbb{C}^*\times \Bbb{C}^*$,  then $S_\Gamma$ is  biholomorphic  to a complex torus.

\item  If $\Omega_{\Kul}(\Gamma) = \Bbb{C}^*\times (\Bbb{H}^+\cup \Bbb{H}^-)$,  then $S_\Gamma$ is  either $M$ or $M\sqcup M$ where $M$ is a  Inoue surface and the map $S_\Gamma \twoheadrightarrow \Omega_{\Kul}(\Gamma)/\Gamma$ is a covering with no ramification points.

\item  If $\Omega_{\Kul}(\Gamma) =  D\times  \Bbb{C}^*$,  then $S_\Gamma$ has countably many components with at least one of them being compact, and each connected component of $S_\Gamma$ is 
 an  elliptic surface with an affine structure.

\item  If $\Omega_{\Kul}(\Gamma) = \Bbb{H}^2_{\Bbb{C}}$,  then  $S_\Gamma$ is   a compact complex hyperbolic manifold.
\end{enumerate}

\end{Thm}

Finally, we  describe the type of groups one gets in each case. For this we need to introduce some subgroups of the affine group regarded as: 
$$\,\Aff(\mathbb{C}^2)=\{g\in \PSL(3,\mathbb{C}):  g(\mathbb{C}^2)=\mathbb{C}^2\}\,,$$
where $ \mathbb{C}^2=\{[z:w:1]: z,w\in \mathbb{C}\}$.
A  group  
 $\Gamma \subset \PSL(3,\mathbb{C})$  is  {\it affine}
if it is conjugate to a subgroup of $\Aff(\mathbb{C}^2)$.

We say that an affine group  $\Gamma$ is 
\textit{controllable}  if there exist a line $\ell$ and a point
$p\notin\ell $ which are invariant under the action of $\Gamma$.
The group   $\Gamma\vert_{\ell}=G$ is called {\it the control group} and
$K=\{h\in \Gamma:h(x)=x \textrm{ for all }x\in \ell\}$ is the
\textit{kernel} of $\Gamma$ (see Subsection \ref{ss:cont} for 
examples).

We now define eight different subgroups of  $\Aff(\mathbb{C}^2)$ that appear in the sequel. First we have:
\begin{small}
\begin{equation} \label{e:formasol}
\begin{array}{c}
Aut(\mathbb{C}\times\mathbb{C}^*)=\{g\in \PSL(3,\mathbb{C}):g(\mathbb{C}\times\mathbb{C}^*)=\mathbb{C}\times\mathbb{C}^*\};\\\\
Aut(\mathbb{C}^*\times\mathbb{C}^*)=\{g \in \PSL(3,\mathbb{C}):g(\mathbb{C}^*\times\mathbb{C}^*)=\mathbb{C}^*\times\mathbb{C}^*\}.
\\
\end{array}
\end{equation}
\end{small}
\noindent
These can be regarded as the subgroups of $\PSL(3,\C)$ of maps that leave invariant two lines and three lines, respectively. Now consider the solvable subgroups:

\begin{small}
\begin{equation}\label{e:formasol2}
\begin{array}{c}
Sol^4_0=\left \{\left(\begin{array}{lll} \lambda & 0 & a\\ 0 &
\vert\lambda\vert^{-2} & b\\ 0 & 0 & 1 \end{array}\right): (\lambda, a,
b)\in \mathbb{C}^*\times \mathbb{C}\times \mathbb{R}\right \};\\\\
Sol^4_1=\left \{\left (\begin{array}{lll} \varepsilon  & a & b\\ 0 &
\alpha & c\\ 0 & 0 & 1 \end{array}\right): \alpha , a, b,c\in
\mathbb{R}, \alpha>0, \varepsilon=\pm 1\right \};
\\\\
Sol^{\prime\,4}_1=\left \{\left( \begin{array}{lll} 1  & a & b+i\,log \alpha\\ 0 & \alpha & c\\ 0 & 0 & 1 \end{array}\right): \alpha , a, b,c\in \mathbb{R}, \alpha>0 \right\}.
\end{array}
\end{equation}
\end{small}
These appear in the classification of the Inoue surfaces. Finally consider the groups: 
\begin{small}
\begin{equation}  \label{e:formasol3}
\begin{array}{c}
A_1=\left \{ \left (\begin{array}{lll} 1 & 0 & b\\ 0 & a & 0\\ 0 & 0
& 1 \end{array}\right ): (a, b)\in \mathbb{C}^*\times
\mathbb{C}\right\} \;  ; \quad
 A_2=\left \{\left (\begin{array}{lll} a & b & 0\\ 0 & a & 0\\ 0 &
0 & 1 \end{array}\right ): (a, b)\in \mathbb{C}^*\times \mathbb{C}
\right\} \;.
\\\\

\end{array}
\end{equation}
\end{small}
The first  of these is a direct product $\mathbb{C}^*\times \mathbb{C}$ while   the second is a semi-direct product. 

\vv

We prove:

\begin{Thm}\label{t:main5}

Let $\Gamma\subset \PSL(3,\mathbb{C})$ be an infinite quasi-cocompact group which is not virtually cyclic.

\begin{enumerate}

\item If $\Omega_{\Kul}(\Gamma)=\Bbb{C}^2$, then $\Gamma$ is affine and it is a finite extension of a unipotent group.

\item If $\Omega_{\Kul}(\Gamma) =\Bbb{C}^*\times \Bbb{C}$, then $\Gamma$ is a finite extension of a group isomorphic to $\mathbb{Z}\oplus \mathbb{Z}\oplus\mathbb{Z}$, which up to a projective conjugation is contained in either $A_1$ or
$A_2$.

\item If $\Omega_{\Kul}(\Gamma) =\Bbb{C}^*\times \Bbb{C}^*$, then $\Gamma$ is a finite extension of a group isomorphic to $\mathbb{Z}\oplus \Bbb{Z}$ which is
contained in $Diag(3,\Bbb{C})$ up to projective conjugation, where the latter is
the group of diagonal matrices.

\item If $\Omega_{\Kul}(\Gamma) = \Bbb{C}^*\times (\Bbb{H}^+\cup \Bbb{H}^-)$, then $\Gamma$ is a finite extension of the fundamental group
of an Inoue Surface. In particular $\Gamma$ belongs  to either  $Sol^4_0, Sol^4_1$ or $Sol^{\prime 4}_1$.
\item If $\Omega_{\Kul}(\Gamma) = D\times \Bbb{C}^*$, then $\Gamma$ is
controllable with
infinite kernel and quasi-cocompact control group $\Sigma$ such that $\Omega(\Sigma)=D$.

\item If $\Omega_{\Kul}(\Gamma) = \Bbb{H}^2_{\Bbb{C}}$, then $\Gamma$ is contained in $\PU(2,1)$ up to projective conjugation.

\end{enumerate}

\end{Thm}

To state our last theorem we recall from \cite{pablo2} (see also \cite{CNS}) that a complex  homothety is an element in $\PSL(3,\C)$ that has a lift to $\SL(3,\C)$ whose normal Jordan  form is of the type
$$\begin{pmatrix}
a & 0 & 0\\
0 & a & 0\\
0 & 0 & a^{-2} 
\end{pmatrix} \;,$$
where $\vert a\vert\neq 1 $. We have:

\begin{Thm}\label{t:main6}
Let $\Gamma\subset PSL(3,\C)$ be a group. Then
$\Gamma$ is  virtually cyclic and  quasi-cocompact  if and only if  $\Gamma$
acts properly discontinuously on some  $\Gamma$-invariant {\rm (connected)} domain $\Omega$ such that
    $\Omega/\Gamma$  has a possibly ramified  finite covering which is  a  Hopf surface. Furthermore, in this case one has that the sets $\Omega_{Kul}(\Gamma)=Eq(\Gamma)$ coincide, and the following statements are equivalent:
    \begin{enumerate}
\item \label{i:h3} The orbifold $\Omega_{Kul}(\Gamma)/\Gamma$ is compact.
\item \label{i:h4} $\Gamma$ contains a complex homothety.
\item \label{i:h2} The set $\Omega_{Kul}(\Gamma)$ is the largest open set on which $\Gamma$ acts properly discontinuously.
\end{enumerate}
\end{Thm}

The following are two immediate consequences of these theorems:


\begin{Cor}
Let $\Gamma \subset \PSL(3,\C)$ be quasi-cocompact.   If no component of 
 $\Omega_{\Kul}/\Gamma$  can be covered by a compact elliptic surface, then  $\Omega_{\Kul}/\Gamma$ is compact and:
 
\begin{enumerate}
 \item If  $\Gamma$ is not the fundamental group of an Inoue surface, then $\Omega_{\Kul}/\Gamma$ is connected.
 \item If  $\Gamma$ is  the fundamental group of an Inoue surface, then $\Omega_{\Kul}/\Gamma$ can be either connected, or else consists of 
 two biholomorphic  connected components. 

    \end{enumerate}
\end{Cor}


\begin{Cor}
There is not an analogue
of Bers simultaneous uniformization theorem for groups of $\PSL(3,\Bbb{C})$ acting on $\Bbb{P}^2_{\Bbb{C}}$.
\end{Cor}

Recall that this theorem of  Bers says that if we have two distinct complex structures on a compact Riemann surface of genus $\ge 2$, then their union can be 
 simultaneously uniformized  by a group of M\"obius transformations (with real coefficients). 
\vv

The idea of the proof of Theorem \ref{t:main1} is the following. If $\Gamma \subset \PSL(3,\C)$ is quasi-cocompact, then it acts cocompactly on some open invariant region $\Omega$. The quotient $M:=\Omega/\Gamma$ is a compact orbifold with a projective structure and $\Omega$ is a {\it divisible set}.
Then from the work of  Kobayashi, Ochiai, Klingler  and others we get the classification of  the  domains that appear as divisible sets in $\PC^2$, these are (Theorem \ref{c:coneccion}):  $\mathbb H^2_{\C}$,
$\C^2$, $\C^2 \setminus \{0\}$, $\C\times \C^*$,  $\C^*\times \C^*$,  $\C\times \mathbb H$ or 
 $\C^* \times D$,  where $D$ is a hyperbolic domain in $\PC^1$. 
 Furthermore (Theorem \ref{t:pi_1-orbifolds}), one has that there is always a developing pair $(\mathcal D, \mathcal H)$ for $M$ such that  
 the image of the orbifold fundamental group of $M$ under the holonomy is the isotropy group $Iso(\Omega,\Gamma)$. 
 
 It is then easy to show that when $\Omega$ is $\mathbb H^2_{\C}$ the group $\Gamma$ is complex hyperbolic, and 
 when  $\Omega$  is $ \C^2$, $\C^2 \setminus \{0\}$, $\C\times \C^*$ or  $\C^*\times \C^*$ then
  $\Gamma$  is  virtually affine. When $\Omega = \C^* \times D$ we show that $\Gamma$ leaves invariant the line $\{0\} \times \ell$ where $\ell$ is the projective line determined by $D$. Hence $\Gamma$ is affine.
  The most difficult case to deal with is when
   $\Omega =\C \times \mathbb H$. This is   a complex cone over  $\mathbb H$ with deleted top that we denote by $p$.  
  In this case we show that $p$ is $\Gamma$-invariant. This gives an action of  $\Gamma$ on the pencil $\mathcal P$ of projective lines passing through $p$. In particular,  the induced action of $Iso(\Omega,\Gamma)$ on the pencil $\mathcal P$ determines a subgroup $\Sigma$ of $\PSL(2,\R)$
  which preserves $\mathbb H$.  We prove that $\Sigma$ is  non-discrete  and it  has a fixed point $q$ in $\partial \mathbb H \subset \PC^1$. Using the classification of the Lie subgroups of $\PSL(2,\C)$ that we give below, we conclude that $q$ is a fixed point for  whole group $\Gamma$ acting on the pencil $\mathcal P$. This implies that $\Gamma$ leaves invariant a line in $\PC^2$ and therefore it is affine. 
  
    \vv
The proofs of Theorems 2-5 have all Theorem \ref{t:main1} as a key-ingredient. For instance, 
Theorem \ref{t:main2} uses Theorem \ref{t:main1}  together with  the quasi-minimality theorem in \cite{pablo1}   (Proposition \ref {p:pkg} below) and  an extension of the work of Kobayashi, Ochiai, Klingler and others for compact surfaces with a projective structure (Theorem \ref{t:klingler}).
 Theorem \ref{t:main3} is proved using Theorems \ref{t:main1} and \ref{t:main2} and Theorem \ref{t:klingler}, and so on. 
  Theorems 1 to 6 are proved in Section \ref{s: proofs}  using the content of sections 1 to 8.  

\vv
 We refer to \cite{CNS} for a  thorough discussion of the subject of
  complex Kleinian groups and 
  limit sets for discrete subgroups of $\PSL(n+1,\C)$. 

\vskip.4cm

\section{Preliminaries on non-discrete subgroups of
$\PSL(2,\C)$}\label{s:equicontinuo}  \label{s:dejedis}

In the sequel we must consider  subgroups of $\PSL(2,\C)$ which may not be discrete, so  
we begin this section with some definitions and  basic properties of these groups that we need. Here we only outline the main ideas and results;  the details are in an appendix at the end.

 Let $\SL(2,\C)$ be the group of $2 \times 2$ matrices with complex entries and determinant 1, and set
 $\PSL(2,\C) := \SL(2,\C)/\pm Id \,$. This is the  group of holomorphic automorphisms of  the projective line $\PC^1$. Recall that  $\PC^1$ is biholomorphic to the Riemann sphere $\s^2 \cong \widetilde \C := \C \cup \infty$, which can be regarded as being the sphere at infinity of the real hyperbolic space $\mathbb{H}^3$. These identifications, and some extra work, lead to the following well-known theorem:

\begin{theorem} There are canonical group isomorphisms:
$$\PSL(2,\C) \cong \hbox{M\"ob}(2,\C) \cong Iso_+ \mathbb{H}^3 \,,$$
where $\hbox{M\"ob}(2,\C) $ is the group of all M\"obius transformations acting on the Riemann sphere $\s^2 $ 
and $Iso_+ \mathbb{H}^3$  is the group of all orientation preserving isometries of  $\mathbb{H}^3$.
\end{theorem}

The group $\PSL(2,\C)$ contains the orthogonal group $\SO(3)$ as the stabilizer subgroup for  the action of $\PSL(2,\C)$ on $\mathbb{H}^3$.

We recall that if $\Gamma$ is a discrete subgroup of $\PSL(2,\C)$, then the limit set of $\Gamma$, usually denoted $\Lambda(\Gamma)$ or simply $\Lambda$, is the set of accumulation points of all orbits of points in $\mathbb{H}^3$.

Since the group is discrete and the action on  $\mathbb{H}^3$ is by isometries, one has that $\Lambda$  necessarily is contained in the sphere at infinity, which can be identified with $\s^2 $. We know (see for instance \cite{maskit}) that if 
$\Lambda$ has finite cardinality, then it consists of at most two points and the group is called {\it elementary}. Otherwise, that is for {\it non-elementary groups}, $\Lambda$ is a perfect set, and it is the set of accumulation points of all orbits.

The complement $\Omega := \s^2 \setminus \Lambda$ is the set of points where the action is discontinuous, and this is also the largest invariant set where the action is properly discontinuous. Moreover, $\Omega$ is also the equicontinuity set of $G$ (see Definition \ref{d:equi}).

In the sequel we need to make similar considerations for subgroups of $\PSL(2,\C)$ which may not be discrete. Notice that  if the group is not discrete,  there is not a discontinuity region.

\begin{definition} Identify  $\PSL(2,\C)$ with ${\rm Iso}_+ \mathbb{H}^3$ and think of $\s^2$ as the  sphere
at infinity of $\mathbb{H}^3$.  Let $\Gamma$ 
be a (discrete or not) subgroup of $\PSL(2,\C)$. Its {\it limit set in the sense of Greenberg}  (see
\cite {greenberg}), denoted $\Lambda_{\Gr}(\Gamma)$, is defined to
be the intersection  of $\s^2$ with the set of accumulation points of all
 orbits of points in $\mathbb{H}^3$.
\end{definition}

Of course this is the usual limit set for discrete groups. 

The theorem  below is to due Greenberg in
 \cite{greenberg}   (Theorem 1 and Proposition 12):

\begin{theorem} \label{t:greenberg}
Let $G$ be a connected Lie subgroup of  $\PSL(2,\C)$.
Then one of the following assertions is satisfied:
\begin{enumerate}
\item The elements  of $G$ have a common  fixed point  in
$\mathbb{H}^{3}$ and $G$ itself is conjugate to a Lie subgroup of
$\SO(3)$. 

\item The elements of $G$ have a common fixed point in
$\mathbb{P}^1_{\mathbb{C}} \cong \s^2$.

 \item There exists a hyperbolic 
line {\rm (a geodesic)}  $\ell\subset \mathbb{H}^{3}$ which  is $G$-invariant.

\item There exists a hyperbolic plane $L\subset \mathbb{H}^{3}$  which  is $G$-invariant. \item $G=\PSL(2,\C)$.
\end{enumerate}
Also, if 
$Card(\Lambda_{\Gr}(\Gamma))\geq  2$ then $\Lambda_{\Gr}(\Gamma)$
is the closure of the set of  fixed points of  loxodromic elements.
\end{theorem}

\begin{definition}\label{d:equi}
Given a subgroup $\Gamma \subset \PSL(2,\C)$,
its {\it equicontinuity region}, denoted $\Eq(\Gamma)$, is  
the set of   points $z\in \mathbb{P}^1_\mathbb{C}$ for which there
exists an open neighborhood $U$ of  $z$   such that $\Gamma\vert_U$  is
a normal family.
\end{definition}

\begin{remark} \label{r:in}
Notice that this set has the following properties:
\begin{enumerate}
\item If
$\gamma\in \Gamma$ is not elliptic, then
$Fix(\gamma)\subset \mathbb{P}^1_{\mathbb{C}}\setminus \Eq(\Gamma)$ .
\item
$\Eq(\Gamma)=\Eq(\overline{\Gamma})$. 
\item $\Eq(\Gamma)$ is an  open
$\Gamma$-invariant set.
\end{enumerate}
\end{remark}

The following are  examples of non-discrete subgroups of $\PSL(2,\C)$ that will appear in the sequel:

\begin{example} \label{e:dih}
The {\it  infinite dihedral group} $Dih_\infty$. This consists of the group $Rot_\infty \cong \SO(2)$ generated by all rotations around the origin,  to which we add the generator given by  the involution $h(z)=-z$. It satisfies $\Eq(Rot_{\infty})=\Eq(Dih_\infty)=\mathbb{P}^1_{\mathbb{C}} \setminus
\Lambda_{\Gr}(\Gamma)=\s^2$.
\end{example}

\begin{example}
The {\it special orthogonal group} $\SO(3)$ can be embedded in $\PSL(2,\C)$ as follows:
\[
\SO(3)=
\left \{\frac{az-\bar{c}}{cz+\bar{a}}\in M\ddot{o}b(\hat{\mathbb{C}})
:\vert a \vert^2 +\vert  c\vert^2 =1
\right \}.
\]
This is a purely elliptic group diffeomorphic to $\mathbb{P}^3_{\mathbb{R}}$, which satisfies $\Eq(\SO(3))=\s^2\setminus \Lambda_{\Gr}(\Gamma)=\Bbb{P}^1_{\C}$.
\end{example}

\begin{example} \label{e:punto} The group ${\rm Epa}(\mathbb{C})$ of all  affine M\"obius transformations which are either parabolic or elliptic. It satisfies $\Eq({\rm Epa}(\mathbb{C}))=\s^2 \setminus
\Lambda_{\Gr}(\Gamma)=\mathbb{C}$.
\end{example}

\begin{example} \label{e:2puntos}  The group $M\ddot{o}b(\mathbb{C}^*)$ of all   M\"obius transformations that leave invariant $\mathbb{C}^*$; it
 satisfies  $\Eq(M\ddot{o}b(\mathbb{C}^*))=\s^2\setminus 
\Lambda_{\Gr}(\Gamma)=\mathbb{C}^*$.
\end{example}

\begin{example}  \label{e:cir}
The group $M\ddot{o}b(\overline{\mathbb{R}})$ of all M\"obius transformations  that  leave  invariant  the circle  $\overline{\mathbb{R}}:= \mathbb{R}\cup\{\infty\}$; it satisfies  $\Eq(M\ddot{o}b(\overline{\mathbb{R}}))=\mathbb{C}\setminus \mathbb{R}=\s^2\setminus \Lambda_{\Gr}(\Gamma)$.
\end{example}

Now let us  give another description of  $\SO(3)$. For this, given an integer $p<0$,  let $\tau_p$ be the M\"obius transformation defined by  $\tau_p(z)=\frac{z-p}{z-1} $. Notice that $\tau_p(0) = p$, $\tau_p(\infty) = 1$ and $\tau_p$ is an involution, {\it i.e.},  $\tau_p= \tau_p^{-1}$.

\begin{definition}\label{e:cr}
{\it The Chinese rings group} $ Cr(p)$ is:
$$ Cr(p)=\overline{\langle Rot_\infty \,, \,
\tau_p \,Rot_\infty \, \tau_p\rangle} \;.$$
\end{definition}

\vv

We have (see the appendix for details):
 
 \begin{proposition} \label{p:cr1}     
\begin{enumerate}
\item \label{p:pgro}  The group $Cr(-1)$ is isomorphic to $\SO(3)$. In fact for each $p<0$,   $Cr(p)$ is
conjugate to $\SO(3)$.

\item \label{p:min1} The $Cr(-1)$-orbit of each  each $z\in \s^2$ is the whole Riemann sphere.

\item \label{p:esim1} For each $p<0$ there exist 
 $\gamma_p\in Cr(p)$ such that $\gamma_p$ has infinite order. Moreover, the fixed points of   $\gamma_p$, $z_1, z_2 \in \C$, are additive inverses, {\it i.e.}, 
$z_1 + z_2 = 0$.

\end{enumerate}
\end{proposition}

\begin{definition}
We say that a subgroup $\Gamma\subset \PSL(2,\Bbb{C})$ is {\it elementary} if its
equicontinuity set omits  at most 2
points in $\Bbb{P}^1_\Bbb{C}$.
 \end{definition}

The following theorems summarize the basic properties  we need about the equicontinuity region for non-discrete groups (we refer to the appendix for the proofs).

\begin{theorem} \label{c:5} \label{p:vacio1} Let $\Gamma\subset \PSL(2,\C)$
be a subgroup, then:
\begin{enumerate}
\item \label{i:el1} $\Eq(\Gamma)=\s^2$ if and only
if $\Gamma$ is  either  finite or   conjugate to a subgroup of   $Cr(-1)$ or $Dih_\infty$.
\item \label{i:el2}
 $\Eq(\Gamma)$ is $\mathbb{C}$, up to a projective transformation,   if
and only if $\Gamma$ is conjugate to a subgroup $ \Gamma_*$ of
${\rm Epa}(\mathbb{C})$ such that $\overline{\Gamma_*}$ contains a
parabolic element.
\item \label{i:el3} $\Eq(\Gamma) $ is  $\mathbb{C}^*$, up to a
projective transformation, if and only if $\Gamma$ is
conjugate to a subgroup $ \Gamma_*$ of $M\ddot{o}b(\mathbb{C}^*)$
such that
 $\Gamma_*$ contains a loxodromic element.
\end{enumerate}
\end{theorem}

\begin{definition}
Define the set of exceptional points of $\Gamma$ as $$Ex(\Gamma)=\{z\in
\s^2\setminus \Eq(\Gamma):\overline{\Gamma z}\neq
\s^2\setminus \Eq(\Gamma)\}.$$ 

\end{definition}

\begin{theorem}\label{t: eq2}
Let $\Gamma\subset \PSL(2,\C)$ be  a non-discrete infinite closed group, then:

\begin{enumerate}
\item $\Eq(\Gamma)=\s^2\setminus \Lambda_{\Gr}(\Gamma)$.

\item $\Gamma$ is purely elliptic if and only if
$\Eq(\Gamma)=\s^2$.

\item If
${H}\subset \Gamma$ is an infinite normal subgroup such that $Card(
\Lambda_{\Gr}({H}))= 2,0 $ and $H$ is not conjugate to a subgroup
of $\SO(3)$, then $\Gamma$ is elementary.

\item If $\mathcal{C}\neq Ex(\Gamma)$ is a closed $\Gamma$-invariant set, then $\Lambda_{\Gr}(\Gamma)\subset \mathcal{C}$.
\end{enumerate}
\end{theorem}

\begin{theorem}\label{t:eq3} 
Let $\Gamma\subset \PSL(2,\C)$ be  a non-discrete  closed group which is non-elementary, then:

\begin{enumerate}
\item If   $\Eq(\Gamma) \ne \emptyset$, then  $\Lambda_{\Gr}(\Gamma)$
is a circle in $\mathbb{P}^1_\mathbb{C}$.

\item 
$\Gamma$ contains  loxodromic elements  and 
 $\s^2\setminus \Eq(\Gamma)$ is the closure of the loxodromic fixed points.

\item 
One has $Card( Ex(\Gamma))<2$.
\end{enumerate}
\end{theorem}

Finally, the following theorem says that these groups also have the usual convergence properties, which are very useful.  

\begin{theorem} \label{c:16}  \label{trampa}
\begin{enumerate}
\item  If  $(\gamma_m)\subset \PSL(2, \mathbb{R})$ is a sequence, then  there exists a subsequence of $(\gamma_m)$, still denoted $(\gamma_m)$, such that either $(\gamma_m)$ converges uniformly to a M\"obius transformation $\gamma\in \PSL(2, \mathbb{R})$ or there exists a point $z\in \overline{\mathbb{R}}$ such that $(\gamma_m)$ converges uniformly to the constant function $z$ on compact sets of $\mathbb{H}^2$.
\item If  $\Gamma\subset \PSL(2,\C)$ is a subgroup and we have a sequence
$(\gamma_n)_{n\in\mathbb{N}}\subset
\Gamma$ such that  $\gamma_n \xymatrix{ \ar[r]_{n \rightarrow  \infty}&}
g$ uniformly on compact sets of $\Eq(\Gamma)$, where $g:\Eq(\Gamma)\rightarrow \s^2$,  then either $g$ is  a
constant function $z$ with  $z \in \Lambda_{\Gr}(\Gamma)$ or else $g\in
\PSL(2,\C)$ with $\gamma_n \xymatrix{ \ar[r]_{n \rightarrow
\infty}&} g$ uniformly on $\s^2$.

\end{enumerate}
\end{theorem}

\begin{proof} We only prove the first statement, since this essentially implies the second. 
Assume that the map $g$ is not a M\"obius transformation. Then the convergence property for isometries or real hyperbolic spaces implies that there exist points $x,y\in \Bbb{C}$ (which may coincide) and a subsequence of $(\gamma_n)$, that we still denote $(\gamma_n)$ for simplicity, such that  $\gamma_n \xymatrix{ \ar[r]_{m \rightarrow  \infty}&} y $ uniformly on compact sets of $\C$. Hence, just as in \cite{CS}, we have  $g=y$ and  $y\in \Bbb{C}\setminus Eq(\Gamma)$. 
\end{proof}

\vskip.4cm


\section{Complex Kleinian groups: definitions and some basic properties} \label{s:recall}
We recall that the complex projective plane
$\mathbb{P}^2_{\mathbb{C}}$ is $$\Bbb{P}^2_\Bbb{C}:=(\mathbb{C}^{3}\setminus \{0\})/\mathbb{C}^*
,$$
 where $\mathbb{C}^*$ acts on $\mathbb{C}^3\setminus\{0\}$ by the usual scalar
 multiplication. This  is   a  compact connected  complex $2$-dimensional
 Riemannian  manifold, naturally equipped with the Fubini-Study  metric.

Let $[\mbox{ }]:\mathbb{C}^{3}\setminus\{0\}\rightarrow
\mathbb{P}^{2}_{\mathbb{C}}$ be   the quotient map. If
$\beta=\{e_1,e_2,e_3\}$ is the standard basis of $\mathbb{C}^3$, for simplicity we
will write $[e_j]=e_j$ and if $w=(w_1,w_2,w_3)\in
\mathbb{C}^3\setminus\{0\}$ then we   write   $[w]=[w_1:w_2:w_3]$.
Also, a set  $\ell\subset \mathbb{P}^2_{\mathbb{C}}$ is said to be a
complex line if $[\ell]^{-1}\cup \{0\}$ is a complex linear
subspace of dimension $2$. Given  $p,q\in
\mathbb{P}^2_{\mathbb{C}}$ distinct points,     there exists a unique
complex line passing through  $p$ and $q$; such line will be
denoted by $\overleftrightarrow{p,q}$. If $\ell_1,\,
\ell_2$ are different complex lines then 
$\ell_1\cap \ell_2$ contains exactly one point.

 Taking   $\mathbb{H}^2_{\mathbb{C}}=\{[a:b:c]\in \mathbb{P}^2_{\mathbb{C}}:
 \vert a\vert+\vert b\vert <\vert c \vert \}$, it is not hard to show that
 $\mathbb{H}^2_{\mathbb{C}}$ is biholomorphic to the unitary ball in
 $\mathbb{C}^2$,  and its boundary $\partial{\mathbb{H}^2_{\mathbb{C}}} $ is diffeomorphic
 to the $3$-sphere. Notice that  for each point $p$  in $\partial \mathbb{H}^2_{\mathbb{C}}$
 there exists exactly one complex line tangent to
 $\partial \mathbb{H}^2_{\mathbb{C}}$ passing through $p$.

Consider the action of $\mathbb{Z}_{3}$ (viewed as the cubic roots of
the unity) on  $\SL(3,\mathbb{C})$ given by the usual scalar
multiplication, then
$\PSL(3,\mathbb{C}):= {\rm SL}(3,\mathbb{C})/\mathbb{Z}_{3}$ is a Lie group
whose elements are called projective transformations.  Let
$[[\mbox{  }]]:\SL(3,\mathbb{C})\rightarrow \PSL(3,\mathbb{C})$ be   the
quotient map and  $\gamma\in \PSL(3,\mathbb{C})$. We say that an element   $\widetilde\gamma\in
\SL(3,\mathbb{C})$  is a  lift of
$\gamma$ if   $[[\widetilde\gamma]]=\gamma$.  We also use the notation $(\gamma_{ij})$ to denote elements  in $\SL(3,\Bbb{C})$. 

One can show that
$ \PSL_{3}(\mathbb{C})$ is a Lie group  that acts  transitively,
effectively and by biholomorphisms  on $\mathbb{P}^2_{\mathbb{C}}$
by $[[\gamma]]([w])=[\gamma(w)]$, where $w\in
\mathbb{C}^3\setminus\{0\}$ and    $\gamma\in \SL(3,\mathbb{C})$. It is an exercise 
 to show that   projective transformations  take complex
lines into complex lines.

 \vv
  
We now consider    a group $G$ acting on a space $X$, $g\in G$ and $A\subset
X$ a subset. We define the {\it isotropy group} of $A$ as $Isot(A,G)=\{g\in G:g(A)=A\}$;   by $GA$ we denote  the orbit of $A$ under $G$ and by $Fix(g)$
the set of fixed points of $g$. We say  that the action of
$G$ is locally faithful  if whenever $f,g\in G$    agree in some
open set then one has $f=g$ on $X$.

\begin{definition} \label{d:lim}
 Let $\Gamma\subset   \PSL(3,\mathbb{C})$ be a discrete subgroup. We  define
 (following Kulkarni  \cite{kulkarni}):

\begin{enumerate}
\item $L_0(\Gamma)$  is the closure  of  the points in
$\mathbb{P}^2_{\mathbb{C}}$ with infinite isotropy group. 
\item
$L_1(\Gamma)$ is the closure of the set  of cluster (or accumulation) points  of
$\Gamma z$  where  $z$ runs  over  $\mathbb{P}^2_{\mathbb{C}}\setminus
L_0(\Gamma)$. Recall that $q$ is a cluster point  for  $\Gamma K$,
where $K\subset \mathbb{P}^2_{\mathbb{C}}$ is a non-empty set, if there exist a sequence
$(k_m)_{m\in\mathbb{N}}\subset K$ and a sequence of distinct elements
$(\gamma_m)_{m\in\mathbb{N}}\subset \Gamma$ such that
$\gamma_m(k_m)\xymatrix{ \ar[r]_{m \rightarrow  \infty}&} q$.

\item  $L_2(\Gamma)$ is  the closure of cluster  points of $\Gamma
K$  where $K$ runs  over all  the compact sets in
$\mathbb{P}^2_{\mathbb{C}}\setminus (L_0(\Gamma) \cup L_1(\Gamma))$.

\item The  \textit{Kulkarni limit set}  of   $\Gamma$  is:  $$\Lambda_{\Kul} (\Gamma) = L_0(\Gamma) \cup
L_1(\Gamma) \cup L_2(\Gamma).$$ 

\item The \textit{Kulkarni 
region of discontinuity} of $\Gamma$ is:
$$\Omega_{\Kul}(\Gamma) = \mathbb{P}^2_{\mathbb{C}}\setminus
\Lambda_{\Kul}(\Gamma).$$
\end{enumerate}
 \end{definition}

Recall that if $G$ is a Lie group acting on a smooth manifold $N$, we say that 
the action is {\it properly discontinuous} on a $G$-invariant set $U$  if for every compact set $K $ in $U$ one has that there exist only finitely many elements $g \in G$ such that $g(K) \cap K \ne \emptyset$.
We know from \cite{kulkarni} that the action of $\Gamma$ on $\Omega_{\Kul}(\Gamma)$ is properly discontinuous. We also have:

\begin{lemma}\label{i:pk5} Let    $\Gamma \subset \PSL(3,\C)$  be a complex  Kleinian group. Then
 the equicontinuity region   of $\Gamma$ is contained in the Kulkarni region of discontinuity: $Eq(\Gamma) \subset
\Omega_{Kul}(\Gamma)$. 
\end{lemma}
\begin{proof}
By Corollary 3.4 in \cite{CS}, it follows that $\Gamma$ acts properly discontinuously on $Eq(\gamma)$. Moreover, the same result yields that for every compact  $K\subset Eq(\Gamma)$ the clusters points of $\Gamma K$ lie on $(L_0(\Gamma)\cup L_1(\Gamma))\cap P^2\setminus Eq(\Gamma)$. Then by Proposition \ref{p:pkg} we conclude that $Eq(\Gamma)\subset \Omega_{Kul}(\Gamma)$.
\end{proof}

\vv
Following  \cite{sv1, CNS} we have:

\begin{definition} 
A discrete group $\Gamma \subset \PSL(3,\C)$ is a \textit{complex  Kleinian group}
 if there exists a non-empty  open invariant set $\Omega$ where the action is properly discontinuous.
 \end{definition}

\begin{proposition}{\rm (See \cite{pablo1})} \label{p:pkg} 
Let    $\Gamma \subset \PSL(3,\C)$  be a complex  Kleinian group . Then:

\begin{enumerate}
\item \label{i:pk1} $\Gamma$ is  countable. 
\item \label{i:pk2}
$\Lambda_{\Kul}(\Gamma),\,L_0(\Gamma),\,L_1(\Gamma),\,L_2(\Gamma)$
are  invariant sets. 

\item \label{i:pk4} If  $\mathcal{C}\subset\mathbb{P}^2_{\mathbb{C}}$ is 
a closed $\Gamma$-invariant set such that  for every compact set $K\subset
\mathbb{P}^2_{\mathbb{C}}- \mathcal{C}$, the set of cluster points
of  $\Gamma K$ is contained in $(L_0(\Gamma)\cup L_1(\Gamma))\cap \mathcal{C}$, then
$\Lambda_{\Kul}(\Gamma)\subset \mathcal{C}$.

\end{enumerate}
\end{proposition}

 We recall  that the  Chen-Greenberg limit set of $\Gamma \subset \PU(2,1)$, denoted $\Lambda_{CG}(\Gamma)$, is  the set of cluster points of an orbit
$\Gamma z$, where $z$ is any point in $\mathbb{H}^2_{\mathbb{C}}$. One has (see  \cite{cg}):

\begin{theorem}[Chen-Greenberg] Let $\Gamma \subset \PU(2,1)$ be a discrete group, then:
\begin{enumerate}
\item $\Lambda_{CG}(\Gamma)$ does not depend on the choice of $z$.
\item The cardinality  $Card( \Lambda_{CG}(\Gamma))$ 
   is  $\,0,\,1,\,2,\,$ or $\infty$.
\item If $Card (\Lambda_{CG}(\Gamma))=\infty$, then $\overline
{\Gamma z}= \Lambda_{CG}(\Gamma)$ for every $z\in
\Lambda_{CG}(\Gamma)$.
\end{enumerate}
\end{theorem}


There are several types of complex Kleinian groups which are specially relevant for this article and we introduce them now.

\begin{definition}\label{definitions}
Let $\Gamma$ be a  complex Kleinian subgroup of $\PSL(3,\C)$.  Then:

\begin{enumerate}
\item $\Gamma$ is {\it complex hyperbolic} if it is conjugate in $\PSL(3,\C)$ to a subgroup of $\PU(2,1)$.

\item  $\Gamma$ is {\it affine} if it is conjugate in $\PSL(3,\C)$ to a subgroup that leaves invariant an affine chart, {\it i.e.}, a copy of $\C^2$  embedded in $\PC^2$. The group is {\it  virtually  affine} if it has a finite index subgroup which is affine.

\item  An affine group  $\Gamma$ is called {\it  controllable}  if there exist a line $\ell$ and a point
$p\notin\ell $ which are invariant under the action of $\Gamma$. In this case 
the group   $\Gamma\vert_{\ell}=G$ is called the control group. This is a subgroup of $\PSL(2,\C)$ which can be non-discrete, even though  $\Gamma$ is discrete.  {\rm (See Lemma (\ref{l:lic}) and \cite[Chapter 5]{CNS}.)}

\item $\Gamma$ is {\it  elementary} if its Kulkarni limit set 
$\Lambda_{\Kul}(\Gamma)$ consists of finitely many projective lines and finitely many points away from these lines. {\rm (See \cite[Chapter 6]{CNS}.)}

\item $\Gamma$ is {\it  quasi-cocompact} if there exists an open $\Gamma$-invariant subset of $\PC^2$ where the action is properly discontinuous and the quotient $\Omega/\Gamma$ is compact. In this case we say that $\Gamma$ {\it acts cocompactly} on $\Omega$.

\end{enumerate}
\end{definition}

 In the specially interesting case of complex hyperbolic groups, one has the following Theorem of \cite{pablo1}.

\begin{theorem} [J. P. Navarrete] \label{t:pfuch}
Let  $\Gamma\subset \PU(2,1)$ be a discrete group.  One has:
\begin{enumerate}
\item $\Lambda_{\Kul}(\Gamma)=\bigcup_{p\in \Lambda_{CG}(\Gamma)}
T_p$, where $T_p$ is the unique complex projective line tangent to
$\partial \mathbb{H}^2_{\mathbb{C}} \cong \s^3$ at $p$.
\item If  $Card (\Lambda_{CG}(\Gamma))=\infty$, then
$\Omega_{\Kul}(\Gamma) $ is the largest open set on which  $\Gamma$
acts    properly discontinuously, and  $\Omega_{\Kul}(\Omega)$ coincides with the equicontinuity set of $\Gamma$.
\end{enumerate}
\end{theorem}

To close this section we provide the following description on the limit set of cyclic groups.

\begin{theorem} [J. P. Navarrete, see \cite{pablo2}] \label{t:pciclic}
Let $\gamma\in \PSL(3,\mathbb{C})$ and $\widetilde\gamma\in
\SL(3,\mathbb{C})$ be a lift of $\gamma$. The limit set in the sense
of Kulkarni for the cyclic group generated by $\gamma$ {\rm (denoted
$\langle\gamma\rangle$)}, in terms of the Jordan's normal form of $\widetilde
\gamma$,  is given by:
\medskip

\begin{tiny}
\begin{center}
\begin{tabular}{|c|l|c|c|c|}
\hline $
\begin{array}{c}
\\
\textrm{Normal Form}\\
\textrm{of } \widetilde\gamma\\
\\
\end{array}
$ & $
\begin{array}{c}
\textrm{Condition over the}\\
\textrm{eigenvalues} \end{array} $ &  $L_0(\langle\gamma\rangle)$&
$L_1(\langle\gamma\rangle)$ & $L_2(\langle\gamma\rangle)$ 
\\
\hline $
\begin{array}{l}
 \\
\left (
\begin{array}{lll}
1  & 1 & 0\\
0  & 1 & 1\\
0  & 0 &1\\
\end{array}
\right)\\\\
\end{array}
$ & 1 \textrm{ is the only eigenvalue} &$\{e_1\}$ & $\{e_1\}$ & $\overleftrightarrow {e_1,e_2}$ 
\\
\hline 
\multirow{5}{*}[-8mm]{ 
$ 
\left (
\begin{array}{lll}
\lambda_1  & 0            & 0\\
0            &  \lambda_2 & 0\\
0            &  0            & \lambda_3\\
\end{array}
\right) $ } & $
\begin{array}{l}
\\
\lambda_1^n=\lambda_2^n=\lambda_3^n=1\\
\textrm{ for some } n
\end{array}
$ 
& $ \emptyset  $ & $\emptyset$  & $\emptyset$
\\
\cline{2-5} & $
\begin{array}{l}
 \\
\vert \lambda_1\vert =\vert \lambda_2\vert =1 \textrm{ and}\\
 \lambda_3^n\neq 1 \textrm{ for all }  n.\\
\\
\end{array}
$  & unimportant & $ \mathbb{P}^2_{\mathbb{C}} $ & $ \emptyset$
\\
\cline{2-5} & $
\begin{array}{l}
\\
\vert \lambda_3\vert\neq \vert\lambda_1\vert=\vert\lambda_2\vert,\\
\vert\lambda_2 \vert\neq 1 \textrm{ and }\\
\lambda_2^n =\lambda_1^n \textrm{ for some } n\in \mathbb{N} \\
\\
\end{array}$ 
& $ \overleftrightarrow{e_1,e_2}\cup\{e_3\} $ 
& $ \overleftrightarrow{e_1,e_2}\cup\{e_3\} $ 
& $ \overleftrightarrow{e_1,e_2}\cup\{e_3\} $
\\
\cline{2-5} & $
\begin{array}{l}
\\
\vert \lambda_3\vert\neq \vert\lambda_1\vert=\vert\lambda_2\vert,\\
\vert\lambda_2 \vert\neq 1 \textrm{ and }\\
\lambda_2^n \neq \lambda_1^n \textrm{ for all } n\in \mathbb{N} \\
\\
\end{array}$ 
& $\{e_1,e_2,e_3\} $ 
& $ \overleftrightarrow{e_1,e_2}\cup\{e_3\} $ 
& $ \overleftrightarrow{e_1,e_2}\cup\{e_3\} $
\\
\cline{2-5} & $
\begin{array}{l}
\\ \\
\vert \lambda_1\vert< \vert\lambda_2\vert<\vert\lambda_3\vert\\
\\ \\
\end{array}
$ 
&
$\{e_1,e_2,e_3\} $ 
& 
$\{e_1,e_2,e_3\} $
& $
\overleftrightarrow{e_1,e_2}\cup \overleftrightarrow{e_3,e_2}
$
\\
\hline
\multirow{3}{*}[-3mm]{ $ \left (
\begin{array}{lll}
\lambda_1    & 1            & 0\\
0            &  \lambda_1 & 0\\
0            &  0              & \lambda_1^ {-2}\\
\end{array}
\right) $ } & $
\begin{array}{l}\\
 \lambda_1^n= 1 \textrm{ for some } n\in \mathbb{N}\\
\\
\end{array}
$
&
$\overleftrightarrow{e_1,e_3} $
&
$
\{e_1\}
$
&
$
\{e_1\}
$
\\
\cline{2-5}
&$
\begin{array}{l}\\
\vert \lambda_1\vert = 1 \textrm{ and } \lambda_1^n\neq 1\\
\textrm{ for all } n \in\mathbb{N}
\\\\
\end{array}
$
&
$\{e_1,e_3\} $
&
$
\overleftrightarrow{e_1,e_3}
$
&
$
\{e_1\}
$
\\
\cline{2-5} & $
\begin{array}{l}\\ \\
\vert \lambda_1\vert\neq 1\\ \\ \\
\end{array}
$ & $\{e_1,e_3\}$
&
$\{e_1,e_3\} $
&
$
\overleftrightarrow{e_1,e_2}\cup \overleftrightarrow{e_3,e_1}
$
\\
\hline
\end{tabular}

\end{center}
\end{tiny}
\end{theorem}



\vskip.4cm

\section{Compact  complex projective 2-dimensional orbifolds} \label{s:corb}
 Recall that a  {\it complex projective structure} on a complex manifold $M$ of dimension $n$ is a maximal atlas for $M$
  modeled on open subsets of the projective
space $\mathbb {P}^n_{\mathbb {C}}$, so that for any two overlapping
charts, the corresponding change of coordinates is restriction of an element in $\PSL(n+1,\C)$.

More generally, we want to consider complex projective  orbifolds, and we refer to \cite{choi, thurston} for basic material 
on this topic.  Like a manifold, an orbifold is specified
by local conditions, an {\it orbifold atlas};  instead of being locally modeled on
open subsets of $\R^n$, an orbifold is locally modeled on quotients
of open subsets of $\R^n$ by finite group actions.  More precisely, 
if $G$ is a Lie group acting effectively, transitively and locally
 faithfully on a smooth manifold $X$, a {\it  $(G,X)$-orbifold} means a  topological Hausdorff
space $M$, called the underlying space, with a
countable basis and equipped with a collection $\{\widetilde{U}_i,\Gamma_i,
\phi_i,U_i\}_{i\in I}$, where the $\{U_i\}$  are an open cover of  $M$, and 
the $\widetilde {U}_i$, called {\it folding charts}, are 
 open subsets of $ X$.
 For each $\widetilde {U}_i$ there exists  a finite group $\Gamma_i\subset G$ acting
on $\widetilde {U}_i$, and a homeomorphism $\phi_i:U_i\rightarrow
\widetilde{U}_i/\Gamma_i$, called a {\it folding map}. These charts
must satisfy a certain compatibility condition (as they do for
manifolds).

The orbifold $M$  is said to be  {\it good} if  there exists a covering orbifold  map
$p:\widetilde M\rightarrow M$ such that $\widetilde M$ is a  manifold. The orbifold is {\it very good} if the manifold $\widetilde M$ actually is compact.  The following theorem is due to Thurston in  \cite {thurston}:  

\begin{theorem}\label{t:thurston} If
$G$ is an analytic group of diffeomorphisms of a manifold $X$, then every $(X, G)$-orbifold is good. Moreover, there exists a simply connected manifold $\widetilde M$ and an $(X, G)$-covering orbifold map which is  unique up to equivalence.
Then {\rm (by globalizing the coordinate charts and a little more work)} one gets a {\it developing  map}:
$$\mathcal D: \widetilde M \to \PC^2 \,,$$
and  a {\it holonomy homomorphism} $\mathcal H: \pi_1^{Orb}(M) \to G$, where $\pi_1^{Orb}(M)$ is the orbifold  fundamental group of $M$. 
\end{theorem}

The simply connected manifold   $\widetilde M $ is  called {\it the universal orbifold cover} of $M$,  and 
the pair $(\mathcal D, \mathcal H) $ is called a {\it developing pair} for the $(X, G)$-orbifold $M$. 
\medskip

In  the case we envisage here  we have  $X= \mathbb{P}^2_{\mathbb{C}}$ and $G =  \PSL(3,\mathbb{C})$. 
One gets  a {\it holonomy homomorphism} $\mathcal H: \pi_1^{Orb}(M) \to \PSL(3,\mathbb{C})$ and one has:

\begin{lemma}\label{c:vego}
Let $M$ be a compact
$(\mathbb{P}^2_{\mathbb{C}}, \PSL(3,\mathbb{C}))$-orbifold, then $M$ is a
very good orbifold.
\end{lemma}

\begin{proof}
 Since  $M$ is compact,  there exists a metric $d$ on the universal covering  orbifold $\widetilde M$ of $M$ which is compatible with its topology. Thus   $(\widetilde M, d )$ is a length space, which is    connected  and geodesically
 complete, see \cite{borzellino, rat}. Thus  $\pi_1^{Orb}(M)$    is  a finitely generated  subgroup 
 of isometries (see \cite{rat}). If
$(\mathcal{D},\mathcal{H})$ is a developing pair for $M$,   by
Selberg's Lemma (see \cite{rat}), $\mathcal{H}(\pi_1^{Orb}(M))$
has a   finite index normal subgroup   $G$ which is torsion free.
Therefore   $\widetilde G=\mathcal{H} ^{-1}(G)$ is a  normal  subgroup of
$\pi_1^{Orb}(M)$ with finite index.

We claim that  the group $\widetilde
G$ acts freely on $\widetilde M$, for otherwise,   if  $x\in M$ is  such that
$Isot(x, \widetilde G)$ is non-trivial, then $Isot(x,\widetilde G)\subset
Ker(\mathcal{H})$ and there exists an open neighborhood $W$ of $x$
which is $Isot(x,\widetilde G)$-invariant and such that
$\mathcal{D}\vert_W$  is  injective. Hence
$\mathcal{D}(g(z))=\mathcal{D}(z)$ with  $g(z)\in W$ for all $z\in
W$ and any $g\in Isot(x, \widetilde G)$. Thence $Isot(x,\widetilde
G)=\{Id\}$ which is a contradiction.

It follows     that  $N=\widetilde
M/\widetilde G$ is a
$(\mathbb{P}^2_{\mathbb{C}},\PSL(3,\mathbb{C}))$-manifold and
$\Gamma=\pi_1^{Orb}(M)/\widetilde G$ is a finite group that acts on
$N$ by  $\widetilde G g (\widetilde G x) =\widetilde G gx$ and  $\Gamma\subset Aut_{(\mathbb{P}^2,\PSL(3,\mathbb{C}))}(N)$. Let  $\varphi:
\widetilde M/\pi_1^{Orb}(M)\longrightarrow N/\Gamma$ be defined by
$\varphi( \pi_1^{Orb}(M)x)= \Gamma(\widetilde G x)$, clearly    $\varphi$
is an $(\mathbb{P}^2_{\mathbb{C}},
\PSL(3,\mathbb{C}))$-equivalence. Since  $Card( \Gamma)
<\infty$ is finite  and $N/\Gamma$ is  compact, it follows  that    $N$
is compact.
\end{proof}

Now  recall from  the introduction that 
 $\Aff(\mathbb{C}^2)$ is the affine group regarded as  $\{g\in \PSL(3,\mathbb{C}):  g(\mathbb{C}^2)=\mathbb{C}^2\}$. In the introduction we defined also the groups $A_1$, $A_2$, $Sol^4_0$, $Sol^4_1$ and $Sol^{\prime 4}_1$ that apear in the next theorem.
Using  Lemma \ref{c:vego}, the results of 
Kobayashi, Ochiai and Inoue \cite{KO, KOI}, Klingler \cite{klingler, klingler2}, Mok and Yeung \cite {Mok1, Mok2, MY} \textit{et al} for
compact $(\mathbb{P}^2_{\mathbb{C}},\PSL(3,\mathbb{C}))$-manifolds   extend easily to compact
$(\mathbb{P}^2_{\mathbb{C}},\PSL(3,\mathbb{C}))$-orbifolds as
follows.

\begin{theorem}  \label{t:klingler}
Let $M$ be a compact
$(\mathbb{P}^2_{\mathbb{C}},\PSL(3,\mathbb{C}))$-orbifold and denote by $\Sigma_M$ its singular set.
Then $M$
is of one of the following 8 types.
\begin{enumerate}

\item \label{i:k1} $\widetilde M=\mathbb{P}^2_{\mathbb{C}}$. Then $\mathcal{D}=Id $ and $\mathcal{H}(\pi^{Orb}_1(M))$ is
finite.

\item \label{i:k2} $\widetilde M=\mathbb{H}^2_{\mathbb{C}}$. Then $\mathcal{D}=Id$  and $\mathcal{H}(\pi_1^{Orb}(M))\subset
\PU(2,1)$.

\item  \label{i:k3}$\widetilde M=\mathbb{C}^2\setminus\{0\}$. Then 
$\mathcal{D}=Id,\, \mathcal{H}(\pi_1^{Orb}(M))\subset
\Aff(\mathbb{C}^2)$ contains a cyclic group  of finite index
generated by a contraction, and $M$ has a  possibly ramified  finite covering  which is  a primary Hopf surface. 

\item \label{i:k4} $\mathcal{D}(\widetilde
M)=\mathbb{C}^*\times\mathbb{C}^*$. Then
$\mathcal{H}(\pi_1^{Orb}(M))\subset Aut(\mathbb{C}^*\times\mathbb{C}^*) \subset  \Aff(\mathbb{C}^2)$ and $M$ has a possibly ramified
finite   covering which is  a surface biholomorphic
to a complex torus. 

\item \label{i:k5} $\mathcal{D}(\widetilde
M)=\mathbb{C}\times\mathbb{C}^*$. Then
$\mathcal{H}(\pi_1^{Orb}(M))\subset Aut(\mathbb{C}\times\mathbb{C}^*) \subset \Aff(\mathbb{C}^2)$ has a
subgroup of finite index contained in   $A_1$ or $A_2$, and $M$  has
a  possibly ramified  finite   covering  which is  a complex a surface
biholomorphic to a  torus. 
\item \label{i:k6} $\mathcal{D}(\widetilde M)=\mathbb{C}^2$. Then
  $\mathcal{H}(\pi_1^{Orb}(M))\subset \Aff(\mathbb{C}^2)$ contains
a unipotent  subgroup of  finite index   and  $M$ has a  possibly ramified  finite
covering  which is   a surface biholomorphic to a
complex torus or a primary Kodaira surface.

\item  \label{i:k7} 
$\widetilde M=\mathbb{H}\times \Bbb{C}=\mathcal{D}(\widetilde M)$.  Then   $\mathcal{H}(\pi_1^{Orb}(M))$ contains
a   subgroup of  finite index   which is contained in either $Sol^4_0, Sol^4_1$ or $Sol^{\prime 4}_1$ and 
$M$ has a possibly ramified   finite covering which is an Inoue  Surface.

\item \label{i:k8} $\widetilde M$ is biholomorphic (but not projectively equivalent)  to
$\mathbb{H}\times \Bbb{C}$. Then there exist
$A,B:\mathbb{H}\longrightarrow \mathbb{C}$   holomorphic maps and
$\mu\in \mathbb{C}^*$ such that  $\mathcal{D}(\widetilde
M)(w,z)=(A(z)e^{w\mu},B(z)e^{w\mu})$. If $\Sigma_M$ is empty, then
$\mathcal{H}(\pi_1^{Orb}(M))\subset \Aff(\mathbb{C}^2)$,
$\pi_1^{Orb}(M)\subset Bihol(\mathbb{C}\times \mathbb{H})$ contains a
subgroup of finite index $\Xi$ which admits the presentation:
$$\langle a_1,b_1,\ldots,a_g,b_g,c,d: c,d \textrm{ are in the center  and
}\Pi_{i=1}^g[a_i,b_i]=c^r\rangle\;,$$ where $2\leq g$, $r\in \mathbb{N}$.  In this case
$M$ has a possibly ramified  finite covering  which is an elliptic affine
surface.
\end{enumerate}
\end{theorem}
  
  We remark that it can happen that the singular set be empty, {\it i.e.}, $\Sigma_M=\emptyset$. In this case $M$ is a projective manifold and $\pi_1^{Orb}(M)$ is the usual fundamental group.

\vv We close this section with 
the following description of the  Inoue   surfaces given by C. T. C. Wall in \cite{wall}, that we use in the sequel.

 \begin{theorem} \label{t:sol40} 

 Let $M$ be an Inoue surface. Then $M$ admits a projective structure and one of the following assertions is satisfied:   

\begin{enumerate}
\item There 
exists $A\in \SL(3,\mathbb{Z})$ having eigenvalues
$\alpha,\beta,\overline{\beta}$ with $\alpha>1$, $\beta\neq
\overline{\beta}$,  a  real eigenvector  $(a_1,a_2,a_3)$
belonging to $\alpha$ and an eigenvector $(b_1,b_2,b_3)$ belonging
to $\beta$, such that  $\pi_1(M)$  is  generated by:
\[
\begin{array}{l}
\gamma_0(w,z)=(\alpha w,\beta z),\\
\gamma_i(w,z)=(w+a_i,z+b_i)\,, \, i=1,2,3.
\end{array}
\]
In this case $(\Bbb{H}\times \Bbb{C}))/\Gamma$ is a 3-torus bundle over a circle and $\Gamma$ belongs to $Sol^4_0$. 

\item There exists  $N\in \SL_2(\mathbb{Z})$
having real  eigenvalues $\alpha,\alpha^{-1}$, 
real eigenvectors   $(a_1,a_2)$,  $(b_1,b_2)$ for $\alpha,\alpha^{-1}$ respectively,  an integer $r$ and complex numbers  $t,c_1,c_2$, such that
$\pi_1(M)$ is generated by:
\[
\begin{array}{l}
\gamma_0(w,z)=(\alpha w,\beta z+t),\\
\gamma_i(w,z)=(w+a_i,z+b_iw+c_i)\, , \, i=1,2,\\
\gamma_3(w,z)=(w,z+r^{-1}(b_1a_2-b_2a_1)).
\end{array}
\]
In this case $(\Bbb{H}\times \Bbb{C}))/\Gamma$ is a bundle over a circle and the fiber itself is a circle bundle over a 2-torus and $\Gamma$ belongs to $Sol^{\prime 4}_1$.

\item There  exists
 $N\in GL(2,\mathbb{Z}) $
having  real eigenvalues $\alpha,-\alpha^{-1}$,
real eigenvectors   $(a_1,a_2)$,  $(b_1,b_2)$ for $\alpha,\alpha^{-1}$ respectively,  an integer $r$ and  complex numbers $c_1,c_2$,  such that
$\pi_1(M)$ is generated by: 
\[
\begin{array}{l}
\gamma_0(w,z)=(\alpha w,-z)\\
\gamma_i(w,z)=(w+a_i,z+b_iw+c_i)\,, \, i=1,2,\\
\gamma_3(w,z)=(w,z+r^{-1}(b_1a_2-b_2a_1)).
\end{array}
\]
In this case $(\Bbb{H}\times \Bbb{C}))/\Gamma$ is the quotient of the manifold above  by an unramified $\Bbb{Z}_2$-action and $\Gamma$ belongs to $Sol^4_1$.
\end{enumerate}
\end{theorem}



\vskip.4cm

\section{Projective Structures  and   Projective Groups} \label{s:puente}

In this section we prove:

\begin{theorem}\label{t:pi_1-orbifolds}
Let $\Gamma\subset \PSL(3,\mathbb{C})$ be a group acting properly  
discontinuously on a non-empty $\Gamma$-invariant   (connected) domain
$\Omega$ in $\PC^2$. Then:
\begin{enumerate}
\item $\Gamma$ can be realized as the image of the holonomy map for some appropriate 
 developing pair
$(\mathcal{D},\mathcal{H})$ for $M=\Omega/\Gamma$;  and

\item The orbifold fundamental group of $M$ satisfies $\Gamma=\pi_1^{Orb}(M)/\pi_1(\Omega)$. Moreover,
 $\pi_1^{Orb}(M) $ is an extension of $\Gamma$ by 
the fundamental group $\pi_1(\Omega)$.
 
\end{enumerate}
\end{theorem}

This theorem is a consequence of the two lemmas below:

\begin{lemma}\label{l:puente}
Let $\Gamma\subset \PSL(3,\mathbb{C})$ be a group acting properly  
discontinuously over a non empty, $\Gamma$-invariant   domain
$\Omega$. Then there exists a developing pair
$(\mathcal{D},\mathcal{H})$ for $M=\Omega/\Gamma$ such that
$\mathcal{D}(\widetilde M)=\Omega$ and
$\mathcal{H}(\pi_1^{Orb}(M))=\Gamma$.
\end{lemma}

\begin{proof} Step 1. -Construction of $\mathcal{D}$ - Let
$P:\widetilde M\longrightarrow M$ be the universal covering orbifold
map, $q:\Omega\longrightarrow M$   the quotient map,   $m\in
X_M\setminus\Sigma_M$, $\widetilde m\in \widetilde M$ and  $x\in \Omega$ such that
$P(\widetilde m)=q(x)=m$. By \ref{t:thurston} there exists a
$(\mathbb{P}^2_{\mathbb{C}},\PSL(3,\mathbb{C}))$-covering  map
$\hat{\mathcal{D}}:(\widetilde M,\widetilde m)\longrightarrow (\Omega,x)$
such that $q\circ \hat{\mathcal{D}} =P$. 
Let $i:\Omega\longrightarrow\mathbb{P}^2_{\mathbb{C}}$  be the
inclusion map  and define $\mathcal{D}=i\circ \hat{\mathcal{D}}$. This is a
$(\mathbb{P}^2_{\mathbb{C}},\PSL(3,\mathbb{C}))$-map.\vv

Step 2. -Construction of $\mathcal{H}$- Let $g\in \pi_1^{Orb}
(M)$. Since  $q(\mathcal{D}(g(\widetilde m)))=q(x)$,  there exists
$\hat{g}\in \Gamma$ such that $\hat{g}(x)=\mathcal{D}(g(\widetilde
m))$. Since $Isot(x,\Gamma)$ is trivial  we conclude  that
$\hat{g}$ is unique. Define $\mathcal{H}:\pi_1^
{Orb}(M)\longrightarrow \Gamma$ by $\mathcal{H}(g)=\hat g$.\vv

 Step 3. - For all $g\in \pi_1(M)$ we have  that  $\mathcal{D}\circ g= \mathcal{H}(g)\circ \mathcal{D}$- Let
$g\in \pi_1^ {Orb}(M)$, then by   \ref{t:thurston} there exists a
$(\mathbb{P}^2_{\mathbb{C}},\PSL(3,\mathbb{C}))$-map  $S:(\widetilde M,\widetilde
m)\longrightarrow(\widetilde M,\widetilde m)$ such that $\mathcal{D}\circ
g\circ S= \mathcal{H}(g)\circ \mathcal{D}$. Therefore  the following
diagram commutes and  from \ref{t:thurston} we obtain   $S=Id_{\widetilde M}$.
\[ 
\xymatrix{
(\widetilde M,\widetilde m) \ar[dr]^P \ar[rr]^S \ar[dd]^{\hat{\mathcal{D}}} & &(\widetilde M,\widetilde m) \ar[dl]_P \ar[dd]_g \\
                                     & (M,m)& \\
(\Omega,x ) \ar[dr]_{\mathcal{H}(g)}\ar[ur]_q &&(\widetilde M,g(\widetilde m))\ar[ul]^P \ar[dl]^{\hat{\mathcal{D}}}\\
                                     & (\Omega,\hat{\mathcal{D}}(g(x))).\ar[uu]_q &
}
\]

Step 4. -$\mathcal{H}$ is a group morphism- Let $g,h\in
\pi_1^{Orb}(M)$, then
$\mathcal{H}(g)\circ\mathcal{H}(h)(x)=\mathcal{H}(g)(\mathcal{D}(h(\widetilde
m)))=\mathcal{D}(g(h(\widetilde m)))$. That is  $\mathcal{H}(g\circ
h)=\mathcal{H}(g)\circ\mathcal{H}(h)$.\vv

Step 5. -$\mathcal{D}(\pi_1^{Orb}(M))=\Gamma$-  Let $g\in \Gamma$
and $z\in\widetilde M $ such that $\mathcal{D}(z)=g(x)$. By  Theorem
\ref{t:thurston}  there exists a
$(\mathbb{P}^2_{\mathbb{C}},\PSL(3,\mathbb{C}))$-equivalence
$\hat{g}:(\widetilde M,\widetilde m) \longrightarrow (\widetilde M,z)$ such
that  $\mathcal{D}\circ \hat g=g\circ \mathcal{D}$. Therefore    the
following diagram commutes and we have $\hat{g}\in \pi_1^{Orb}(M)$.
\[
\xymatrix{
(\widetilde M,\widetilde m) \ar[dr]^P \ar[rr]^{\hat{g}} \ar[dd]^{\hat{\mathcal{D}}} & & (\widetilde M, z) \ar[dl]_P \ar[dd]_{\hat{\mathcal{D}}} \\
                                     & (M,m)& \\
(\Omega,x ) \ar[ur]_q \ar[rr]^g & & ( \Omega,g(x))\ar[ul]^q.\\
}
\]
\end{proof}

\begin{lemma} \label{c:kerhol}
Let $\Gamma,\,\Omega,\, \widetilde m,\,  x,\, m,\,P,\,q,\,\mathcal{D}$
and $\mathcal{H}$ be  as in Lemma \ref{l:puente}, then we have the
following exact sequence of groups:
\[
\xymatrix{
0 \ar[r] & \pi_1(\Omega) \ar[r]^i & \pi_1^{Orb}(M) \ar[r]^{\mathcal{H}} & \Gamma \ar[r] & 0,
}
\]
where $i$ is the  inclusion induced by the following commutative
diagram:
\[
\xymatrix{
\widetilde M\ar[ddr]_p \ar[dr]^{\mathcal{D}}\ar[rr]^g  & & \widetilde M\ar[dl]_{\mathcal{D}}\ar[ddl]^p\\
                                                   &\Omega\ar[d]^q&\\
                                                   &M&\\
}
\]
\end{lemma}

\begin{proof}
Since $\mathcal{ H}$ is an epimorphism we only have to show that
$Ker(\mathcal{H})=\pi_1(\Omega)$.
Let $g\in
\pi_1(\Omega)$, then   $\mathcal{D}(g(\widetilde
m))=\mathcal{D}(\widetilde m)=Id(x)$. By the definition of $
\mathcal{H}$ we conclude that $\pi_1(\Omega)\subset
Ker(\mathcal{H})$.

Conversely,  let $g\in
Ker(\mathcal{H})$. Then $\mathcal{D}\circ g=\mathcal{D}$. Since
$\mathcal{D}$ is a covering and $\widetilde M$ is simply connected  we get 
$Ker(\mathcal{H})\subset \pi_1(\Omega)$.
\end{proof}



\vskip.4cm

\section{Some Properties of Complex Kleinian Groups} \label{s:propklei}
\subsection{A Characterization of finite groups}

In this section we prove: 

\begin{proposition} \label{c:clasfin}
Let $\Gamma\subset \PSL(3,\mathbb{C})$ be a discrete group. Then the
following conditions are equivalent:
\begin{enumerate}
\item $\Gamma$ is finite.
\item The order  $o(\gamma)$ is finite for all $\gamma\in \Gamma$.
\item $\Gamma$ acts  properly  discontinuously on $\mathbb{P}^2_{\mathbb{C}}$.
\end{enumerate}
\end{proposition}

The proof of (\ref {c:clasfin}) is based on
 Selberg's lemma and it  takes the rest of this subsection. 
  

\begin{lemma}\label{l:friendly}
Let  $(A_m)\subset PSL(3,\C)$  be a  sequence of finite groups such that  each $A_m$  properly contains $A_{m-1}$. Assume there exists  $k_0 \in \Bbb{N}$ and a sequence $(B_m)$  such that each $B_m$ is a commutative normal subgroup of 
$A_m$ and 
$
Card(A_m/B_m)=k_0  \; \textrm{for every   } m\in \Bbb{N}.$
Then  $ \langle \bigcup _{m\in \Bbb{N} } A_{m}\rangle$ contains an infinite  commutative  group.
\end{lemma}

\begin{proof}
Set   $n_m=max\{o(g): g\in
B_m\}$ where    $o(g)$  is the order of $g$.  We have the following possibilities:


{\it Case} 1.-  The sequence   $(n_m)_{m\in \mathbb{N}}$ is not bounded- Then we can find a subsequence of  $(A_m)_{m\in\mathbb{N}}$,  also denoted   $(A_m)_{m\in\mathbb{N}}$,  satisfying   $n_1>k_0$ and
$k_0n_m<n_{m+1}$ for every  $m$.  For each  $m$ let  $\gamma_m\in
B_m$ be such that $o(\gamma_m)=n_m$. Then
$\gamma_m^{k_0}\in\bigcap_{m\leq j} B_j$, 
$\gamma_l^{k_0}\neq \gamma_j^{k_0}$  whenever $l\neq j $ and
 $\gamma_j^{k_0} \gamma_l^{k_0}=\gamma_l^{k_0} \gamma_j^{k_0}$    for every  $l,j$. Hence the group  $\langle\gamma_m^{k_0}:m\in \mathbb{N}\rangle$ 
 is an infinite  commutative  subgroup of   $ \langle \bigcup _{m\in \Bbb{N} } A_{m}\rangle$.

\vv

{\it Case} 2.- The sequence   $(n_m)_{m\in \mathbb{N}}$ is  bounded-  In this case there exist   $c_0>1$ and  a subsequence of $(A_m)_{m\in\mathbb{N}}$, still denoted $(A_m)_{m\in\mathbb{N}}$,\
such that   
$n_m=c_0$ for all $m$. 
We claim that in this case there exists a sequence   $(w_n)\subset G$  and a sequence of subsequences of 
 $(A_m)$,  denoted $((A_{nm}))$, so that if we consider the corresponding sequence of sequences 
 of  $(B_m)$,  denoted   $((B_{nm}))$, then we have:

\begin{enumerate}
\item  $(A_{i+1,m})\subset (A_{im})$ for all $i$;
\item the inequality $Card(B_{j1}) >k_0c_0^{3+j} $ holds;  and
\item 
$w_j\in \bigcap_{m\in \Bbb{N}} B_{jm}\setminus \langle Id, w_1,\ldots, w_{j-1}\rangle$.
\end{enumerate}

\vv \noindent
The construction of these sequences is by an inductive process. 
 Assume first    $Card( B_1)>k_0c_0^3 $. For each 
$m>1$ we define a map:
\[
\phi_{1,m}:B_1\longrightarrow A_m/B_m \;.
\]
\[
l\mapsto
B_m l
\]
 Since we have   $Card( B_1) > Card(  A_m/B_m) =k_0$ we conclude that  $\phi_{1,m}$ is not injective. Hence there is an element  $w_1\neq id$ and a subsequence  $(A_{1m})_{n\in\mathbb{N}}\subset
(A_m)_{n\in\mathbb{N}}$, with its corresponding  subsequence  $(B_{1m})\subset (B_m)$, such that   $Card(B_{11}) >k_0c_0^4 $ and
$w_1\in \bigcap_{m\in \mathbb{N}} B_{1m}$.\\

Now consider, by induction, that one has the following sequence of transformations:
\[
\phi_{j+1,m}: B_{j1}/\langle Id, w_1,\ldots, w_{j}\rangle\longrightarrow A_{jm}/B_{jm}
\]
\[
\langle Id, w_1,\ldots, w_{j}\rangle l\mapsto B_{jm}l.
\]    
Since one has  $Card(B_{j1})) >k_0c_0^{3+j}$ we conclude that   $\phi_{j+1,m}$ is not injective. Hence there exist an element $w_{j+1}\notin  \langle Id, w_1,\ldots, w_{j}\rangle$,  a subsequence  $(A_{j+1,m})\subset
(A_{jm})$ and the corresponding   subsequence  $(B_{j+1,m})\subset (B_{jm})$ such that  $Card(B_{j+1,1})) >k_0c_0^{j+4} $ and 
$w_{j+1}\in \bigcap_{m\in \mathbb{N}} B_{j+1,m}$.

This proves the above claim. It is now clear that the group
  $\langle w_m:m\in \mathbb{N}\rangle$ is an infinite  commutative subgroup of $ \langle \bigcup _{m\in \Bbb{N} } A_{m}\rangle$.
\end{proof}

For the next lemma we need the following theorem due to Jordan (see  \cite[Theorem 8.29 ]{rag}):

\begin{theorem} \label{t:jor}
For any  $n\in \mathbb{N}$  there exists an integer  $S(n)$ with the
following property: if   $G \subset {\rm GL}(n,\mathbb{C})$ is any finite
subgroup, then $G$ admits   an abelian normal   subgroup $N$ such
that  $Card( G)\leq S(n)Card ( N ) $.
\end{theorem}

\begin{lemma} \label{c:con}
 Let $G$ be an infinite  countable  subgroup of $ {\rm GL}(3,\mathbb{C})$. Then  there exists an
 infinite commutative subgroup $N$ of $G$.
\end{lemma}
\begin{proof}
Suppose the lemma is false, so that  each element in $G$ has finite order. Then by Selberg's lemma we have that  $G$  must have an infinite number of generators, say  $\{\gamma_m\}_{m\in \mathbb{N}}$. 
For each 
$m\in \Bbb {N}$ define
$\,A_m=\langle\gamma_1,\ldots \gamma_m\rangle.$
Again by Selberg's lemma we have that each of these groups 
 $A_m$ is  finite. Then
 Theorem \ref{t:jor} implies that each  $A_m$ has a commutative normal subgroup  $N(A_m)$ with finite index  and such that 
 $$Card(A_m)\leq S(3)Card( N(A_m))\,.$$ 
 It is clear that we can construct a subsequence 
 of  $(A_m)_{m\in\mathbb{N}}$,  also denoted  $(A_m)_{m\in\mathbb{N}}$,  which satisfies:
\[
Card(A_m)=k_0 Card( N(A_m))\, \textrm{ for some } k_0\in
\mathbb{N}\, \textrm{ and every  } m\in \Bbb{N}.
\]
The result now follows from Lemma \ref {l:friendly}.
\end{proof}


\begin{lemma} \label{l:com}
Let $N\subset \SL(3,\mathbb{C})$ be a commutative subgroup where every
element is diagonalizable. Then there exists $\tau \in \SL(3,\mathbb{C})$
such that every element in $ \tau N\tau^{-1}$ is  a diagonalizable
matrix.
\end{lemma}
\begin{proof}
 Let  $g=(\gamma_{ij})\in N\setminus \{Id\} $; we may assume that
 $g$ is a diagonal matrix with    $\gamma_{11}\neq \gamma_{33}\neq \gamma_{22} $.
  Now let $h=(h_{i,j})\in N\setminus \{Id\}$ be  a non-diagonal matrix. By
  comparing the  coefficients  in the equation  $gh=hg$ we deduce
  $h_{13}=h_{23}=h_{31}=h_{32}=0$ and $\gamma_{11}=\gamma_{22}$. Set
$
\widetilde h=
\left(
\begin{array}{ll}
h_{11} & h_{12} \\
h_{21} & h_{22} \\
\end{array}
\right) $, then there exists     $k\in \SL(2,\mathbb{C})$ such that
$k^{-1}\widetilde hk$ is a diagonal matrix with distinct eigenvalues.
Set $ C= \sqrt{det(K)^{-1}} \left(
\begin{array}{ll}
k & 0\\
0 & 1
\end{array}
\right), $ then    $C^{-1}gC=g$. To conclude observe that
comparing both sides of the equations   $xg=gx,\,
xc^{-1}hc=c^{-1}hcx$ for each   $x\in C^{-1}NC$, one  can deduce
that every  element in  $C^{-1}NC$  is a diagonal  matrix.
\end{proof}

\begin{lemma} \label{c:cinf}
Let $\Gamma\subset \PSL(3,\mathbb{C}) $ be an  infinite discrete  group.
Then there exists an element $\gamma\in \Gamma$ with infinite order.
\end{lemma}

\begin{proof} By Lemma \ref{c:con},  $\Gamma$ contains  an infinite commutative
subgroup  $N$.  If  $o(\gamma)<\infty $ for every  $\gamma\in N$,
then every  element in $\Gamma$ has a lift which is  diagonalizable. By Lemma \ref{l:com} there exists an element $\tau \in
\PSL(3,\mathbb{C})$ such that $\tau N \tau^{-1}$ is a group where
every element has a lift which  is a diagonal matrix   whose eigenvalues are roots
of the unity. Therefore  $N$ is non-discrete, which is a contradiction.
\end{proof}

This completes the proof of Proposition \ref{c:clasfin}.

\vskip.3cm

\subsection{Complex lines and projective Groups} The following result can be
proved by  standard arguments, see \cite{kulkarni}:

\begin{lemma}  \label{l:limite} Let $\Gamma\subset \PSL(3,\mathbb{C})$ be a group acting properly discontinuously on  an open set $\Omega$. Then $L_0(\Gamma)\cup
L_1(\Gamma)\subset \mathbb{P}^2_{\mathbb{C}}\setminus  \Omega$. Moreover,  for every
compact set $K\subset \Omega$ the set of cluster
points of $\Gamma K$ is contained in $\mathbb{P}^2_{\mathbb{C}}\setminus 
\Omega$.
\end{lemma}

One has:

\begin{proposition} \label{l:line}
Let $\Gamma\subset \PSL(3,\mathbb{C})$ be an infinite  group  acting properly  
discontinuously on  an open set $\Omega$.  Then $\mathbb{P}^2_{\mathbb{C}}\setminus  \Omega$ contains at least one complex line $\ell$.
\end{proposition}
\begin{proof} By Lemma \ref{c:cinf} there exists an element $\gamma\in \Gamma$
with infinite order. Let $\widetilde\gamma\in \SL(3,\mathbb{C})$ be a lift of $\gamma$.
By the normal Jordan form theorem it is enough to consider  the following 4  cases :\vv

i) There exists a lifting $\widetilde\gamma$ of $\gamma$  given by:
\begin{equation}  \label{e:fornojor1}
\widetilde{\gamma}=
\left(
\begin{array}{llc}
1 & 1 & 0\\
0 & 1 & 1\\
0 & 0 & 1\\
\end{array}
\right),
\end{equation}
Notice  that one has:
\begin{equation}\label{e:juni} 
\widetilde{\gamma}^n=
\left(
\begin{array}{llc}
1 & n & {n(n-1)}/{2}\\
0 & 1 & n\\
0 & 0 & 1\\
\end{array}
\right).
\end{equation}
We claim that  the line $\overleftrightarrow{e_1,e_2} $ is contained in 
$\mathbb{P}^2_{\mathbb{C}}\setminus  \Omega$.  Otherwise there exists
$z\in\mathbb{C}^*$ such that $[z:1:0]\in \Omega$. 
For every   $\varepsilon\in
\mathbb{C}$, define 
 $$a_n(\varepsilon):=\left [z:1:
\frac{2(\varepsilon - n)}{n(n-1)}\right ]\in 
\Omega.$$ 
It is clear that the sequence $\{a_n(\varepsilon)\}$ converges to 
 $ [z:1:0]$ as $n$ tends to $\infty$. This implies that there exists $m(\epsilon)\in \Bbb{N}$ such that  $a_n(\epsilon)$ is in $\Omega$ whenever  $n>m(\epsilon)$. Hence  the set $$K(\epsilon)=\{a_n(\epsilon): n>m(\epsilon) \}\cup   \{[z:1:0]\}$$ is a compact subset of $\Omega$. Then Lemma \ref{l:limite} applied to  $K(\epsilon)$ and $\Omega$ yields that the cluster points of 
  $ \gamma^n(K)$ are contained in $\PC^2\setminus \Omega$. Now a straightforward computation shows that  $\gamma^n(a_n(\epsilon))$ converges to $[-z-\epsilon:1:0]$ and therefore
    $[-z-\epsilon:1:0]\in \PC^2\setminus \Omega$.  Since  $\epsilon$ is arbitrary  we get  $W=\{[w:1:0]:w\in \C\}\subset \PC^2\setminus \Omega $. And since   $\PC^2\setminus \Omega$ is a closed set, we get $\overline W=\overleftrightarrow{e_1,e_2} \subset \PC^2\setminus \Omega$, which contradicts our assumptions. Thence  $\mathbb{P}^2_{\mathbb{C}}\setminus  \Omega$ contains the line $\overleftrightarrow{e_1,e_2} $,  proving \ref{l:line} in this case.
 
    \vv

ii) Assume now that   $\widetilde\gamma$ is given by:
\[
\widetilde{\gamma}=
\left(
\begin{array}{lll}
\lambda & 1  & 0\\
0         &\lambda        & 0\\
0         & 0               & \lambda\\
\end{array}
\right),
\]
for some $\lambda\in \mathbb{C}$. Taking $\widetilde\gamma^{-1}$ if necessary,  we can  assume that $\vert \lambda\vert<1$.

Now we claim that
$\overleftrightarrow{e_1,e_2}\subset \mathbb{P}^2_{\mathbb{C}}\setminus  \Omega$
or $\overleftrightarrow {e_1,e_3}\subset \mathbb{P}^2_{\mathbb{C}}\setminus 
\Omega$. If $\overleftrightarrow {e_1,e_2}\not\subset
\mathbb{P}^2_{\mathbb{C}}\setminus  \Omega$ (the other case being similar) then  there exists $z\in\mathbb{C}^*$ such
that $[z:1:0]\in \Omega$. Observe that  for each $w\in \mathbb{C}^*$
we have  $$[wz:w: n\lambda^{3n-1}] \xymatrix{ \ar[r]_{n \rightarrow
\infty}&}[z:1:0]\,,$$ and 
$$\gamma^{n}\left (\left [z:1: \frac{n\lambda^{3n-1}}{w}\right ]\right)=\left
[\frac{wz\lambda}{n}+w: \frac{w\lambda }{n}:1\right ]\xymatrix{
\ar[r]_{n \rightarrow  \infty}&} [w:0:1]\,,$$ for all $ w\in
\mathbb{C}^*.$ Then   we  get
$\overleftrightarrow{e_1,e_3}\subset\mathbb{P}^2_{\mathbb{C}}\setminus  \Omega$ as   in the previous case. The other case is similar.\vv

iii) The third case is when $\widetilde{\gamma}=(\gamma_{ij})$ is a diagonal matrix
with $\vert \gamma_{11}\vert<\vert \gamma_{22}\vert<\vert
\gamma_{33}\vert $ and $ \gamma_{11} \gamma_{22} \gamma_{33}=1$. We
claim that $\overleftrightarrow {e_1,e_2}\subset
\mathbb{P}^2_{\mathbb{C}}\setminus  \Omega$ or $\overleftrightarrow
{e_2,e_3}\subset \mathbb{P}^2_{\mathbb{C}}\setminus  \Omega$. Suppose that there
exists $z\in\mathbb{C}^*$ such that $[z:1:0]\in \Omega$, then
$$[zw\gamma_{33}^n:w\gamma_{33}^n : \gamma_{22}^n]\xymatrix{ \ar[r]_{n
\rightarrow  \infty}&} [z:1:0]\,,$$
 for each $w\in \mathbb{C}^*$, and
 \[ \gamma^{n}([zw\gamma_{33}^n:w\gamma_{33}^n : \gamma_{22}^n])=
[zw\gamma_{11}^n\gamma_{22}^{-n}:w :1]\xymatrix{ \ar[r]_{n \rightarrow
\infty}&} [0:w:1]\,,\] for all $w\in \mathbb{C}^*.$
 Therefore $\overleftrightarrow{e_2,e_3}\subset\mathbb{P}^2_{\mathbb{C}}\setminus  \Omega$.\\

iv) Finally, if $\widetilde{\gamma}$ has another normal Jordan form,  Theorem
\ref{t:pciclic} ensures that $L_0(\gamma)\cup L_1(\gamma)$ contains a complex
line.
\end{proof}

\begin{theorem} \label{c:maxc} \label{l:invlin} 
Let $\gamma\in \PSL(3,\mathbb{C}) \setminus Id$ and let $\widetilde\gamma$ be a lift of
$\gamma$. The maximal open sets where  $\langle\gamma\rangle$ acts properly  
discontinuously and the set of $\langle\gamma\rangle$-invariant lines  are given by in the tables below {\rm (in terms of 
the Jordan normal form of
$\widetilde\gamma$):}

\begin {center}

\begin{tiny}
\begin{center}
\begin{tabular}{|c|l|l|}
\hline $
\begin{array}{c}
\\
\textrm{Normal Form}\\
\textrm{of } \widetilde\gamma\\
\\
\end{array}
$ & $
\begin{array}{c}
\textrm{Condition over the }
\lambda\textrm{'s} \end{array} $ & Maximal Regions of Discontinuity\\
\hline \multirow{3}{*}{ 
$ 
\left (
\begin{array}{lll}
\lambda_1  & 0            & 0\\
0            &  \lambda_2 & 0\\
0            &  0            & \lambda_3\\
\end{array}
\right) $ } & $ \vert \lambda_1\vert<\vert\lambda_2\vert<\vert\lambda_3\vert 
$ 
& 

$
\begin{array}{c}
\\
\mathbb{P}^2_{\mathbb{C}}\setminus (\overleftrightarrow{e_1,e_2}\cup \{e_3\});\\\\
\mathbb{P}^2_{\mathbb{C}}\setminus (\overleftrightarrow{e_3,e_2}\cup \{e_1\}) \\\\
\end{array}
$\\
\cline{2-3} & Other Case
&
$
\begin{array}{l}
\\
\Omega_{\Kul}(\langle\gamma\rangle)\\
\\
\end{array}$
\\
\hline
\multirow{3}{*}{ $ \left (
\begin{array}{lll}
\lambda    & 1            & 0\\
0            &  \lambda & 0\\
0            &  0              & \lambda^ {-2}\\
\end{array}
\right) $ } & $\vert \lambda \vert\neq 1 
$
&
$
\begin{array}{l}
\\
\mathbb{P}^2_{\mathbb{C}}\setminus (\overleftrightarrow{e_1,e_2}\cup \{e_3\});\\\\
\mathbb{P}^2_{\mathbb{C}}\setminus \overleftrightarrow{e_3,e_1} \\\\
\end{array}
$
\\
\cline{2-3}
&
Other Case 
&
$
\begin{array}{l}
\\
\Omega_{\Kul}(\langle\gamma\rangle)\\
\\
\end{array}
$
\\
\hline
$ 
\begin{array}{l}
 \\\\
\left (
\begin{array}{lll}
1  & 1 & 0\\
0  & 1 & 1\\
0  & 0 &1\\
\end{array}
\right)\\\\
\end{array}
$ &  & $
\begin{array}{l}
\\
\Omega_{\Kul}(\langle\gamma\rangle)\\
\\
\end{array}$
\\
\hline 

\end{tabular}

\end{center}
\end{tiny}

\end{center}

\begin{tiny}
\begin{center}
\begin{tabular}{|c|l|l|}
\hline $
\begin{array}{c}
\\
\textrm{Normal Form}\\
\textrm{of } \widetilde\gamma\\
\\
\end{array}
$ & $
\begin{array}{c}
\textrm{Condition over the }
\lambda\textrm{'s} \end{array} $ & Invariant Lines\\
\hline 
\multirow{3}{*}{ 
$ 
\left (
\begin{array}{lll}
\lambda_1  & 0            & 0\\
0            &  \lambda_2 & 0\\
0            &  0            & \lambda_3\\
\end{array}
\right) $ } & $ \lambda_1\neq \lambda_2\neq \lambda_3 \neq \lambda_1
$ 
& 

$
\begin{array}{l}\\
\{\overleftrightarrow{e_1,e_2},\, \overleftrightarrow{e_1,e_3}.\,
\overleftrightarrow{e_3,e_2} \}\\
\\
\end{array}
$\\
\cline{2-3} & $
\lambda_1=\lambda_2\neq \lambda_3
$ &
$\begin{array}{l}
\\
\{\overleftrightarrow{e_1,e_2}\}\cup\{\overleftrightarrow{e_3,p} \; \big \vert \;
p \in \overleftrightarrow{  e_1,e_2} \}\\
\\
\end{array}$
\\
\hline
\multirow{3}{*}{ $ \left (
\begin{array}{lll}
\lambda    & 1            & 0\\
0            &  \lambda & 0\\
0            &  0              & \lambda^ {-2}\\
\end{array}
\right) $ } & $\lambda^3\neq 1 
$
&
$
\begin{array}{l}
\\
\{\overleftrightarrow{e_1,e_3},\, \overleftrightarrow{e_1,e_2} \}\\\\
\end{array}
$
\\
\cline{2-3}
&
$\lambda^3=1$
&
$
\begin{array}{l}
\\
\{\overleftrightarrow{p,e_1}\vert p\in \overleftrightarrow { e_3,e_2}\}
\\\\
\end{array}
$
\\
\hline 
$ 
\begin{array}{l}
 \\\\
\left (
\begin{array}{lll}
1  & 1 & 0\\
0  & 1 & 1\\
0  & 0 &1\\
\end{array}
\right)\\\\
\end{array}
$ &  & $
\begin{array}{l}
\\
\{\overleftrightarrow{e_1,e_2}\}
\\
\\
\end{array}
$
\\
\hline
\end{tabular}

\end{center}
\end{tiny}

\end{theorem}



\begin{proof} Notice that
by the normal Jordan form theorem  we  only need to consider the three
cases indicated in the left side of  these tables. We look first at the maximal regions of discontinuity, {\it i.e.}, the first table. This is  based on Navarrete's work in \cite{pablo2}. 


Consider an element  $\gamma \in PSL(3,\C)$ of infinite order and let  $\Omega \subset \Bbb{P}^2_{\C} $ be an open invariant set where $G=\langle \gamma \rangle$ acts  properly discontinuously. We know from Theorem \ref{t:pciclic}  that  $\Lambda_{Kul}(G)$ contains  either one or two  lines, or else it is all of $\PC^2$. We consider each of these possible cases independently. Notice that 
$\Lambda_{Kul}(G) = \PC^2$   if and only if $\gamma$ has a lift in ${\rm \SL}(3,\C)$ whose normal form is
\[
\left (
\begin{array}{lll}
\alpha &0 & 0\\
0 &\beta & 0\\
0 &0 & \tau\\
\end{array}
\right )
\]
with $\vert \alpha\vert=\vert \beta\vert=\vert \tau\vert=1 $, so we focus in the remaining cases.

\vv
{\bf Case 1.}  The Kulkarni set $\Lambda_{Kul}(G)$ has exactly one line $\ell$. In this case Theorem \ref{t:pciclic} says that $\gamma$ has a lift in ${\rm \SL}(3,\C)$ whose normal form is:

\[
\left (
\begin{array}{lll}
1 &1 & 0\\
0 &1 & 1\\
0 &0 & 1\\
\end{array}
\right )
\;, \; \hbox{ or } \;
\left (
\begin{array}{lll}
\alpha &0 & 0\\
0 &\alpha & 0\\
0 &0 & \alpha^{-2}\\
\end{array}
\right )
\hbox{ with } \vert \alpha \vert\neq 1 \;, \hbox { or }  \;
\left (
\begin{array}{lll}
\alpha &1 & 0\\
0 &\alpha & 1\\
0 &0 & \alpha^{-2}\\
\end{array}
\right )
\, \hbox { 
with } \vert \alpha \vert = 1 \,.\]
In all these cases we know from the proof of Proposition \ref {l:line}  that the line $\ell$ must be contained in 
 $\PC^2\setminus \Omega$. Now we observe that $B=\Lambda_{Kul}(G)\setminus \ell$ is either  empty or it contains one point. If $B$ is empty, we have finished the proof. Otherwise, if $B$ is a point, then Theorem \ref{t:pciclic}  implies that this point is in  $L_0(G)$, which is contained in $\PC^2\setminus \Omega$, which proves the theorem in this case.

\vv

{\bf Case 2.}  If  $\Lambda_{Kul}(G)$ has two distinct  lines $\ell_1$ and $\ell_2$ then we know that the lifts of $\gamma$ are of the remaining types in Theorem \ref{t:pciclic}. Then from the proof of Proposition \ref {l:line}  
we know that if 
 $\ell_2$ is not contained in $\PC^2\setminus \Omega$, so $\ell_1$  necessarily is in $\PC^2\setminus \Omega$, and viceversa. 
 The proof of Proposition \ref {l:line}  also  shows that for every compact set $K \subset  \ell_2 \setminus L_0(G)$   the set of cluster points of the orbit  $\Gamma K$ is contained in $\ell_1$, and viceversa. This completes the proof of this part of Theorem \ref{c:maxc}.

\vv

Now we look at the second table, {\it i.e.}, we determine the invariant lines in each case.

\vv

{\bf Case} 1. -$\widetilde\gamma$ is diagonalizable- Let us 
consider the following possibilities:\vv

Possibility 1. -Every eigenvalue has multiplicity 1- Let $\{u,v,w\}$ be
the set of fixed points of $\gamma$ and $\ell$  an invariant
line  under $\gamma$. We claim that
$\ell$ is one of the lines $\overleftrightarrow{u,v} $, $\overleftrightarrow{w,v}$, 
or $\overleftrightarrow{u,w}$.  Assume on the contrary  that
$\ell$ is different from these lines, so  $\ell$ contains
exactly one fixed point, say $u$.   Let
$K=\overleftrightarrow{v,w}$.  Then $K\cap \ell$ contains exactly one point, say $*$, and $*\neq u$. Then   $*$
is a fixed point, which is a contradiction.\vv

Possibility 2. -One eigenvalue has multiplicity 2-. After conjugating with a linear  transformation we can assume that $\{e_1,e_2,e_3\}$ is a basis of eigenvectors for $\widetilde\gamma$ such that $e_1$ and $e_2$ have the same eigenvalue. Then it follows  from Theorem
\ref{t:pciclic}  that
$Fix(\gamma)=\overleftrightarrow{e_1,e_2}\cup \{e_3\}$. Thus every
line  that contains $e_3$ is invariant. Let $\ell$ be an invariant
line  such that $e_3\notin \ell$ and $c\in \ell$, then
$\gamma(c)=\gamma(\overleftrightarrow{e_3,c}\cap \ell)=c$. Thus
$c\in \overleftrightarrow{e_1,e_2}$ and therefore
$\ell=\overleftrightarrow{e_1,e_2}$.\\

{\bf Case} 2. -$\widetilde\gamma$ has at most two  linearly
independent eigenvectors-  In this case  $\widetilde\gamma$  has the following
Jordan's normal form:
\[
\left(
\begin{array}{lll}
\lambda & 1 & 0\\
0& \lambda & 0\\
0 & 0 & \lambda^{-2} \\
\end{array}
\right).
\]

Let us consider the following possibilities:

Possibility 1. -$\lambda^3 \neq 1$- Let  $\ell$ be an invariant line.
We claim that $\ell= \overleftrightarrow{e_1,e_2}$ or $\ell=
\overleftrightarrow{e_1,e_3} $. Assume  that $\ell$ is different
from these lines, so there exists  $[w]=[z_1:z_2:z_3]\in \ell$ such that
$z_2z_3\neq 0$. Consider the  equation $$ 0=\alpha_1 w+\alpha_2
\widetilde\gamma(w)+\alpha_3 \widetilde\gamma^2(w).$$ 
This equation is equivalent to the system:
\[
\begin{array}{l}
\alpha_1 +\alpha_2\lambda^{-2} +\alpha_3\lambda^{-4}=0\\
\alpha_1 +\alpha_2\lambda +\alpha_3\lambda^2=0\\
\alpha_2+2\alpha_3\lambda=0.\\
\end{array}
\]
Since the determinant of the system  is
$(\lambda^{-2}-\lambda)^2\neq 0$,  we get $\alpha_1=\alpha_2=\alpha_
3=0$. Therefore $[w]$, $[\widetilde\gamma(w)]$ and $[\widetilde
\gamma^2(w)]$ are not contained in a complex line, which is a
contradiction.\vv

Possibility 2. -$\lambda^3= 1$- Let  $\ell$ be an invariant line.  Assume
that $\ell\neq \overleftrightarrow{e_1,e_2}$ and $\ell\neq
\overleftrightarrow{e_1,e_3} $. Then  there exists a point
$[w]=[z_1: z_2: z_3]\in \ell $ such that $z_2 z_3\neq 0$.
Consider the equation $$ 0=\alpha_1 w+\alpha_2\widetilde
\gamma(w)+\alpha_3 e_1.$$  This 
equation is equivalent to the system:
\[
\begin{array}{l}
\alpha_2 z_{2} +\alpha_3=0\\
\alpha_1+\alpha_2=0,\\
\end{array}
\]
which has the non-trivial solution
$\alpha_1=z_2^{-1},\,\alpha_2= -z_2^{-1},\, \alpha_3=1$. Hence 
$[w]$, $[\widetilde\gamma(w)]$ and $e_1$ lie  in the same line. Since
$\ell$ is invariant we conclude that $e_1\in \ell$.\vv

{\bf  Case} 3.-  $\widetilde\gamma$ has the normal form:
\[
\left(
\begin{array}{lll}
1 & 1 & 0\\
0 & 1 & 1\\
0 & 0 & 1\\
\end{array}
\right).
\]
Let $\ell$ be  an invariant line  and assume  that $\ell\neq
\overleftrightarrow{e_1,e_2}$. Then there exists
$[w]=[z_1:z_2:z_3]\in \ell$ with $z_3\neq 0$. Since  the
equation $\alpha_1 w+\alpha_1
\widetilde\gamma(w)+\alpha_3\widetilde\gamma^2(w)=0$  is   equivalent to
the system:
\[
\begin{array}{l}
\alpha_1 +\alpha_2 +\alpha_3=0\\
\alpha_2 +2\alpha_3=0\\
\alpha_3=0 \;,\\
\end{array}
\]
and  such system has only the trivial solution, we conclude that
$[w], [\widetilde\gamma(w)], [\widetilde\gamma^2(w)]$ are not
contained in a complex line, which  contradicts the initial
assumption.
\end{proof}

\vskip.4cm

\subsection{Controllable Groups} \label{ss:cont}

Recall that controllable groups were defined in \ref{definitions}. The simplest type of such groups are the suspensions:

\begin{example}
[Suspension with a group]  Let $\Gamma\subset \PSL(2,\C)$ be a discrete
group with non-empty discontinuity region,  $G\subset \mathbb{C}^*$ a
discrete group and $i:\SL(2,\C)\rightarrow \SL(3,\mathbb{C})$ 
the inclusion given by: $i(h)= \left (
\begin{array}{ll}
h & 0\\
0 & 1
\end{array}
\right ).$ The suspension of $\Gamma$ with respect to $G$, denoted
$Susp(\Gamma,G)$, is defined by:
\[
Susp(\Gamma,G)=
\left
\langle
\{i(h):h\in [\Gamma]^{-1}\},
\left \{
\left (
\begin{array}{lll}
g & 0 & 0\\
0 & g & 0\\
0 & 0 & g^{-2}\\
\end{array}
\right ):g\in G
\right\}
\right
\rangle,
\]
where $[ \;]$ is the projection $\SL(2,\C) \to \PSL(2,\C)$.
Observe that if  $G=\{\pm 1\}=\mathbb{Z}_2$, then 
 $Susp(\Gamma,\mathbb{Z}_2)$ coincides with the double suspension of
$\Gamma$ defined in \cite{pablo2},  and   when $[\Gamma]^{-1}_2$
contains a subgroup $\widetilde\Gamma$ for which  $[\widetilde
\Gamma]=\Gamma$, then  $i(\widetilde\Gamma)$ coincides with the
suspension of $\Gamma$ defined in \cite{sv1}. Notice that the line 
$\overleftrightarrow{e_1,e_2}$ is invariant. Let 
$\Lambda(\Gamma)$  be the usual
limit set for  the  action of $Susp(\Gamma,G)$ on
$\overleftrightarrow{e_1,e_2}$, which coincides with the
action of $\Gamma$ on $\mathbb{P}^1_{\mathbb{C}}$.  Then, just as in 
 \cite{pablo2}, one has:
\[
\Lambda_{\Kul}(Susp(\Gamma,G))=
\left \{
\begin{array}{ll}
\bigcup_{p\in \Lambda(\Gamma)}\overleftrightarrow{p,e_3} & \textrm{if } G
\textrm{ is finite;}\\
\bigcup_{p\in \Lambda(\Gamma)}\overleftrightarrow{p,e_3}\cup
\overleftrightarrow{e_1,e_2} & \textrm{if } G \textrm{ is infinite.}
\end{array}
\right.
\]
\end{example}

Notice that if $ \Gamma\subset \PSL(3,\mathbb{C})$ is a subgroup, $p$ is a $\Gamma$-invariant point in  $\mathbb{P}^2_{\mathbb{C}}$  and $\ell$ is a complex projective  line not containing $p$, then we have a 
natural projection map
$\pi_{p,\ell}:\mathbb{P}^2_{\mathbb{C}}- \{p\}\longrightarrow
\ell$  given by $\pi_{p,\ell}(x)=\overleftrightarrow{x,p}\cap \ell$. In other words, for every $x \in \PC^2 \setminus \{p\}$ we have a well-defined projective line $\overleftrightarrow{x,p}$, and this line meets $\ell$ in a unique point, which by definition is $\pi_{p,\ell}(x)$.
We observe that $\pi_{p,\ell}$ induces a 
map
$\Pi_{p,\ell}:\Gamma\longrightarrow  Bihol(\ell) $
given by $\Pi(g)(x)=\pi(g(x))$. For simplicity we set $\pi = \pi_{p,\ell}$ and $\Pi = \Pi_{p,\ell}$.
We have:

\begin{lemma} \label{l:control}
Let   $ \Gamma\subset \PSL(3,\mathbb{C})$ be a subgroup, $p\in \mathbb{P}^2_{\mathbb{C}}$
such that  $\Gamma p=p$ and $\ell$ a complex line not containing $p$.
Define $\pi$ and $\Pi$ as above.
Then:
\begin{enumerate}
\item \label{i:con1} $\pi$ is a holomorphic  function. 
\item \label{i:con2} $\Pi$   is a group morphism. 
\end{enumerate}
\end{lemma}

\begin{proof}  Let us show (\ref{i:con1}). If $p=e_3$ and $\ell=\overleftrightarrow{e_1,e_2}$,  then
$\pi([z:w:x])=[z:w:0]$ which
is  holomorphic. If this is not the case, take $g\in
\PSL(3,\mathbb{C})$ such that $g(p)=e_3$ and
$g(\ell)=\overleftrightarrow{e_1,e_2}$. Then for each 
$h\in\Gamma$ and each $ x\in \ell$, one has:

\begin{equation} \label{e:idline}
\begin{array}{l}
\overleftrightarrow{h(x),p}=\overleftrightarrow{\Pi_{\ell,p}(h(x)),p} \;;  \hbox{ and}\\
\overleftrightarrow{h(x),p}=h(\overleftrightarrow{x,p}) \;.
\end{array}
\end{equation}
Hence:
\begin{equation} \label{e:hol}
\pi(x) =g^{-1}( \pi_{e_3,\overleftrightarrow{e_1,e_2}}( g(x))) \,,
\end{equation}
so $\pi$ is holomorphic.

Now we show (\ref{i:con2}).  Step 1.- $\Pi(\Gamma)\subset Bihol(\ell)$- It is enough to observe
that $\Pi_{}(g)=\pi\circ g\vert_{\ell}$.

\medskip

Step 2.-$\Pi$ is a group morphism.  From equation
\ref{e:idline} we get : 
$$\Pi(g)\circ \Pi(h)(x)=\overleftrightarrow{g(\Pi(h)(x)),p}\cap
\ell=g(\overleftrightarrow{(h(x),p})\cap \ell=\Pi(g\circ
h)(x),$$
thus proving the statement.
\end{proof}

Now we have:

\begin{theorem}\label{l:control2}  Let   $ \Gamma\subset \PSL(3,\mathbb{C})$ be a discrete subgroup, $p\in \mathbb{P}^2_{\mathbb{C}}$
such that  $\Gamma p=p$ and $\ell$ a complex line not containing $p$. Define $\Pi$ as above.
\begin{enumerate}
\item \label{i:con3}
If   $ Ker(\Pi)$ is finite and $\Pi(\Gamma)$ is discrete, then
$\Gamma$ acts properly   discontinuously on $\Omega:=   \bigcup_{z\in
\Omega (\Pi (\Gamma))}\big(\overleftrightarrow{z,p} \setminus  \{p\}\big)$, where
$\Omega(\Pi(\Gamma)$ denotes the discontinuity set of
$\Pi(\Gamma)$.
 \item \label{i:con4} If  $\Pi(\Gamma)$
is non-discrete and  $\ell$ is invariant, then $\Gamma$ acts  properly
discontinuously on $\Omega= \Big( \bigcup_{z\in
\Eq(\Pi(\Gamma))} \overleftrightarrow{z,p}\Big) \setminus (\ell\cup\{p\} )$.
\end{enumerate}
\end{theorem}

\begin{proof}  
Let us show (\ref{i:con3}).  Let $K\subset \Omega$ be a compact set and define 
$\,K(\Gamma)=\{\gamma\in\Gamma: g (K)\cap K\neq \emptyset\}\,.$ 
Assume that  $K(\Gamma)$ is infinite. Let  $(\gamma_n)_{n\in
\mathbb{N}}$ be an enumeration of $K(\Gamma)$. Since $Ker(\Pi)$ is
finite, there exists a
subsequence of $(\gamma_n)$, still denoted $(\gamma_n)$,  such that   $\Pi(\gamma_k)\neq \Pi(\gamma_l)$ whenever  $l\neq k$.
Therefore $$\{\Pi(\gamma_n):n\in \mathbb{N}\}\subset\{g\in \Pi(\Gamma):
g(\pi(K))\cap \pi(K)\neq \emptyset\},$$  which is a contradiction. Thus $\Gamma$ acts discontinuously on $\Omega$.\vv

Let us show (\ref{i:con4}). 
Assume that $\Gamma$ does not 
act properly  discontinuously. After conjugating with a projective transformation, if necessary, we can take   $p=e_3$ and  $l=\overleftrightarrow{e_1,e_2}$.  Then there exist
$k=[z:h:w]$, $q\in \Omega$, $(k_n)\subset \Omega$ and
$(\gamma_n=(\gamma^{(n)}_{ij})\subset
\Gamma$ a sequence of  distinct elements  such that $k_n\xymatrix{
\ar[r]_{n \rightarrow  \infty}&} k$ and $\gamma_n(k_n)\xymatrix{
\ar[r]_{n \rightarrow  \infty}&} q $.  By Corollary \ref{c:16} there exist a subsequence of $(\gamma_n)$, still denoted $(\gamma_n)$, and  a holomorphic map
$f:\Eq(\Pi(\Gamma))\longrightarrow \overline{\Eq(\Pi(\Gamma))}$ such
that $\Pi(\gamma_n)\xymatrix{ \ar[r]_{n \rightarrow \infty}&} f$
uniformly on compact sets of $\Eq(\Pi(\Gamma))$. Moreover, by the convergence property, either $f\in
Bihol(\ell)$   or else $f$ is a constant function $c\in
\partial \Eq(\Pi(\Gamma))$.  Since $\pi(\gamma_n(k_n))$ tends to
$ \pi(q)\in \Eq(\Pi(\Gamma))$
as $n $ tends to $\infty$, we conclude that  $f$ is non-constant. Therefore
there exist $\gamma_{11},\gamma_{12},\gamma_{21},\gamma_{22}\in \mathbb{C}$ such that $\gamma_{11}\gamma_{22}-\gamma_{21}\gamma_{12}=1$
and $$f[z:w:0]=[\gamma_{11}z+\gamma_{12}w:\gamma_{21}z+\gamma_{22}w:0].$$
 Since
$\gamma^{(n)}_{13}=\gamma^{(n)}_{23}=\gamma^{(n)}_{31}=\gamma^{(n)}_{32}=0$,
$\gamma_{33}^{(n)}(\gamma^{(n)}_{11}
\gamma^{(n)}_{22}-\gamma^{(n)}_{12}\gamma^{(n)}_{21})=1$, there exists a     subsequence of $(\gamma^{(n)}_{ij}) $, still denoted
$(\gamma^{(n)}_{ij}) $, such that:
$
\gamma_{ij}^{(n)}\sqrt{\gamma_{33}^{(n)}}\xymatrix{ \ar[r]_{n \rightarrow \infty}&} \gamma_{ij} \;.
$

{\it Claim} 1.- The  sequence  $(\gamma_{33}^{(n)})$ is bounded, for otherwise  there exists  a subsequence of $(\gamma_{33}^{(n)})$, still denoted $(\gamma_{33}^{(n)})$, such that   $\gamma^{(n)}_{33} \xymatrix{ \ar[r]_{n \rightarrow  \infty}&} \infty $ and $\gamma_n(k_n)\xymatrix{ \ar[r]_{n \rightarrow  \infty}&}
[0:0:1]$, which is a contradiction. Thence $(\gamma_{33}^{(n)})$ is bounded.\vv

{\it Claim} 2.- The  sequence  $(\gamma_{33}^{(n)})$ is a bounded distance away from $0$,  for otherwise  there exists  a   subsequence of $(\gamma_{33}^{(n)})$, still denoted 
$(\gamma_{33}^{(n)})$, such that  $\gamma^{(n)}_{33}\xymatrix{ \ar[r]_{n \rightarrow  \infty}&} 0$. Then  $\gamma_n(k_n)\xymatrix{ \ar[r]_{n
\rightarrow  \infty}&}
[\gamma_{11}z+a_{12}h:\gamma_{21}z+\gamma_{22}h:0]$, which is a contradiction. \vv

 It follows from the previous claims that   there exists  a subsequence of $(\gamma_{33}^{(n)})$, still denoted  $(\gamma_{33}^{(n)})$, and  $\gamma_{33}\in \mathbb{C}^*$,  such that 
$\gamma_{33}^{(n)}\xymatrix{ \ar[r]_{n \rightarrow  \infty}&}\gamma_{33}$. Hence:
\[
\left [
\begin{array}{lll}
\gamma_{11}^{(n)} & \gamma_{12}^{(n)} & 0\\
\gamma_{21}^{(n)} & \gamma_{22}^{(n)} & 0\\
0            & 0            & \gamma_{33}^{(n)}\\
\end{array}
\right ]
\xymatrix{ \ar[r]_{n \rightarrow  \infty}&}
\left [
\begin{array}{lll}
\gamma_{11}\sqrt{\gamma^{-1}_{33}} & \gamma_{12}\sqrt{\gamma^{-1}_{33}} & 0\\
\gamma_{21}\sqrt{\gamma^{-1}_{33}} & \gamma_{22}\sqrt{\gamma^{-1}_{33}} & 0\\
0                   & 0                   & \gamma_{33}\\
\end{array}
\right ]\in \PSL(3,\mathbb{C}).
\]
This is a contradiction since $\Gamma$ is discrete.
\end{proof}

\begin{remark} Notice that in Theorem \ref{l:control2}  we did not consider the case
when  $\Pi(\Gamma)$ is discrete and  $ Ker(\Pi)$ is infinite. That case does not  appear in the context of this article, though it  is interesting for other reasons  and it  is part of a forthcoming article.
\end{remark}

\vskip.3cm

\subsection{Quasi-cocompact Groups}

Recall from Definition \ref{definitions} that a group  acts cocompactly on an invariant set $\Omega \subset \PC^2$ if the quotient $\Omega/ \Gamma$ is compact. The group $\Gamma$ is quasi-cocompact if there exists an open invariant set $\Omega \subset \PC^2$ where the action is properly discontinuous and cocompact.

\begin{theorem} \label{c:coneccion}
Let  $\Gamma\subset \PSL(3,\mathbb{C})$ be an infinite  group acting  cocompactly
on an open proper subset $\Omega$ of $\PC^2$.  Then, up to projective equivalence,  $\Omega$ is of one of the following 7 types:

\begin{enumerate}

\item $\Omega= \mathbb{H}_{\mathbb{C}}^2$,
\item $\Omega=\mathbb{P}^2_{\mathbb{C}}\setminus (\overleftrightarrow{e_1,e_2}\cup
\{e_3\})$,
\item $\Omega=\mathbb{P}^2_{\mathbb{C}}\setminus  (\overleftrightarrow{e_1,e_2}\cup
\overleftrightarrow{e_1,e_3}\cup \overleftrightarrow{e_3,e_2})$,
\item $\Omega=\mathbb{P}^2_{\mathbb{C}}\setminus (\overleftrightarrow{e_1,e_2}\cup
\overleftrightarrow{e_1,e_3})$,
\item  $\Omega=\mathbb{P}^2_{\mathbb{C}}\setminus 
\overleftrightarrow{e_1,e_2}$,
\item $\Omega=\bigcup_{z\in\mathbb{H}}
\overleftrightarrow{z,e_1}\setminus \{e_1\}$,
\item There exists a domain  $D$ in $\overleftrightarrow{e_1,e_2}$ omitting at least 3 points  such that:  $$\Omega=\bigcup_{z\in D}
\overleftrightarrow{z,e_3}\setminus (\{e_3\}\cup
\overleftrightarrow{e_1,e_2}).$$
\end{enumerate}
\end{theorem}

\begin{proof} By Theorems \ref{t:thurston}, \ref{t:klingler} and Lemma
\ref{l:puente},  cases 1-6 correspond  to $\Omega=
\,\mathbb{H}_{\mathbb{C}}^2,$\, $ \mathbb{C}^2\setminus  \{0\},\,
\mathbb{C}^*\times\mathbb{C}^*,\, \mathbb{C}^*\times\mathbb{C},\,
\mathbb{C}^2, \,\mathbb{C}\times\mathbb{H}$. In the latter case, the
developing map $ \mathcal{D}:\mathbb{C}\times
\mathbb{H}\longrightarrow \mathbb{P}^2_{\mathbb{C}}$ is a covering transformation  given by
$(z,w)\mapsto [A(w):B(w):e^{-\mu z}] $ where $\mu \in
\mathbb{C}^*$ and $A,B:\mathbb{H}\longrightarrow\mathbb{C}$ are
holomorphic maps. 

We notice that
$A, B$ are not constants, since otherwise,  for each $z\in\mathbb{C}$  we have that $\mathcal{D}(\{z\}\times \mathbb{H})$ is a point, which is not possible because   
$\mathcal{D}$ is a local homomorphism. Also notice that 
 $A,B$ do not have a common zero, for otherwise,  if  $w_0\in\mathbb{H}$ is
 a common zero, then $\mathcal{D}(\mathbb{C}\times \{w_0\})=\{e_3\}$, which is a contradiction, since $\mathcal{D}$ is a local homeomorphism.\vv

From the previous claims follows that  $D=\{[A(z):B(z):0]: z\in \mathbb{H}\}$ is an open set on  $\overleftrightarrow{e_1,e_2}$. Moreover,  $\mathcal{D}(\mathbb{C}\times \mathbb{H})$ is biholomorphic to $D\times \mathbb{C}^*$. If $D$ is non hyperbolic, then its  universal covering  is either $\mathbb{C}$ or $\mathbb{P}^1$. Therefore the universal covering  of  $D\times \mathbb{C}^*$ is either $\mathbb{C}\times \mathbb{C}$ or $\mathbb{P}^1_\mathbb{C} \times \mathbb{C}$. Since $\mathcal{D}$ is a covering and $\mathbb
{H}\times \mathbb{H}$ is simply connected, it follows that $\mathbb
{H}\times \mathbb{H}$ is biholomorphic to either  $\mathbb{C}\times \mathbb{C}$ or $\mathbb{P}^1_\mathbb{C} \times \mathbb{C}$, which is not possible. Thus  $D$ is  as stated.
\end{proof}

\begin{definition}[Kulkarni \cite{kulkarni}]
Let $\Gamma\subset \PSL(3,\Bbb{C})$ be a group acting properly 
discontinuously on a domain $\Omega$ and $R$ a subset of $\Omega$. Then $R$ is said 
to be a {\it pre-fundamental region} for the action of $\Gamma$ on $\Omega$ if
the set $\{\gamma\in \Gamma:\gamma R\cap R \}$ is finite and $
\Gamma R=\Omega$. 
\end{definition}

\begin{corollary} \label{l:comp}  Let $\Gamma\subset \PSL(3,\Bbb{C})$
be a group acting properly discontinuously on an open set  $\Omega\subset \Bbb{P}^2$, then:
\begin{enumerate}
\item \label{i:co1}
If $\Gamma$ acts  on  $\Omega$ with compact quotient,  then $\Gamma$ is finitely  generated and:
\begin{enumerate}
\item  it has a compact pre-fundamental region for its action  on $\Omega$;
\item   every point  $p\in \partial \Omega$ 
is a cluster point of $\Gamma R$.
\end{enumerate}
\item  \label{i:co2}
If $\Omega $ is connected, $\Omega/\Gamma$ is compact  and  $\Omega_0\subset \Omega$ is a subdomain 
such that $\Omega_0/Isot(\Omega_0,
\Gamma)$ is compact, then $\Omega=\Omega_0$.
\item  \label{i:co3}
If $\Omega/\Gamma$ is connected  and $\Omega_0$ 
 is a connected component of 
$\Omega$, then $\Omega_0/Isot(\Omega_0,
\Gamma)=\Omega/\Gamma$.
\end{enumerate}
\end{corollary}

\begin{proof} That  $\Gamma$ is finitely generated follows from the proof of  Lemma \ref{c:vego}. From \cite{kulkarni} we know that there are always pre-fundamental regions, and it is clear that in our there exists a compact one, that we denote by $R$. Now let 
$p\in \partial \Omega $,  $n\in \Bbb{N}$ and  denote by  $B(p,n)$ the ball in $\PC^2$ with center at $p$ and radius $n^{-1}$
with respect to the 
Fubini-Study metric. Then there  exist  $k_n\in R$ and  $\gamma_n\in \Gamma$ such that $\gamma_n(k_n)\in B(p,n)$. This proves statement (1).

Assume now that statement (2)  does not hold. Since 
 $\Omega$ is connected we have that there exists  $q \in \partial \Omega_0\cap \Omega$.   Then by statement (1)  there  exist  sequences $(k_n)\subset \Omega_0$ and  $(\gamma_n)$ such that  $k_n\xymatrix{
\ar[r]_{m \rightarrow  \infty}&} k $ with $k\in \Omega_0$ and $\gamma_n(k_n)\xymatrix{
\ar[r]_{m \rightarrow  \infty}&} q $,  with $(\gamma_n)$ a sequence of distinct elements. This is a contradiction, so we get statement (2).

Now we prove 
(3).  Define a map   $\rho:\Omega_0/Isot(\Omega_0,\Gamma)\rightarrow \Omega/\Gamma$ by  $Isot(\Omega_0,\Gamma)x\mapsto  \Gamma x $. This map is obviously continuous and injective, so  we only need to show that $\rho$ is surjective. For this, let  $x\in \Omega/\Gamma$, then there exists $y \in \Omega_0$ such that  $x=\Gamma y$. Let $\Omega_1$
be the connected component of $\Omega_0$ that contains the point 
 $y$. Since $\Omega/\Gamma$ is connected we have that there exists
 $\gamma_1\in \Gamma$ such that $\gamma_1 y\in \Omega_0$, and the proof of (\ref{l:comp}) is complete.
\end{proof}




\vskip.3cm


\section{Quasi-cocompact Groups Arising From Inoue  Surfaces}\label{s:inoue}

As before, we set
$  \mathbb{H}^+=\{[0:z:1]:Re(z)> 0\}$.
In this section we prove:

\begin{theorem}\label{t: Inoue} 
Let  $\Gamma\subset \PSL(3,\C)$ be  discrete and $\Omega$ an open invariant set in $\PC^2$ where 
 $\Gamma$ acts  properly  discontinuously.  Assume there exists  a connected component $\Omega_0 \subset  \Omega$ such that
 $\Omega_0/Isot(\Omega_0,\Gamma)$ is compact and it has a possibly ramified finite cover $S$ which is an Inoue surface. Then:
\begin{enumerate}

\item \label{t: Inoue1}
The set $\Omega_0$ is  $\Omega_0 =\C \times \Bbb{H}^+$.
\item \label{t: Inoue2} The Kulkarni limit set is   $\Lambda_{Kul}=L_{0}(\Gamma)\cup L_1(\Gamma)=\partial \Omega_0$. This set splits $\PC^2$ in two connected components which are projectively equivalent and isomorphic to ${\mathbb H}^+ \times \C$.
\item\label{t: Inoue3}  The set $\Omega$ is contained in the Kulkarni set  $\Omega_{Kul}(\Gamma)$, which  is the largest open invariant set where the action is properly discontinuous.  
\item\label{t: Inoue4} The quotient  $\Omega_{Kul}(\Gamma)/\Gamma$ is compact and it is  either  $S$ or else consists of $S$  and another disjoint biholomorphic copy of $S$.
\end{enumerate}
 
\end{theorem} 


Notice that if $\Omega_0$ is as above, then  statement (7) in Theorem \ref{t:klingler}  already implies  $\Omega_0 =\C \times \Bbb{H}^+$,  hence up to conjugation  we can assume 
$$\Omega_0=\bigcup_{z\in
\mathbb{H}^+}\overleftrightarrow{e_1,z} \setminus \{e_1\}\; .$$

Let
$G\subset \PSL(3,\mathbb{C})$ be  a  subgroup acting properly  discontinuously on $\Omega_0$ with compact quotient $N=\Omega_0/G$. We set $
\; \ell=\overleftrightarrow{e_3,e_2} \; \; \hbox{and}  \; \,p=e_1.$
 Then:

\begin{lemma}
\label{i:lic1} The point $p$ is $G$-invariant.
\end{lemma}

\begin{proof}
  Let $\ell_1,\,
\ell_2\subset \partial \Omega_0\subset \mathbb{P}^2_{\mathbb{C}}\setminus
\Omega$ be distinct complex lines  and $\gamma\in G$. Then
$\{e_1\}=\ell_1\cap \ell_2 $ and  $\gamma(\ell_1),\,\gamma(\ell_2)\subset
\mathbb{P}^2_{\mathbb{C}}\setminus  \Omega$ are different complex lines which satisfy $\gamma(\ell_i)\cap \bigcup_{z\in
\mathbb{H}^+}\overleftrightarrow{e_1,z}\neq \emptyset $, for $i=1,2$. Thus  $e_1\in \gamma(\ell_1\cap \ell_2)$ and therefore  $\gamma(e_1)=e_1$.
\end{proof}

 This lemma implies that  we have a projection $\pi_{p,\ell}: \PC^2 \setminus \{p\} \to \ell$ and a homomorphism 
  $\Pi_{p,\ell}:G\rightarrow Bihol(\ell)$  as in Lemma
\ref{l:control}. Again, for simplicity we set $\pi:= \pi_{p,\ell}$ and $\Pi:= \Pi_{p,\ell}$.

\vv

 Let  $(\mathcal{D},\mathcal{H})$ be a developing pair  for $N$. Then we know from statement (7) in Theorem \ref{t:klingler}
that $\mathcal{D}(\Omega_0)=\Omega_0$,  $\mathcal{H}(\pi_1(N))=G$ and $N$ is an Inoue surface. Furthermore, we know from Theorem \ref{t:klingler}  that up to conjugation the group $G$ is contained in $ Sol^4_{0}, Sol^4_{1}$ or $Sol^{\prime 4}_0$.
Thus we have:

\begin{lemma}\label{r:inoue} 
Let $\Omega_0 $  and 
$ G$ be as above. Then  $G$
 is  a 
subgroup of either $ Sol^4_{0}$, $Sol^4_{1}$ or $Sol^{\prime 4}_0$.
\end{lemma}

We set: $ \mathbb{R}= \{[0:z:1]: Im(z)=0 \}$ 
and
$
 \overline{\mathbb{R}}=\mathbb{R}\cup\{\infty\}$ with $ \infty=e_2$, and 
recall that the Greenberg limit set $\Lambda_{\Gr}$ of possibly non-discrete  subgroups of $\PSL(2,\R)$  was defined in Section \ref{s:dejedis}.

\begin{lemma} \label{l:lic} With the hypothesis and notation of Theorem \ref{t: Inoue}, set
 $\Gamma_0 = Isot(\Omega_0,\Gamma)$. Then:

\begin{enumerate}
\item \label{i:lic2} $\Lambda_{\Gr}(\Pi(\Gamma_0))=\overline{\mathbb{R}}$.
\item \label{i:lic3} The point  $\infty$ is $\Pi(\Gamma_0)$-invariant, {\it i.e.}, $\Pi(\Gamma_0)e_2=e_2$.
\item \label{i:lic4} $\Pi(\Gamma_0)$ is non-discrete.
\item \label{i:lic5} $\Lambda_{\Gr}(\Pi(\Gamma))=\overline{\mathbb{R}}$.
\item \label{i:lic6} $\overleftrightarrow{e_1,e_2}$ is $\Gamma$-invariant.

\end{enumerate}
\end{lemma}

\begin{proof}

Let us prove (\ref{i:lic2}).  Observe that $\pi(\Omega_0)=\mathbb{H}^+$, thus $\mathbb{H}^+$ is  $\Pi(\Gamma_0)-$invariant and   therefore $\Pi(\Gamma_0)$ is contained in $ \PSL(2,\mathbb{R})$ (up to conjugation).
Now let  $x\in
\mathbb{R}\subset \ell$, where $\ell$ is as above,  and let $K$ be   a fundamental
region for the action of $\Gamma_0$ on $\Omega_0$. By
 Corollary \ref{l:comp} it follows that   $x$  is a cluster point of $\Gamma_0\overline{K}$.
Thus   $\pi(x)=x$ is a cluster point of
$\Pi(\Gamma_0)(\pi(\overline{K}))$. Therefore there exist two subsequences $(\gamma_m)\subset \Pi(\Gamma_0)$  and $(k_m)\subset \pi(\overline{K})$ such that $\gamma_m(k_m)\xymatrix{ \ar[r]_{n \rightarrow  \infty}&} x$.  By Corollary  \ref{trampa} there exists a subsequence  of  $(\gamma_m)$,  still denoted $(\gamma_m)$, which converges uniformly on compact sets of  $\mathbb{H}^+$ to a function $\gamma$. Moreover, $\gamma$ is  either an element of $\PSL_2(\mathbb{R})$ or  a constant in $\overline{\mathbb{R}}$. Since $\gamma_m(k_m)\xymatrix{ \ar[r]_{n \rightarrow  \infty}&} x$ and $x\in \overline{\mathbb{R}}$, it follows that $\gamma_m$ converges uniformly to $x$ on compact sets of $\mathbb{H}^+$.  By Corollary \ref{c:16}  we get  that $x\in \Lambda_{\Gr}(\Pi(\Gamma_0))$, and  since $\mathbb{H}^+$  is $\Pi(\Gamma_0)$-invariant we have  $\Lambda_{\Gr}(\Pi(\Gamma_0))=\overline{\mathbb{R}}$.

\vv

Let us prove (\ref{i:lic3}). By Theorem \ref{t:klingler}, $\Gamma_0$ contains a normal subgroup $\Gamma_1$ with finite index  such that $\Gamma_1$ leaves invariant  $\overleftrightarrow{e_1,e_2}$. Thus $\Pi(\Gamma_1)$ is a normal subgroup  of $\Pi(\Gamma_0)$ with finite index, which leaves  $\infty$ invariant.  Hence  $\Pi(\Gamma_0)(\infty)$ is finite and  $\{\infty\}$ is in the exceptional set of $\Pi(\Gamma_0)$ (see Definition \ref{d:exceptional}). Since the exceptional set  has at most one element, by Corollary \ref{c:11}, it follows that  $\Pi(\Gamma_0)(\infty)=\infty$.

\vv

The proof of (\ref{i:lic4}) follows easily from  parts (\ref{i:lic2}) and (\ref{i:lic3}) of the present lemma.
\vv

Let us prove (\ref{i:lic5}).
 Since $
\Pi(\Gamma(\Omega))$ is an open set which omits
$\overline{\mathbb{R}}$ and contains  $\mathbb{H}^+$ we conclude that
$\Lambda_{\Gr}(\Pi(\Gamma))\neq \ell$. Moreover, since
$\Lambda_{\Gr}(\Pi(\Gamma_0))\subset \Lambda_{\Gr}(\Pi(\Gamma))$
with $\Lambda_{\Gr}(\Pi(\Gamma_0))$  being a circle in $\ell$,  it follows that
$\Lambda_{\Gr}(\Pi(\Gamma))=\overline{\mathbb{R}}$.

\vv
Finally, let us prove (\ref{i:lic6}). From the previous claim we see that $\Pi(\Gamma_0)$ is a subgroup of $\Pi(\Gamma_0)$ with index at most 2. From claim (3), it follows that $\infty \in Ex(\Pi(\Gamma) )$, which implies that $\Pi(\Gamma)\infty=\infty$. Thus $\overleftrightarrow{e_1,e_2}$ is $\Gamma$-invariant.
\end{proof}

\begin{lemma} \label{l:kernel}  Either
$Ker(\Pi\vert_{\Gamma_0})$ is trivial or for  every element
$\gamma\in Ker(\Pi\vert _{\Gamma_0})$ there exists  $\tau(\gamma)\neq 0$ such that the following is a lift of $\gamma$:
\[
\tilde\gamma=
\left (
\begin{array}{lll}
1 & 0 & \tau(\gamma)\\
0 & 1 & 0\\
0 & 0 & 1\\
\end{array}
\right).
\]

\end{lemma}

\begin{proof} If $Ker(\Pi\vert_{\Gamma_0})$ is non-trivial, then
$Ker(\Pi\vert_{\Gamma_0})$ is infinite. Hence, by Lemma \ref{c:cinf}, there exists an element $\gamma_0\in Ker(\Pi\vert_{\Gamma_0})$ with infinite order.
Let $\widetilde{\gamma}_0$ be a lift of $\gamma_0$. By Lemma \ref{l:puente}, it follows that up to  a projective transformation   $\Gamma_0$  is a  subset of either $Sol^4_0$,
$ Sol_1^4$ or $Sol_1^{\prime\,4}$. Therefore there are $\varepsilon =e^{2\pi i
\theta}$ and $ \beta,\, \iota \in\mathbb{C}$, such that:
\[
\widetilde{\gamma}^n_0=
\left (
\begin{array}{lll}
\varepsilon^n & \beta \sum_{j=0}^{n-1}\varepsilon^j &
\iota\sum_{j=0}^{n-1}\varepsilon^j\\
0 & 1 & 0\\
0 & 0 & 1\\
\end{array}
\right).
\]

{\it Claim} 1. -$\varepsilon=1$-
 Since  $\gamma_0$ has infinite order, it follows  that either  $\varepsilon=1$ or $\theta\in \mathbb{R}\setminus
 \mathbb{Q}$. Assuming that  $\theta\in \mathbb{R}\setminus  \mathbb{Q}$, it follows that $\widetilde\gamma_0$
 is  diagonalizable with unitary eigenvalues. Hence
Theorem  \ref{t:pciclic} ensures that
 $L_0(\gamma_0)\cup L_1(\gamma_0)=\mathbb{P}^2_{\mathbb{C}}$ which is not possible, so $\varepsilon=1$.\\

{\it Claim} 2.- $\beta=0$- Recall  that  $\Gamma_0$ is a subgroup of either  $Sol^4_0$,
$ Sol_1^4$ or $Sol_1^{\prime\,4}$. Observe that if  $  \Gamma_0\subset Sol_0^4 $ then   obviously  $\beta =0$. Let us show the claim in the case $\Gamma_0\subset Sol_1^4 $; the reader should observe that  the  same arguments work in the case $Sol_1^{\prime\,4}$.   By (\ref{i:lic2}) of Lemma \ref{l:lic} one has  that   $\Lambda_{\Gr}(\Pi(\Gamma_0))=\overline{\Bbb{R}}$. Thus by Corollary \ref{l:10} follows that $\Pi(\Gamma_0)$ contains hyperbolic elements. Hence  there exists
$\tau\in \Gamma_0$ with a lift $(\tau_{i,j})\in \SL(3,\Bbb{C})$ such that
$\tau _{22}>1$. To conclude, observe:
\begin{scriptsize}
\[
\tau^n \gamma_0\tau^{-n}=
\left [
\begin{array}{llc}
1 & \tau_{22}^{-n}\beta & \iota-\tau_{23}\tau_{22}^{-1}\beta\sum_{j=0}^{n-1}\tau_{22}^{-j}\\
0 & 1 & 0\\
0 & 0 & 1\\
\end{array}
\right]
\xymatrix{ \ar[r]_{n \rightarrow  \infty}&}
\left [
\begin{array}{lll}
1 & 0 & \iota-\tau_{23}\beta(\tau_{22}-1)^{-1}\\
0 & 1 & 0 \\
0 & 0 & 1\\
\end{array}
\right].
\]
\end{scriptsize}
Since $ \Gamma_0$ is discrete,   we deduce $\beta=0$.
\end{proof}

\begin{lemma} \label{l:r0}
 Let $r \in\Bbb{R} $, $z\in \Bbb{C}$, $(x_n)\subset \Omega_0$ be a convergent sequence  whose limit lies in $\Omega_0$ and  $(\gamma_n)\subset \Gamma_0$ a sequence of distinct elements such that $\gamma_n(x_n)\xymatrix{ \ar[r]_{n \rightarrow  \infty}&}[z:r:1]$. If $(\gamma^{(n)}_{ij})\in Sol^4_0\cup Sol^4_1\cup Sol^{\prime 4}_0$ is a lift of $\gamma_n$, it follows that
$
\gamma^{(n)}_{22}\xymatrix{ \ar[r]_{n \rightarrow  \infty}&} 0$ and $ \gamma^{(n)}_{23}\xymatrix{ \ar[r]_{n \rightarrow  \infty}&} r.
$

\end{lemma}

\begin{proof}
Since $\gamma_n(x_n)\xymatrix{ \ar[r]_{n \rightarrow  \infty}&}[z:r:1]$, it follows that $\Pi(\gamma_n(\pi(x_n)))\xymatrix{ \ar[r]_{n \rightarrow  \infty}&}[0:r:1]$. Hence, since $\Pi( \gamma_n)$ leaves invariant $\pi(\Omega_0)=\Bbb{H}$,  Corollary \ref{trampa} and the   last convergence yield that
$
 \Pi(\gamma_n)\xymatrix{ \ar[r]_{n \rightarrow  \infty}&}  r $  uniformly
on compact sets of   $\pi(\Omega_0).$
Finally, since   $(\gamma^{(n)}_{ij})\in Sol^4_0\cup Sol^4_1\cup Sol^{\prime 4}_0$ is a lift of $\gamma_n$, the previous convergence implies:
\[
\tau^{(n)}_{22}\xymatrix{ \ar[r]_{n \rightarrow  \infty}&} 0\;,\, \hbox{and } \tau^{(n)}_{23}\xymatrix{ \ar[r]_{n \rightarrow  \infty}&} r.
\]
\end{proof}

\begin{proposition}\label{p:sol4}
  If $\Gamma_0\subset  Sol^4_0$, then $\partial\Omega_0 \subset L_0(\Gamma)$.
\end{proposition}

\begin{proof} Let $[z:r:1]\in \partial \Omega_0$, then $z\in\Bbb{C}$, $r\in \Bbb{R}$  and  by Corollary  \ref{l:comp}, there  exist a convergent sequence
$(x_n)\subset \Omega_0$  whose limit lies in $\Omega_0$, and  $(\gamma_n)\subset \Gamma_0$ a sequence of distinct elements such that $\gamma_n(x_n)\xymatrix{ \ar[r]_{n \rightarrow  \infty}&}[z:r:1]$. 

Now let $\widetilde\gamma=(\gamma^{(n)}_{ij})\in Sol^4_0$ be a lift of $\gamma_n$ and  let  $\widetilde{\gamma}_0,\widetilde{\gamma}_1,\, \widetilde{\gamma}_2,\,\widetilde{\gamma}_3$ be the generators  given by Theorem \ref{t:sol40}. Using the  notation from Theorem \ref{t:sol40}, let $i\in\{1,2,3\}$ be such that $b_i\neq 0$. Define
$\tau_{n,k,l,m,\tilde n}=\widetilde \gamma_n (\widetilde{\gamma}_0^k\widetilde{\gamma}_i\widetilde{\gamma}_0^{-k})^l(\widetilde{\gamma}_0^m\widetilde{\gamma}_i\widetilde{\gamma}_0^{-m})^{\widetilde n}$, then:
\[
\tau_{n,k,l,m,\tilde n}=
\left (
\begin{array}{lll}
\gamma_{11}^{(n)} & 0                          & \gamma_{11}^{(n)}b_i(\tilde n\beta^m+l\beta^k)+\gamma_{13}^{(n)}\\
0           & \vert\gamma_{11}^{(n)}\vert^{-2} & \vert\gamma_{11}^{(n)}\vert^{-2} a_i(\tilde n\vert\beta\vert^{-2m}+l\vert\beta\vert^{-2k})+\gamma_{23}^{(n)}\\
0 & 0                  & 1 \\
\end{array}
\right ).
\]
Now the conclusion follows from the following claim:\vv

{\it Claim}. There exist  sequences $(k_n),\,(m_n),(l_n),\,(\tilde{n}_n) \subset \Bbb{N}$ and $(p_n)\subset \Bbb{P}^2_{\Bbb{C}}$  such that $p_n\in Fix(\tau_{n,k_n,l_n,m_n,\widetilde{n}_n }) $ and $p_n\xymatrix{ \ar[r]_{n \rightarrow  \infty}&}[z:r:1]$.  Since $\beta^{n},\beta^{n+1}$ are linearly independent we have
 \begin{equation}\label{e:descomposicion}
b_i^{-1} (z_0-z )=r_n\beta^n+s_n\beta^{n+1},
\end{equation}
 for some $r_n,s_n\in \mathbb{R}$. Taking the inverse of $\widetilde{\gamma}_0$ if necessary, we may assume that  $\vert  \beta\vert<1 $, then    equation \ref{e:descomposicion} yields:
\begin{equation}\label{e:fixos}
b_i([r_n]\beta^n+[s_n]\beta^{n+1})\xymatrix{ \ar[r]_{n \rightarrow
 \infty}&} z_0-z,
\end{equation}
 where   $[x]$   denotes the integer part of $x\in \Bbb{R}$. On the other hand,  by lemma \ref{l:r0}, it follows that there exists a sequence of $(\gamma_n)$, still denoted $(\gamma_n)$, such that:
\begin{equation}\label{e:fix1}
\gamma_{22}^{(n)}a_i([r_n]\vert\beta\vert^{-2n}+[s_n]\vert\beta
\vert^{-2(n+1)})+\gamma_{23}^{(n)}\xymatrix{ \ar[r]_{n \rightarrow
 \infty}&}r.
\end{equation}
Since   $\gamma_n(x_n)\xymatrix{ \ar[r]_{n \rightarrow \infty}&}[z:r:1]$ and $\gamma_{11}^{(n)}\xymatrix{ \ar[r]_{n \rightarrow
 \infty}&}\infty$, it follows that:
\begin{equation}\label{e:fix2}
\frac{\gamma_{13}^{(n)}}{\gamma_{11}^{(n)}}\xymatrix{ \ar[r]_{n \rightarrow
 \infty}&}-z_0.
\end{equation}
Set $k_n=n+1, \, l_n=[s_n],\, m_n=n,\, \tilde{n}_n=[r_n]$. A simple calculation shows that
$$
p_n=
\left[
\frac{
\gamma_{11}^{(n)}b_i(\tilde n\beta^m+l\beta^k)+\gamma_{13}^{(n)}
}
{
1-\gamma_n^{(n)}
}:
\frac{
\vert\gamma_{11}^{(n)}\vert^{-2} a_i(\tilde n\vert\beta\vert^{-2m}+l\vert\beta\vert^{-2k})+\gamma_{23}^{(n)}
}
{
1+\vert\gamma_n^{(n)}\vert^{-2}
}:1
\right ]
$$
is  fixed by $\tau_{n,k_n,l_n,m_n,\widetilde{n}_n} $. From Lemma \ref{l:r0} and equations \ref{e:fix1} and \ref{e:fix2}  follows that $p_n\xymatrix{ \ar[r]_{n \rightarrow  \infty}&}[z:r:1]$, which concludes the proof.
\end{proof}

\begin{corollary} \label{r:equi}
If $\Gamma_0\subset Sol^4_0$, then  $\Eq(\Gamma)=\emptyset$.
\end{corollary}

\begin{proposition}\label{p:sol4-1}
If $\Gamma_0\subset  Sol^4_1$, then $\partial\Omega_0 \subset L_1(\Gamma)$.
\end{proposition}

\begin{proof}Let $[z:r:1]\in \partial \Omega_0$,  then $z\in\Bbb{C}$, $r\in \Bbb{R}$. Moreover,   by  Corollary  \ref{l:comp}, there  exists a  sequence
$([z_n:w_n:1])\subset \Omega_0$  converging to $[z_0:w_0:1]\in\Omega_0$,  and  $(\gamma_n)\subset \Gamma_0$ a sequence of distinct elements such that $\gamma_n(x_n)\xymatrix{ \ar[r]_{n \rightarrow  \infty}&}[z:r:1]$. Let $(\gamma^{(n)}_{ij})\in Sol^4_1$ a lift of $\gamma_n$,  then by Lemma \ref{l:r0}, it follows that $\gamma^{(n)}_{22}\xymatrix{ \ar[r]_{n \rightarrow  \infty}&}0$ and $\gamma^{(n)}_{23}\xymatrix{ \ar[r]_{n \rightarrow  \infty}&}r$. On the other hand, by equation \ref{e:formasol2}, it follows that:
 \[ \gamma_n(k_n)=[\pm
z_n+\gamma^{(n)}_{12}w_n+\gamma^{(n)}_{13}: \gamma^{(n)}_{22}w_n
+\gamma^{(n)}_{23}:1].
\]

{\it Claim} 1. Either $(\gamma^{(n)}_{12})$ or $(\gamma^{(n)}_{13})$ is bounded-  Assume these sequences  are unbounded. Then there exists a subsequence of $(\gamma_{n})$, still denoted $(\gamma_{n})$, such that  $\gamma^{(n)}_{12},\gamma^{(n)}_{13}\xymatrix{ \ar[r]_{n \rightarrow  \infty}&} \infty $ and therefore either  $(\gamma^{(n)}_{12}/\gamma^{(n)}_{13} )$ or $(\gamma^{(n)}_{13}/\gamma^{(n)}_{12})$ is bounded. Assume without loss of generality that  $(\gamma^{(n)}_{12}/\gamma^{(n)}_{13})$ is bounded, then there exist a subsequence of $(\gamma_{n})$, still denoted $(\gamma_{n})$, and $c\in \Bbb{C}$ such that  $\gamma^{(n)}_{12}/\gamma^{(n)}_{13}\xymatrix{ \ar[r]_{n \rightarrow  \infty}&} c$. If   $cw_o+1\neq 0$, then:
\[
\gamma_n([z_n:w_n:1])=\left [\frac{\pm z_n+\gamma_{12}^{(n)}w_n}{\gamma_{13}^{(n)}}+1:\frac{\gamma_{22}^{(n)}w_n+ \gamma_{23}^{(n)}}{\gamma_{13}^{(n)}}:\frac{1}{\gamma_{13}^{(n)}}\right ]\xymatrix{ \ar[r]_{n \rightarrow  \infty}&}[cw_0+1:0:0],
\]
which is a contradiction. Thus    $c\neq 0$ and  $\gamma^{(n)}_{13}/\gamma^{(n)}_{12}\xymatrix{ \ar[r]_{n \rightarrow  \infty}&} -w_0$. Hence   $Im(\gamma^{(n)}_{13}(\gamma^{(n)}_{12})^{-1})<0$ for $n$ large, which  contradicts equations \ref{e:formasol2} in page \pageref{e:formasol2}. Therefore $Im(\gamma^{(n)}_{13}(\gamma^{(n)}_{12})^{-1})=0$ for all $n\in \mathbb{N}$, which proves the claim.\vv

{\it Claim} 2. It is not possible that $(\gamma^{(n)}_{12})$ be bounded  and $(\gamma^{(n)}_{13})$ be  unbounded- Assume that  $( \gamma^{(n)}_{12})$ is bounded and  $\gamma^{(n)}_{13}$ is unbounded, then there exist a subsequence of $(\gamma_n)$, still denoted $(\gamma_n)$, and $\gamma_{12}\in \Bbb{C}$ such that $\gamma^{(n)}_{12}\xymatrix{ \ar[r]_{n \rightarrow  \infty}&} \gamma_{12}$   and $\gamma^{(n)}_{13}\xymatrix{ \ar[r]_{n \rightarrow  \infty}&} \infty $. Then
\[
\gamma_n([z_n:w_n:1])=\left[\frac{\pm z_n+\gamma_{12}^{(n)}w_n}{\gamma_{13}^{(n)}}+1:\frac{\gamma_{22}^{(n)}w_n+ \gamma_{23}^{(n)}}{\gamma_{13}^{(n)}}:\frac{1}{\gamma_{13}^{(n)}}\right ]\xymatrix{ \ar[r]_{n \rightarrow  \infty}&}[1:0:0] \,,
\]
which is a contradiction. \vv

A  similar argument   yields that it 
is not possible that $(\gamma^{(n)}_{12})$ be  unbounded and $(\gamma^{(n)}_{13})$ be  bounded. By  Claim 1, it follows that there  exist a subsequence of $(\gamma_n)$, still denoted $(\gamma_n)$, and  $\gamma_{12}, \gamma_{13}\in \mathbb{C}$ such that  $\gamma^{(n)}_{12}\xymatrix{ \ar[r]_{n \rightarrow  \infty}&}
\gamma_{12}$ and $\gamma^{(n)}_{13}\xymatrix{ \ar[r]_{n
\rightarrow  \infty}&} \gamma_{13}$. Finally,  define $x=[\mp z-\gamma_{13}-i\gamma_{12}:i:1]$, clearly  $x\in \Omega_0$ and
\[
\gamma_n(x)=[z+i(\gamma^{(n)}_{12}-\gamma_{12})+\gamma_{13}^{(n)}-\gamma_{13}:\gamma^{(n)}_{22}i+\gamma^{(n)}_{33}:1]\xymatrix{ \ar[r]_{n \rightarrow  \infty}&} [z:r:1].
\]
Thence  $[z:r:1]\in L_1(\Gamma)$. \end{proof}

\begin{proposition}\label{p:sol4-2}
If $\Gamma_0\subset Sol^{\prime\,4}_1$, then $\partial\Omega_0 \subset L_1(\Gamma)$.
\end{proposition}
\begin{proof}
Let $\gamma_0\neq id$ be a   generator of $\Gamma_0$  such that $\gamma_0\in Ker(\Pi\vert
_{\Gamma_0})$ is non-trivial. By Lemma \ref{l:kernel}, it follows  that $\gamma$ has a lift $\widetilde\gamma\in \SL(3,\Bbb{C})$ given  by:
\[
\tilde\gamma=
\left (
\begin{array}{lll}
1 & 0 & \gamma_{13}\\
0 & 1 & 0\\
0 & 0 & 1\\
\end{array}
\right).
\]

On the other hand by part (\ref{i:lic2}) of Lemma \ref{l:lic} it follows that there exists a sequence  $(\gamma_n)\subset\Gamma_0$ such that  $(\Pi( \gamma_n))$ is a sequence of distinct elements and $\Pi(
\gamma_n)\xymatrix{ \ar[r]_{n \rightarrow  \infty}&} Id$.  If $(\gamma^{(n)}_{ij})\in Sol^{\prime,4}_0$ is a lift of $\gamma_n$, it follows that
$\gamma_{22}^{(n)}\xymatrix{ \ar[r]_{n \rightarrow  \infty}&} 1,$
$\gamma_{23}^{(n)}\xymatrix{ \ar[r]_{n \rightarrow  \infty}&} 0$. By equations \ref{e:formasol2} we have that
$Im(\gamma_{13}^{(n)})\xymatrix{\ar[r]_
{n\rightarrow \infty}&} 0$. Hence we can assume that there exist  $(l_n)_{n\in
\mathbb{N}}\subset \mathbb{Z}$ such that
$\gamma_{13}^{(n)}+l_n\gamma_{13}\xymatrix{ \ar[r]_{n \rightarrow
\infty}&} c\in \mathbb{C}$. Therefore:
\[
\gamma_0^{l_n} \gamma_n[z:0:1]=
[z+\gamma_{13}^{(n)}+l_n\gamma_{13}:\gamma_{23}^{(n)}:1]\xymatrix{\ar[r]_{n
\rightarrow  \infty}&}[z+c:0:1].
\]
Hence $\overleftrightarrow{e_1,e_3}\subset L_1(\Gamma)$. To
conclude, observe that
$\overline{\Gamma\overleftrightarrow{e_1,e_3}}=\bigcup_{r
\in \Bbb{R}}\overleftrightarrow{e_1,r}$.
\end{proof}

\begin{lemma} \label{l:inouep}
If   $\gamma_0\in Isot (\Omega_0, \Gamma)$ has finite order, then:
\begin{enumerate}
\item \label{i:lio1} $\gamma_0\in Ker(\Pi=\Pi_{e_1,\overleftrightarrow{e_2,e_3}})$.
\item \label{i:lio2}$Fix(\gamma_0)=\ell_{\gamma_0}\cup \{e_1\}$, where $\ell_{\gamma_0}$
is a complex line not containing $e_1$.
\item\label{i:lio3}$Ker(\Pi\vert_{\Gamma_0})$ is infinite.
\end{enumerate}
\end{lemma}

\begin{proof} Let us show (\ref{i:lio1}).  By part (\ref{i:lic3})  of Lemma \ref{l:lic}  it follows that $\Pi(
\gamma_0)\infty=\infty$. By part (\ref{i:lic4})  of Lemma \ref{l:lic}, it follows  that $\Pi({ \gamma_0})(\mathbb{H}^+)=\mathbb{H}^+$. Therefore
$\Pi({\gamma_0})=Id$.

Let us show (\ref{i:lio2}). Since $\gamma_0\in Ker(\Pi)$ has finite order, there exist  $\kappa,\sigma\in
\mathbb{C}$,  $\vartheta=e^{2\pi
i\theta}$, with $\theta\in \mathbb{Q}\setminus \mathbb{Z}$, such that 
\[
\widetilde\gamma_0=
\left(
\begin{array}{lll}
\vartheta & \kappa  & \sigma\\
0      & 1  & 0\\
0      & 0  & 1 \\
\end{array}
\right),
\]
is a lift  of $\gamma_0$.
A simple inspection shows  that $Fix(\gamma_0)=\ell_{\gamma_0}\cup\{p\}$
where $\ell_{\gamma_0}$ is a complex line that does not contain $p$.

Let us show (\ref{i:lio3}). For this   assume on the contrary that
$Ker(\Pi\vert_{\Gamma_0})$ is finite and define:
\[
\begin{array}{l}
Fin(\Gamma)=\langle\{\gamma\in Isot(\Omega_0,\Gamma) : o(\gamma)<\infty \}\rangle,\\
Ab(\gamma_0)=\{h\in Isot(\Omega_0,\Gamma): h\gamma_0=\gamma_0h\}.
\end{array}
\]

{\it Claim} 1.- $Fin(\Gamma)$ is finite-. Let $\tau\in Fin(\Gamma)$. Since $\Gamma_0$ is a normal subgroup of $Isot(\Omega_0,\Gamma_0)$ with finite index, it follows that there exists an integer $n_0$ such that $\tau^{n}\in \Gamma_0$. Thus  $\tau^{n}\in \Gamma_0\cap Ker (\Pi)$. Since we are assuming that $Ker(\Pi\vert_{\Gamma_0} )$ is finite and $\Gamma_0$ is torsion free, it follows that $\tau^{n}=Id$. Then by  Lemma \ref{c:cinf} we get that $Fin(\Gamma)$ is finite.\vv

{\it Claim} 2.- $Ab(\gamma_0)$ is a subgroup of $Isot(\Omega_0,\Gamma)$ with finite index-. Otherwise, let
$(\gamma_i)\subset Isot(\Omega_0,\Gamma)$ be a sequence of distinct elements such that  $$Ab(\gamma_0)\gamma_i\neq
Ab(\gamma_0)\gamma_j \textrm{ whenever } i\neq j.$$ Since
$\{\gamma_i^{-1}{\gamma_0}\gamma_i\}\subset Fin(\Gamma)$, we have that or each $i$   there exist
$i,j\in \mathbb{N}$ with $i\neq j$  such that
$\gamma_i^{-1}{\gamma_0}\gamma_i=\gamma_j^{-1}\gamma_0\gamma_j$. Hence
$(\gamma_j\gamma_i^{-1})\gamma_0(\gamma_j\gamma_i^{-1})^{-1})\gamma_0^{-1}=Id$, which is a contradiction.\vv

From the previous claim  follows that $\Pi(Ab(\gamma_0))$ is a subgroup of  $\Pi(Isot(\Omega_0,\Gamma)$ with finite index. Since $\Pi(Isot(\Omega_0,\Gamma))$ is non-discrete, it follows that $\Pi(Ab(\gamma_0))$ is non-discrete. Since   $Fix(\gamma_0)=\ell_{\gamma_0}\cup p$, it follows that
$\ell_{\gamma_0}$ is $Ab(\gamma_0)$-invariant.  Hence  $\ell_{\gamma_0}\subset L_0(\Gamma)\cup L_1(\Gamma)$, which is a contradiction.
\end{proof}

\begin{lemma} \label{l:kfinito}
 If $\Gamma_0\subset Sol^4_0$, then $Ker(\Pi_{\Gamma_0})$ is trivial.
\end{lemma}

\begin{proof} 
If $Ker(\Pi_{\Gamma_0})$  is non-trivial, then it is 
infinite, so  there exists an element $\gamma\in Ker(\Pi_{\Gamma_0})$  with infinite order. Then, by  Lemma
\ref{l:comp} and Corollary \ref{l:10}  we have  that there exists $\tau\in\pi_1(M)$ such that $\pi(\tau)$ is hyperbolic. Let $(\tau_{ij})\in Sol^4_0$ be a lift of $\tau$ and $(\gamma_{ij})\in Sol^4_0$ be a lift of $\gamma$. Then  by equations \ref{e:formasol2} we can assume that $\tau_{22}<1$. Finally   from  and equations \ref{e:formasol2} follows that
\[
(\tau_{ij})^{n}(\gamma_{ij})(\tau_{ij})^{-n}=
\left(
\begin{array}{llc}
1 & 0 & \gamma_{13}\tau_{22}^{n}\\
0 & 1 & 0\\
0 & 0 & 1 \\
\end{array}
\right)\xymatrix{ \ar[r]_{n \rightarrow  \infty}&} Id.
\]
This is a contradiction because $\Gamma$ is discrete.
\end{proof}

\begin{lemma}  \label{c:inoue}
Let $M$ be an Inoue surface and  $g:M\longrightarrow M$  a
$(\mathbb{P}^2_{\mathbb{C}},\PSL(3,\mathbb{C}))$-equivalence,
different from the identity, with  at least one fixed point.  Then
$g$ has infinite order.
\end{lemma}

\begin{proof}  Assume   there exists  $g:M\longrightarrow M$ a
$(\mathbb{P}^2_{\mathbb{C}},\PSL(3,\mathbb{C}) )$-equivalence which
is not  the identity,   with finite order and  with at least one fixed
point $z\in M$. Let $P:\mathbb{C}\times \mathbb{H}^+\longrightarrow M$ be
the  universal covering map and $x\in P^{-1}(z)$.    By the standard
lifting lemma for covering maps,
there exists a homeomorphism $\hat
{g}:\mathbb{C}\times\mathbb{H}^+\longrightarrow \mathbb{C}\times\mathbb{H}^+$ such
that the following diagram commutes:
\begin{equation} \label{d:orin}
 \xymatrix{
(\mathbb{C}\times \mathbb{H}^+,x) \ar[r]^{\hat{g}} \ar[d]_{P}&(\mathbb{C}\times
\mathbb{H}^+,x)\ar[d]^{P}\\
(M,z) \ar[r]_{g}  & (M,z).
}
\end{equation}
Since  $P$ and $ g$ are
$(\mathbb{P}^2_{\mathbb{C}},\PSL(3,\mathbb{C}))$-maps we deduce
that $\hat{g}$ is a
$(\mathbb{P}^2_{\mathbb{C}},\PSL(3,\mathbb{C}))$-map.  Since
$\mathbb{C}\times \mathbb{H}^+$  has the projective structure induced
by the natural inclusion we conclude that  $\hat{g}$ is the
restriction of an element  $\gamma_g\in \PSL(3,\mathbb{C})$; moreover
$\gamma_g$ has a    lift $\widetilde{\gamma}_g$ given by:
\[
\widetilde\gamma_g=
\left (
\begin{array}{lll}
\vartheta & \kappa  & \sigma\\
0      & 1  & 0\\
0      & 0  & 1 \\
\end{array}
\right) \,.
\]
By diagram \ref{d:orin} we can ensure that  $\hat g$ has
finite order and $\pi_1(M)$ is a normal subgroup with finite index
of $\Sigma=\langle\gamma_g,\pi_1(M)\rangle$. We also know by Theorem \ref {t:sol40} that $ \pi_1(M)$ is a subgroup of either $Sol^4_0$,  $Sol^4_1$ or  $
Sol^{\prime\,4}_1$.
By Lemma  \ref{l:inouep},   $Ker(\Pi\vert_{\pi_1(M)})$ is   infinite, 
 and    Lemma \ref{l:kfinito} implies that $ \pi_1(M)$  is not in  $Sol^4_0$.
It follows that  $ \pi_1(M)$ is a subgroup of either
  $Sol^4_1$ or  $
Sol^{\prime\,4}_1$.

\vv

{\it Claim} 1. Let
$\tau\in  \pi_1(M) $ and $(\tau_{ij})\in Sol^4_1\cup Sol^{\prime\,4}_1 $ a lift of $\tau$, then
$$\tau_{12}=\kappa\frac{\tau_{22}-1}{1-\vartheta}.$$
Let us prove the claim. Set  $h=\tau\gamma_g\tau^{-1}\gamma^{-1}_g$, then $h\in Ker(\Pi)$,  by equations \ref{e:formasol2}, it follows that  $h$ has
a  lift $\tilde h \in Sol^4_1\cup Sol^{\prime\,4}_1$ given by:
\[
\tilde h=
\left(
\begin{array}{llc}
1 & -\kappa+\tau_{22}^{-1}(\tau_{12}(1-\vartheta)+\kappa)   &
\tau_{23} \tau _{22}^{-1}(\tau_{12}(\vartheta-1)-\kappa)+\tau_{13}(1-\vartheta)\\
0 & 1 & 0\\
0 & 0 &1 \\
\end{array}
\right).
\]
Since $ \pi_1(M)$ has finite index in $\Sigma$, it follows that  there exists $n\in \mathbb{N}$ such
that $h^n\in Ker(\Pi\vert_{\pi_1(M)})$. Therefore   Lemma \ref{l:kernel} yields 
   $\tau_{12}=\kappa(\tau_{22}-1)/(1-\vartheta)$, 
which proves our claim.\vv

On the other hand, since $\Pi(\pi_1(M))$ is non-discrete,  there exists a sequence
$(\gamma_n)\subset\pi_1(M)$ such that $(\pi(\gamma_n))$ is a sequence of distinct elements and  $\pi(\gamma_n)\xymatrix{
\ar[r]_{n \rightarrow  \infty}&} Id$. Let $(\gamma_{ij}^{(n)})\in Sol^4_1\cup Sol^{\prime\,4}_1 $ be a lift of $\gamma_n$. By equations \ref{e:formasol2}, it follows that  $\gamma_{22}^{(n)}\xymatrix{
\ar[r]_{n \rightarrow  \infty}&} 1$, 
$\gamma_{23}^{(n)}\xymatrix{
\ar[r]_{n \rightarrow  \infty}&} 0$ and  $Im(\gamma^{(n)}_{13})\xymatrix{ \ar[r]_{n \rightarrow  \infty}&} 0$. Let
$\gamma\in Ker(\Pi\vert_{\pi_1(M)})$ be
an element with infinite order and $(\gamma_{ij})\in Sol^4_1\cup Sol^{\prime\,4}_1$  a lift of $\gamma$. From Lemma \ref{l:kernel} follows  that $\gamma_{13}\in \mathbb{R}\setminus \{0\}$, so there exists a sequence
$(l_n)\subset \mathbb{Z}$ such that $( l_n\gamma
_{13}+Re(\gamma^{(n)}_{13}))$ is bounded. Thus  we can choose subsequences of  $(\gamma_n )$ and $(l_n)$, of $(\gamma_n)$ and $(l_n)$ respectively, still denoted $(\gamma_n )$ and $(l_n)$, and $c\in \Bbb{R}$, such that  $l_n\gamma
_{13}+Re(\gamma^{(n)}_{13})\xymatrix{ \ar[r]_{n \rightarrow
\infty}&} c$. Hence 
\[
(\gamma_{ij})^{l_n}(\gamma_{ij}^{(n)})=
\left(
\begin{array}{llc}
1 & \kappa(\gamma^{(n)}_{22}-1)/(1-\vartheta) & \gamma_{13}l_n+\gamma_{13}^{(n)}\\
0 & \gamma^{(n)}_{22}                   & \gamma^{(n)}_{23}\\
0 & 0      &1 \\
\end{array}
\right) \xymatrix{ \ar[r]_{n \rightarrow
\infty}&}
\left(
\begin{array}{llc}
1 & 0  & c\\
0 & 1  & 0\\
0 & 0  & 1\\
\end{array}
\right).
\]
This is a contradiction because $\Gamma$ is discrete.
\end{proof}

Notice that we  have improved (\ref{i:k7}) of  Theorem \ref{t:klingler} as follows:

\begin{corollary}  \label{c:klingler} Let $M$ be a compact
$(\mathbb{P}^2_{\mathbb{C}},\PSL(3,\mathbb{C}))$-orbifold which is covered by an Inoue surface,   then the  singular locus of $M$ is empty.
\end{corollary}

Now we can easily proof Theorem  \ref{t: Inoue}.  Notice that statement (\ref{t: Inoue1}) is already proved. Statement (\ref{t: Inoue2})
is immediate from (\ref{p:sol4}), (\ref{p:sol4-1}) and (\ref{p:sol4-2}). Statement (\ref{t: Inoue3}) follows from (\ref{t: Inoue2}) and Lemma \ref{l:limite}. For statement (\ref{t: Inoue4}) notice first that $\Omega_0/\Gamma_0$ is compact by hypothesis. Then by \ref{c:klingler} the quotient $\Omega_{Kul}/\Gamma$ is an Inoue surface, and 
statement (\ref{t: Inoue2}) Theorem  \ref{t: Inoue} implies this quotient consists of one or two copies of $\mathbb H^+ \times \C$. 
\qed

Putting together the results of this section we arrive to Theorem \ref{t: Inoue}.

\section{Quasi-cocompact Controllable Groups} \label{s:eliptic}

In this section we prove:

\begin{theorem}\label{t: elafsur} 
Let  $\Gamma\subset \PSL(3,\C)$ be  discrete and $\Omega$ an open invariant set in $\PC^2$ where 
 $\Gamma$ acts  properly  discontinuously.  Suppose   $M=\Omega/\Gamma$ has a possibly ramified finite cover $S$ which contains a elliptic affine Inoue surface. Then:

\begin{enumerate}
\item There exists a fixed point $p$ under the action of $\Gamma$;
\item there exists a projective  line  $\ell$  not containing $p$, invariant under the action of $\Gamma$;
\item  there exists a hyperbolic open set  $D$  in $\ell \cong \PC^1$ such that the Kulkarni
discontinuity region is $D\times \C^*$;
\item the set $\Omega$ is contained in the Kulkarni set  $\Omega_{Kul}(\Gamma)$, which  is the largest open invariant set where the action is properly discontinuous;  
\item
the limit set is   $\Lambda_{Kul}=L_{0}(\Gamma)\cup L_1(\Gamma)$. 
\item The group $\Gamma $ is a controllable group
with infinite kernel
and quasi co-compact 
control group $\Sigma\subset PSL(2,\C)$  such
that  $\Omega( \Sigma) = D$;
\item The quotient  $\Omega_{Kul}(\Gamma)$ has at most countably many connected components and each of these is an  affine surface. Furthermore, each 
 compact connected component is an elliptic 
surface;
\item we have $Eq(\Gamma)=\Omega_{Kul}(\Gamma)$. 
\end{enumerate} 
\end{theorem}

Let us set: $p=e_1$,
$\ell=\overleftrightarrow{e_3,e_2}$ and $D$ is a hyperbolic domain in $\ell \cong \PC^1$.
 We set $\Omega_0=\bigcup_{z\in D}\overleftrightarrow{z,e_1}\setminus (D\cup\{e_1\})$ and  consider a group   $\Gamma\subset \PSL(3,\mathbb{C})$ whose action  on the orbit 
$\Gamma\Omega_0$ is properly discontinuous and co-compact. Notice that $\Omega_0 = D \times \C^*$.

\begin{lemma} \label{i:lp1}
The point $e_1$ is $\Gamma$-invariant.
\end{lemma}

\begin{proof}  Let
$\gamma\in \Gamma$ and
$\ell_1,\ell_2\subset \partial \Omega_0$ complex lines such that
$e_1\in \ell_1\cap \ell_2$ and $\gamma(\ell_1)\neq
\overleftrightarrow{e_3,e_2}\neq  \gamma(\ell_2)$. Let $z_1,z_2\in
D$ and define  $\xi_{ij}=\overleftrightarrow{z_i,e_1}\cap\gamma(\ell_j)$. Then $\xi_{11}\neq \xi_{12}$, $\xi_{21}\neq \xi_{22}$ and  $\xi_{i,j}\in \{e_1,z_j\}$, so  $e_1\in \gamma(\ell_1)\cap \gamma(\ell_2)$. Hence    $\gamma(e_1)=e_1$.
\end{proof} 

As before, this lemma implies that we have a projection $\pi: \PC^2 \setminus \{p\} \to \ell$ and a 
morphism $\Pi:\Gamma\rightarrow Bihol(\ell)$. 
Let 
  $\Gamma_0$ be  a  torsion free subgroup  
of $Isot(\Omega_0,\Gamma)$ with finite index.

\begin{lemma} \label{l:lp}
The following properties hold:

\begin{enumerate}
\item \label{i:lp2} $D$ is $\Pi(\Gamma_0)$-invariant.
\item \label{i:lp3} $\overleftrightarrow{e_3,e_2}\subset L_0(\Gamma)\cup L_1(\Gamma)$.
\item \label{i:lp4} $\ell$ is $\Gamma$-invariant.
\end{enumerate}
\end{lemma}

\begin{proof}
 Let us  show (\ref{i:lp2}). Let $\gamma\in \Gamma_0$ and  $x\in D$.  As $\partial (\Omega_0)$ is $\Gamma_0$-invariant, we have  that either $\gamma(x)\in \bigcup_{z\in
\partial D}\overleftrightarrow{z,e_1}$ or $\gamma(x) \in D$. Let us assume that $\gamma(x)\in \bigcup_{z\in
\partial D}\overleftrightarrow{z,e_1}$. Then, from the $\Gamma_0$-invariance of $e_1$ we get  that $\pi(\gamma(x))\in \partial D$. On the other hand,  as   $e_1$ is $\Gamma_0$-invariant, we have:
\[
\gamma(\overleftrightarrow{e_1,x}\setminus \{x,e_1\})=\overleftrightarrow{e_1,\gamma(x)}\setminus\{\gamma(x),e_1\}.
\]
From the $\Gamma_0$-invariance of $\Omega_0$, we see that $\gamma(\overleftrightarrow{e_1,x})\setminus\{x,e_1\})\subset \Omega_0$. Now it follows easily that $\pi(\gamma(x))\in D$,  which is a contradiction. Thence $\gamma(x)\in D$.

Let us show (\ref{i:lp3}).  From  (\ref{i:lp2})  it follows that $\overleftrightarrow{e_3,e_2}$ is $\Gamma_0$-invariant. If $\Pi(\Gamma_0)$ is non-discrete, then each point is an accumulation of its own orbit and we have $\overleftrightarrow{e_3,e_2}\subset L_1(\Gamma)\cup L_0(\Gamma)$. So let us assume that $\Pi(\Gamma_0)$  is discrete. In this case,  by (\ref{i:con3}) of
Theorem \ref{l:control2} and Proposition \ref{c:clasfin}, we deduce that
$Ker(\Pi) $ is infinite. So Lemma  \ref{c:cinf} implies that there exists
$\gamma\in Ker(\Pi)$ with infinite order. Moreover, as
$\gamma(\ell\cup\{e_1\})=\ell\cup\{e_1\}$, we conclude that $\gamma$ has a lift $\tilde \gamma\in \SL(3,\mathbb{C}) $ given by:
\[
\tilde\gamma=\left(
\begin{array}{ccc}
a^{-2}   & 0  & 0\\
0   & a  & 0\\
0   & 0  & a\\
\end{array}
\right) \textrm{ where } \vert a\vert\neq 1.
\]
Thus     $\overleftrightarrow{e_3,e_2}\subset
L_0(\Gamma).$\vv

Finally let us prove (\ref{i:lp4}). Let $\gamma\in \Gamma$ and $z\in D$, so  $w=\gamma(\ell) \cap \overleftrightarrow{z,e_1}\subset \mathbb{P}^2_{\mathbb{C}} \setminus \Omega$. Since $\overleftrightarrow{z,e_1}\cap \mathbb{P}^2_{\mathbb{C}}\setminus \Omega=\{e_1,z\}$ and  $e_1$ is $\Gamma$-invariant, it follows that $w=z$ and therefore $D\subset \gamma(\ell)$. Hence $\ell= \gamma(\ell)$.
\end{proof}

\begin{lemma} \label{l:elerio}
  $\Pi(\Gamma_0)$ is non-elementary.
\end{lemma}

\begin{proof}
If this is not the case, by subsection  (\ref{ss:cont})   follows  that
the region of discontinuity $\Omega_0$ can be extended to a region of discontinuity $\Omega_1$, which is either   $\mathbb{C}^2\setminus\{0\}$, 
$\mathbb{C}^*\times\mathbb{C}^*$ or
$\Bbb{C}^*\times\mathbb{C}$. On the other hand, let $R$ be a fundamental region for the action of $\Gamma_0$ on $\Omega_0$.  Then each point in $\partial \Omega_0$ is a cluster point of $R$. Since $R\subset\Omega_1$, it follows that $\partial\Omega_0\subset \Bbb{P}^2_{\Bbb{C}}\setminus \Omega_1$. Finally, since  $D$ omits  a least 3 points in $\ell$, we have  that $\Bbb{P}^2_{\Bbb{C}}\setminus\Omega_1$  contains 3 complex lines which are concurrent, which is a contradiction.
\end{proof}

Let  $\widetilde M$ be   the universal covering space of
$M=\Omega_0/\Gamma_0$, $\pi_1^{Orb}(M)$    the orbifold fundamental group of $M$
and $(\mathcal{D},\mathcal{H})$  a  developing pair for $M$. By
(\ref{i:k8}) of Theorem  \ref{t:klingler} we have that:
\begin{enumerate}
\item[i)]  $\widetilde M$ is biholomorphic to $\mathbb{C}\times \mathbb{H}^+$; 
\item[ii)] there exist $h_1,h_2:\mathbb{H}^+\longrightarrow \mathbb{C}$   holomorphic maps,
$\mu\in \mathbb{C}^*$ such that  $\mathcal{D}(\tilde
M)(z,w)=[e^{-z\mu}:h_1(w):h_2(w)]$; 
\item[iii)] there exists  a subgroup $\Xi$ of $\pi_1(M)$ with
finite index  whose center  $Zen(\Xi)$ contains a free abelian
subgroup of rank 2 with generators  $c,d$. 

\end{enumerate}
Thus
for every  $\vartheta\in \pi_1(M)$ there exist $\omega_{\vartheta}\in \PSL(2,\mathbb{R})$ and
$h_{\vartheta},\gamma_{\vartheta}:\mathbb{H}^+\longrightarrow \mathbb{C}$
 holomorphic maps such that  $h_{\vartheta}(w)\neq w$ for all
$w\in \mathbb{H}^+$ and
$
\vartheta(z,w)=(h_{\vartheta}(w)z+\gamma_{\vartheta}(w),\omega_{\vartheta}(w))
$. So we can  define $P_2:\pi_1(M)   \longrightarrow \PSL(2,\mathbb{R})$  by
$P_2(\vartheta)=\omega_\vartheta$;  clearly    $P_2$
 is a group morphism.

\begin{lemma} \label{l:keraf}
  $Ker (\Pi\vert_{\Gamma_0})$ is infinite.
\end{lemma}

\begin{proof}
Let us  proceed by contradiction. Without loss of generality we can 
assume that $\Gamma_0$ is torsion free,  then
$Ker(\Pi\vert_{\Gamma_0})$ is trivial. By Corollary  \ref{l:comp} and (\ref{i:con3}) of Theorem
\ref{l:control2} we deduce that $\Pi(\Gamma_0)$  is non-discrete.\vv

\noindent
{\it Claim} 1. $o(\Pi(\mathcal{H}(c))),\,
o(\Pi(\mathcal{H}(d)))<\infty$- Otherwise, since $c,d\in Zen (\Xi)$, it follows that   $Fix(g)=Fix(h)$ for all
$g,h\in \Pi(\mathcal{H}(\Xi))$. Then  $\Pi(\Gamma_0)$ is elementary, which contradicts Lemma \ref{l:elerio}.\vv

\noindent
{\it Claim} 2.-$o(P_2(c)), o(P_2(d))<\infty$-
Define   $\widetilde{\mathcal{H}}:P_2(\Xi)\to
\Pi(\mathcal{H}(\Xi))$ defined by $\widetilde{\mathcal{H}}(\vartheta)=
\Pi(\mathcal{H}(P_2^{-1}(\vartheta)))$
Let
$[0:k_1:k_2]\in \pi(\Omega_0)$,  then    there  exists
 $w\in \mathbb{H}^+$ such that
 $[0:h_1(w):h_2(w)]=[0:k_1:k_2]$.
Hence
\begin{equation} \label{d:tildeh}
\widetilde{\mathcal{H}}(\vartheta)([0:k_1:k_2])=[0:h_1(
\vartheta(w)): h_2(\vartheta(w))] \,.
\end{equation}
By equation \ref{d:tildeh}, it follows that  $\tilde{\mathcal{H}}$ is a  well defined group morphism. If either $o(P_2(c))=\infty$ or $o(P_2(d))=\infty$, then $P_2(\Xi)$ is commutative and in consequence   $\Pi(\mathcal{H}(\Xi))$ is
commutative. Thus  $\Pi(\mathcal{H}(\Xi))$ is elementary. Since $\Pi(\mathcal{H}(\Xi))$ is a  subgroup of $\Pi(\Gamma_0)$ with finite index, it follows that $\Pi(\Gamma_0)$ is elementary. This contradicts Lemma \ref{l:elerio} and concludes the proof of the claim.\vv

Let $l\in \mathbb{N}$ be such that $c^l,d^l\in Ker(\Pi\circ \mathcal{H})\cap
Ker(P_2)$.  Then
\[
\begin{array}{l}
c^l(z,w)=(h_{c^l}(w)z+\gamma_{c^l}(w),w);\\
d^l(z,w)=(h_{d^l}(w)z+\gamma_{d^l}(w),w).
\end{array}
\]
Since $c^l,d^l$ do not have fixed points, we  conclude that $h_{c^l}=h_{d^l}=1$ and   $\mathcal{H}(c^l)
=\mathcal{H}(d^l)=Id$,  so we deduce
$$
[e^{-\mu z}:h_1(w):h_2(w)]=[e^{-\mu (z+\gamma_{c^l}(w))}:h_1(w):h_2(w)]=
[e^{-\mu (z+\gamma_{d^l}(w))}:h_1(w):h_2(w)].
$$
Therefore  $e^{\mu h_d(w)}=e^{\mu h_c(w)}=1$. Then  there exist
$k,\,n\in\mathbb{Z}\setminus \{0\}$ such that: $$h_c=2\pi ik\mu^{-1},\,h_d=2\pi
in\mu^{-1}.$$  Hence $c^{lk}=d^{ln}$ which is a contradiction.
\end{proof}

\begin{lemma}
\begin{enumerate}
\item  \label{p:ile1} $\Pi(\Gamma)$ is discrete.
\item \label{p:ile2}$\Pi(\Gamma_0)$ acts properly    discontinuously on
$\pi(\Omega_0)$.
\item \label{p:ile3} $\pi(\Omega_0)/\Pi(\Gamma_0)$ is a compact orbifold.
\item \label{p:ile4} $\bigcup_{z\in\Lambda(\Pi(\Gamma))}
\overleftrightarrow{z,p} \subset L_1(\Gamma)$.
\item \label{p:ile5}  We have  $Eq(\Gamma)=\Omega_{Kul}(\Gamma)$.
\end{enumerate}
\end{lemma}

\begin{proof}

Let us show (\ref{p:ile1}). If this is not the case, then there exists a sequence
$(\gamma_n)\subset\Gamma$ such that $(\Pi (\gamma_n))$ is a sequence of distinct elements such that  $\Pi(
\gamma_n)\xymatrix{ \ar[r]_{n \rightarrow  \infty}&} Id$ uniformly on $\ell$. In consequence, if $(\gamma_{i,j}^{(n)})\in \SL(3,\Bbb{C})$ is a lift of $\gamma_n$ we have:
\[
\sqrt{\gamma_{11}^{(n)}}\gamma_{22}^{(n)},\sqrt{\gamma_{11}^{(n)}}\gamma_{33}^{(n)}
\xymatrix{ \ar[r]_{n \rightarrow  \infty}&} 1
\]
\[
1/\gamma_{11}^{(n)},\sqrt{\gamma_{11}^{(n)}}\gamma_{23}^{(n)},\sqrt{\gamma_{11}^{(n)}}\gamma_{32}^{(n)}
\xymatrix{ \ar[r]_{n \rightarrow  \infty}&} 0.
\]
 Now let   $\gamma_0\in
Ker(\Gamma)$ be an element with infinite order and $(\gamma_{ij})\in \SL(3,\Bbb{C})$  a lift of $\gamma_0$. Taking the inverse of $\gamma_0$ if necessary, we can  assume that $\vert \gamma_{11}\vert<1$. Therefore  there exists a sequence   $(l_n)\subset \mathbb{Z}$ such that
the sequence $(\gamma_{11}^{2l_n} \gamma_{11}^{(n)})$ is bounded and bounded away from 0. Thus, taking  subsequences if necessary, we may assume that there exists $h\in \mathbb{C}^*$ such that $\gamma_{11}^{2l_n} \gamma_{11}^{(n)} \xymatrix{ \ar[r]_{n
\rightarrow  \infty}&} h^{2}.$
In consequence:
\begin{scriptsize}
\[
(\gamma_{ij})^{-l_n}(\gamma_{ij}^{(n)})=\left(
\begin{array}{lll}
\gamma_{11}^{2l_n}\gamma_{11}^{(n)} & 0 & 0\\
0 &\gamma_{11}^{-l_n}\gamma_{22}^{(n)} & \gamma_{11}^{-l_n}\gamma_{23}^{(n)} \\
0 &\gamma_{11}^{-l_n}\gamma_{32}^{(n)} & \gamma_{11}^{-l_n}\gamma_{33}^{(n)}\\
\end{array}
\right)
\xymatrix{ \ar[r]_{n \rightarrow  \infty}&}
\left(
\begin{array}{lll}
h^2 & 0      & 0\\
0   & h^{-1} & 0 \\
0   & 0      & h^{-1} \\
\end{array}
\right).
\]
\end{scriptsize}
This  is a contradiction because $\Gamma$ is discrete.\vv

Let us show (\ref{p:ile2}). Since  $\partial(\pi(\Omega_0))$ is closed and
$\Pi(\Gamma_0)$-invariant, we have
$\Lambda(\Pi(\Gamma_0))\subset\partial\pi(\Omega_0)$. Hence
$\pi(\Omega_0)\subset \Omega(\Pi(\Gamma_0))$.\vv

For  (\ref{p:ile3}), let $R\subset\Omega_0$ be a fundamental domain for the action of
$\Gamma_0$ on $\Omega_0$, then
$\overline{\pi(R)}=\pi(\overline{R})\subset\pi(\Omega_0)$ is
compact. Now the assertion follows easily.\vv

Let us show (\ref{p:ile4}). After conjugating with a projective transformation,  we may assume that
$[0:1:1]\in \Lambda(\Pi(\Gamma_0))$. Thus there exists a sequence  $(\gamma_n)\subset \Gamma_0$    such  that $(\Pi(\gamma_n))$ is a sequence of distinct elements
and  $\Pi (\gamma_n)\xymatrix{ \ar[r]_{n \rightarrow  \infty}&} [0:1:1]$ uniformly on compact sets of $\pi(\Omega_0)$. Now let  $z\in \Bbb{C}^*$,
$[0:z_0:1]\in\pi(\Omega_0)$,    $\gamma_0\in
Ker(\Gamma)$ an element with infinite order  and $(\gamma_{ij})\in \SL(3,\Bbb{C})$  a lift of $\gamma_0$.  Taking the inverse of $\gamma_0$ if necessary, we can as assume that $\vert \gamma_{11}\vert>1$. Now, if $(\gamma_{ij}^ {(n)})\in \SL(3,\Bbb{C} )$ is a lift of $\gamma_n$, it follows that  there exists a sequence   $(l_n)\subset \mathbb{Z}$ such that the sequence
$(
(\gamma_{11}^{(n)}\gamma_{11}^{3l_n})/(\gamma_{32}^{(n)}z_0+\gamma_{33}^{(n)}))$ is bounded and bounded away from $0$. Thus, taking subsequences if necessary, we may assume that there exists $c\in \Bbb{C}^*$ such that $(\gamma_{11}^{(n)}\gamma_{11}^{3l_n})/(\gamma_{32}^{(n)}z_0+\gamma_{33}^{(n)}))\xymatrix{ \ar[r]_{n \rightarrow  \infty}&} c$. 
A straightforward calculation shows:
 $$\gamma^{-l(n)}\gamma_n [zc^{-1}:z_0:1]=\left[\frac{\gamma_{11}^{(n)}\gamma_{11}^{3l_n}zc^{-1}}{\gamma_{32}^{(n)}z_0+\gamma_{33}^{(n)}}:\frac{\gamma_{22}^{(n)}z_0+\gamma_{23}^{(n)}}{\gamma_{23}^{(n)}z_0+\gamma_{33}^{(n)}}:1\right ]\xymatrix{ \ar[r]_{n \rightarrow  \infty}&} [z:1:1].$$
Since $[zc^{-1}:z_0:1]\in \Omega_0$,   Lemma \ref{l:limite}  yields $[z:1:1]\in  L_1(\Gamma)$. In consequence $l=\overleftrightarrow{e_1,[0:1:1]}\in L_1(\Gamma)$. To conclude observe that $\overline{\Gamma l}=\bigcup_{z\in\Lambda(\Pi(\Gamma))}
\overleftrightarrow{z,p}$.

To prove  (\ref{p:ile5}) notice  that $\Omega_{Kul}(\Gamma)$ is a $\Gamma$-invariant set  whose complement has 3 lines in general position. Then
\cite[Theorem 3.5]{bcn} implies that $\Omega_{Kul}(\Gamma)\subset Eq(\Gamma)$, and from Lemma \ref{i:pk5}
we get  $Eq(\Gamma)\subset \Omega_{Kul}(\Gamma)$.
\end{proof}





\vskip.4cm

\section{The Elementary Quasi-cocompact Groups} \label{s:elem}

In this section we prove:

\begin{theorem}\label{t: elementary}
Let $\Omega_0$  be $\PC^2$  minus one, two or three projective lines, and let $\Gamma$ be a group acting properly discontinuously on $\Omega_0$  with compact quotient. Then $\Bbb{P}^2_{\Bbb{C}}\setminus \Omega_0\subset L_0(\Gamma)\cup L_1(\Gamma)$ and we have:
\begin{enumerate}
\item If $\Omega_0  =\mathbb{C}^*\times\mathbb{C}^*$, then $\Gamma$ contains a subgroup $\Gamma_0$ with finite
index, isomorphic to $\mathbb{Z}\oplus \mathbb{Z}$,  where each  element has a lift to ${\rm SL}(3,\C)$ which is a
diagonal matrix.
\item If $\Omega_0=\mathbb{C}\times  \mathbb{C}^*$, then $\Gamma$ has a
 subgroup $\Gamma_0\subset \Gamma$
with finite index, isomorphic to $\mathbb{Z}\oplus \mathbb{Z}\oplus\mathbb{Z} $.
\end{enumerate}
\end{theorem}

This theorem follows from lemmas (\ref{c**}),  (\ref{c*}) and  (\ref{c}) below.

\begin{lemma} \label{c**} 
Let $\Gamma$ be as in Theorem \ref{t: elementary} with  $\Omega_0=\mathbb{C}^*\times\mathbb{C}^*$, then:

\begin{enumerate}
\item \label{i:1l**}  $\Bbb{P}^2_{\Bbb{C}}\setminus \Omega_0\subset L_0(\Gamma)\cup L_1(\Gamma)$.
\item \label{i:2l**} $\Gamma$ contains a subgroup $\Gamma_0$ with finite
index, isomorphic to $\mathbb{Z}\oplus \mathbb{Z}$  where each  element has a lift which is a
diagonal matrix.
\item \label{i:3l**} We have  $Eq(\Gamma) = \Omega_0$.
\end{enumerate}
\end{lemma}

\begin{proof}
Let us prove (\ref{i:1l**}). Notice that $\PC^2 \setminus \Omega_0  =  \overleftrightarrow{e_1,e_2} \cup \overleftrightarrow{e_3,e_1} \cup \overleftrightarrow{e_3,e_2}$. Consider the isotropy subgroups of these lines and 
set
$\Gamma_1=Isot(\overleftrightarrow{e_1,e_2},\Gamma)\cap Isot(\overleftrightarrow
{e_3,e_1},\Gamma)\cap Isot(\overleftrightarrow {e_3,e_2},\Gamma)$. Since
$\mathbb{C}^*\times\mathbb{C}^*$ is $\Gamma$-invariant, we deduce that
$\Gamma_1$ is a normal subgroup of $\Gamma$ and $\Gamma/\Gamma_1$
is a subgroup of $S_3$, the group of permutations in three elements. By Selberg's lemma
there exists a  torsion free normal subgroup $\Gamma_0\subset \Gamma_1$
with finite index. Thus  $\Omega_0/\Gamma_0$ is a compact manifold and $\Gamma_0 e_j=e_j\, (1\leq j\leq 3)$. Let
$\Pi_i=\Pi_{e_i,\overleftrightarrow {e_j,e_k}}$ where $j,k\in
\{1,2,3\}\setminus \{i\}$. By a similar argument as the one used in the proof of Lemma \ref{l:elerio} one can show that either $Ker(\Pi_i)$ is infinite or $\Pi_i(\Gamma_0)$ is
non-discrete. We claim that in both  cases we have  $\overleftrightarrow {e_1,e_2}\cup
\overleftrightarrow {e_3,e_1}\cup\overleftrightarrow
{e_3,e_2}\subset L_0(\Gamma)\cup L_1(\Gamma)$. In fact, if $Ker(\Pi_i)$ is infinite then each of these lines is actually contained in 
$L_0(\Gamma)$ because for each line there is an element in the group with infinite order that leaves invariant each point in the line.
If $\Pi_i(\Gamma_0)$ is
non-discrete then the line is in $L_1(\Gamma)$ because the kernel contains a sequence of points converging to the identity.

Now let us prove (\ref{i:2l**}).  Let $\Gamma_0\subset \Gamma$ be   a torsion free
subgroup with finite index and such that the three lines 
$\overleftrightarrow{e_1,e_2}$, $\overleftrightarrow{e_1,e_3}$ and 
$\overleftrightarrow{e_3,e_2}$  are $\Gamma_0$-invariant. Let
$(\mathcal{D},\mathcal{H})$  be a developing pair for
$M= \Omega_0/\Gamma_0$. By Lemma \ref{l:puente}  we can assume that
$\mathcal{D}(\widetilde M)=\Omega_0$ and $\Gamma_0=\mathcal{H}(\pi_1^{Orb}(M))$. By part (\ref{i:k4}) of Theorem \ref{t:klingler} we can assume
that
$\Omega_0/\Gamma_0$ is a  complex torus and
$\pi_1(\Omega_0/\Gamma_0)=\mathbb{Z}^4$. Since  by Theorem \ref{t:pi_1-orbifolds} the holonomy 
$\mathcal{H}:\pi_1(\Omega_0)\longrightarrow \Gamma_0$  is an
epimorphism, it follows  that $\Gamma_0$ contains  a free abelian
subgroup $\check \Gamma$ of finite index and rank $k\leq 4$. In consequence
$\mathcal{H}^{-1}(\check \Gamma)$ is a free
abelian group of rank 4, see \cite[Section 5.4]{rat}. Let
$(\check{\mathcal{D}},\check{\mathcal{H}})$ be the developing pair
for $\check M=\Omega_0/\check \Gamma$ given by Lemma
\ref{l:puente}. Then by Lemma \ref{c:kerhol} we have  the
following exact sequence of groups:
\begin{equation} \label{d:secgc}
\xymatrix{
0 \ar[r] & \mathbb{Z}^2=\pi_1(\Omega_0)  \ar[r]^{\check{q}_*} &
\mathbb{Z}^4=\pi_1(\check M) \ar[r]^{\check{\mathcal{H}}} & \mathbb{Z}^k=\check\Gamma
\ar[r] & 0\\
}
\end{equation}
where $\check q_*$ is the group morphism induced by the quotient
map $\check q:\Omega_0\longrightarrow \check M $.  Since
\ref{d:secgc} is an exact sequence of free abelian groups, we deduce that
$\mathbb{Z}^4=\mathbb{Z}^2\oplus \mathbb{Z}^k$ (see \cite[Chapter 10]{rotman}) so $k=2$.

Now let us prove (\ref{i:3l**}). Recall that $\Omega_{Kul}(\Gamma)$ is a $\Gamma$-invariant set whose complement has 3 lines in general position.  We know from
\cite[Theorem 3.4]{bcn} that  $\Omega_{Kul}(\Gamma)\subset Eq(\Gamma)$, and from  Lemma \ref{i:pk5}
 we get  $Eq(\Gamma)\subset \Omega_{Kul}(\Gamma)$. 
\end{proof}


\begin{lemma} \label{c*}
Let $\Gamma$ be as in Theorem \ref{t: elementary} with $\Omega_0=\mathbb{C}\times\mathbb{C}^*$, then:
\begin{enumerate}

\item \label{i:2c*} $\Gamma$ has a
subgroup $\Gamma_0\subset \Gamma$
with finite index, isomorphic to $\mathbb{Z}\oplus \mathbb{Z}\oplus\mathbb{Z} $.

\item \label{i:1c*} $\Bbb{P}^2_{\Bbb{C}}\setminus \Omega_0\subset L_0(\Gamma)\cup L_1(\Gamma)$.

\item \label{1:3c*} We have:  $Eq(\Gamma)=\Omega_{Kul}(\Gamma)$.

\end{enumerate}
\end{lemma}

\begin{proof} As in the proof above, we can take a torsion free subgroup $\Gamma_0\subset \Gamma$ with finite index such that $\Gamma_0\subset Isot(\overleftrightarrow{e_1,e_2},\Gamma )\cap Isot(\overleftrightarrow{e_1,e_3},\Gamma )$ and $\Gamma_0$ is isomorphic to $\Bbb{Z}\oplus \Bbb{Z}\oplus \Bbb{Z}$. For simplicity
assume
that $\Gamma_0\subset A_2$, where the latter is the affine group defined in the introduction, in page \pageref{e:formasol3}.

Let
$D_i:\Gamma_0\longrightarrow Bihol(\overleftrightarrow{e_1,e_i})\,,$
$i=2,3$, be defined by
$D_i(\gamma)=\gamma\vert_{\overleftrightarrow{e_1,e_i}}$. Clearly the only interesting case is when
$Ker(D_2)$ and $Ker(D_3)$ are both trivial. In this case $D_2(\Gamma_0)$
and $D_3(\Gamma_0)$ are isomorphic to $\mathbb{Z}^3$. On the other
hand, using the definition of the group $A_2$ we get:
\[
\begin{array}{l}
D_2(\gamma_{ij})(z)=z+g^{-1}_{22}\gamma_{12}\,, \hbox{and} \;
D_3(\gamma_{ij})(z)=\gamma_{11}z.
\end{array}
\]

Therefore $D_{2}(\Gamma_0)$ is isomorphic to an additive
subgroup $\mathbb{C}$ and $D_{3}(\Gamma_0)$ is isomorphic to
a multiplicative subgroup $\mathbb{C}^*$. Now recall that
we know (see for instance \cite{rat}) that the additive subgroups of $\C$ and the multiplicative subgroups $\mathbb{C}^*$
which are discrete and commutative, have rank $\le 2$, which is a contradiction.
Therefore $D_2(\Gamma_0)$ and
$D_3(\Gamma_0)$ are non-discrete groups. Thence we get
$\overleftrightarrow{e_1,e_2}\cup \subset L_0(\Gamma)\cup
L_1(\Gamma)$.

Let us show (\ref{1:3c*}). Since $\Gamma_0$ is contained in either  $A_1$ or $A_2$, given a sequence $(\gamma_n)\subset \Gamma_0$ we get that $\gamma_n$ has a lift $\widetilde \gamma_n\in \SL(3,\C)$ determined by
\[
\widetilde{\gamma_n}=
\left (
\begin{array}{lll}
1 & 0 & b \\
0 & a& 0\\
0 & 0 & 1\\
\end{array}
\right )
\]
\, \textrm{ or }
\[
\widetilde{\gamma_n}=
\left (
\begin{array}{lll}
a & b & 0 \\
0 & a & 0\\
0 & 0 & 1\\
\end{array}
\right ) \;.
\]
Now by Corollary 3.4 in \cite{CS}, it is clear that $Eq(\Gamma_0)=\C\times \C^*$. Finally, since $\Gamma_0$ has finite index in $\Gamma$ we conclude $Eq(\Gamma)=Eq(\Gamma_0)=\C\times \C*$ .

\end{proof}


\begin{lemma} \label{c*}  
Let $\Gamma$ be as in Theorem \ref{t: elementary} with  $\Omega_0=\mathbb{C}\times\mathbb{C}^*$, then: 

\begin{enumerate}
\item \label{i:2c*} $\Gamma$ has a
 subgroup $\Gamma_0\subset \Gamma$
with finite index, isomorphic to $\mathbb{Z}\oplus \mathbb{Z}\oplus\mathbb{Z} $.
\item \label{i:1c*}  $\Bbb{P}^2_{\Bbb{C}}\setminus \Omega_0\subset L_0(\Gamma)\cup L_1(\Gamma)$.

\end{enumerate}

\end{lemma}
\begin{proof} As in the proof above, we can take a torsion free subgroup  $\Gamma_0\subset \Gamma$ with finite index such that  $\Gamma_0\subset Isot(\overleftrightarrow{e_1,e_2},\Gamma )\cap   Isot(\overleftrightarrow{e_1,e_3},\Gamma )$ and $\Gamma_0$ is  isomorphic to $\Bbb{Z}\oplus \Bbb{Z}\oplus \Bbb{Z}$. For simplicity 
assume
that $\Gamma_0\subset A_2$, where the latter is the affine group defined in the introduction, in page \pageref{e:formasol3}.
Let
$D_i:\Gamma_0\longrightarrow Bihol(\overleftrightarrow{e_1,e_i})\,,$
$i=2,3$, be defined by
$D_i(\gamma)=\gamma\vert_{\overleftrightarrow{e_1,e_i}}$. Clearly the only interesting case  is when
$Ker(D_2)$  and  $Ker(D_3)$ are both   trivial. In this case  $D_2(\Gamma_0)$
and $D_3(\Gamma_0)$ are   isomorphic to $\mathbb{Z}^3$. On the other
hand, using the definition of the group $A_2$  we get:
 \[
\begin{array}{l}
D_2(\gamma_{ij})(z)=z+g^{-1}_{22}\gamma_{12}\,, \hbox{and} \;
D_3(\gamma_{ij})(z)=\gamma_{11}z.
\end{array}
\]
 Therefore $D_{2}(\Gamma_0)$  is isomorphic to an  additive
 subgroup  $\mathbb{C}$ and  $D_{3}(\Gamma_0)$ is isomorphic to
 a multiplicative subgroup  $\mathbb{C}^*$.  Now recall that 
we know (see for instance   \cite{rat}) that the additive subgroups of $\C$ and the multiplicative subgroups  $\mathbb{C}^*$
 which are discrete and commutative, have rank $\le 2$, which is a contradiction.
 Therefore $D_2(\Gamma_0)$ and
 $D_3(\Gamma_0)$ are non-discrete groups. Thence we get
 $\overleftrightarrow{e_1,e_2}\cup \subset L_0(\Gamma)\cup
 L_1(\Gamma)$.
\end{proof}


\begin{lemma} \label{c}
Let $\Gamma$ be as in Theorem \ref{t: elementary} with   $\Omega_0=\mathbb{C}^2$. Then
 $\Bbb{P}^2_{\Bbb{C}}\setminus \Omega_0\subset L_0(\Gamma)\cup L_1(\Gamma)$. Moreover, if $\Gamma$ is not a finite extension of the fundamental group of a primary Kodaira surface,  then $Eq(\Gamma)=\Omega_0$.
\end{lemma}

\begin{proof}  Assume that $\Bbb{P}^2_{\Bbb{C}}\setminus \Omega_0\nsubseteq L_0(\Gamma)\cup L_1(\Gamma)$. 
 By Selberg's lemma, there exists a 
torsion free normal subgroup $\Gamma_1\subset \Gamma$  with finite index.
Since $\overleftrightarrow{e_1,e_2}$ is $\Gamma_1$-invariant, we
deduce that  $D:\Gamma\longrightarrow
Bihol(\overleftrightarrow{e_1,e_2})$ given by
$D(\gamma)=\gamma\vert_{\overleftrightarrow{e_1,e_2}}$ is a group
morphism. Thus   $Ker(D)$ is trivial and
$D(\Gamma)$ is discrete. In consequence   $D$ is an isomorphism
and every element in $\Gamma_1$ is unipotent (see (\ref{i:k6}) of
Theorem \ref{t:klingler}), so  every element in
$\Gamma_1\setminus \{Id\}$ has a lift with the following normal  Jordan
form:
\[
\left (
\begin{array}{lll}
1 & 1 & 0\\
0 & 1 & 1\\
0 & 0 & 1\\
\end{array}
\right ).
\]
By  Theorem \ref{l:invlin} we deduce that $D(\Gamma_1)$ contains only
parabolic elements. Since $D(\Gamma_1)$ is discrete we conclude
that   $D(\Gamma_1)$ is isomorphic to either $\mathbb{Z}$ or
$\mathbb{Z}\oplus \mathbb{Z}$. Since $D$ is an isomorphism,  
Lemma \ref{l:suwa} below implies that  $\Omega_0/\Gamma_1$ is
non-compact, which is a contradiction.  

Now suppose $\Gamma$ is not a finite extension of the fundamental group of a primary Kodaira surface. 
Then Theorem  \ref{t:klingler} implies that $\Gamma$ is a finite extension of the
fundamental  group of a complex torus $\s^1 \times \s^1 \times \s^1 \times \s^1  $. Hence 
$\Gamma$ has a finite index normal subgroup 
 $\Gamma_0$ which is generated by four translations, which are determined by 4 linearly independent vectors, say 
 $w_1=(x_1,y_1), \ldots, w_4=(x_4,y_4)$. If 
  $(\gamma_n)\subset \Gamma_0$ is a sequence of distinct elements, then 
 each $\gamma_n$ has a lift $\widetilde \gamma_n\in {\rm \SL}(3,\C)$ given by:
\[
\widetilde{\gamma_n}=
\left (
\begin{array}{lll}
1 & 0 &  \sum_{j=1}^4 \alpha_{jn}x_j  \\
0 & 1 & \sum_{j=1}^4 \alpha_{jn}y_j\\
0 & 0 & 1\\ \end{array}
\right ) \;,
\]
with  $ \alpha_{jn}\in \Bbb{N}$. Then 
 Corollary 3.4 in \cite{CS} implies  $Eq(\Gamma_0)=\C^2$. Since $\Gamma_0$ has finite index in $\Gamma$ we conclude $Eq(\Gamma)=Eq(\Gamma_0)=\C^2$.
\end{proof}


\begin{lemma}[Proposition 2 in \cite{suwa}] \label{l:suwa}
Let $F$ be a free abelian group  acting on $\mathbb{C}^2$ freely and
properly discontinuously. If the rank  of  $F$ is $\le 3$, then  the quotient space  of $\mathbb{C}^2/F$ cannot be
compact.
\end{lemma}



\vskip.4cm

\section{Proof of the main theorems}\label{s: proofs}

In this section we prove  Theorems \ref{t:main1}-\ref{t:main5}. Let $\Omega$ be an open set in $\Bbb{P}^2_{\Bbb{C}}$ where $\Gamma$ acts properly discontinuously, 
$p:\Omega\longrightarrow \Omega/\Gamma$ is the quotient map and
$M\subset \Omega/\Gamma$ is a compact connected component. Then
$\hat {\Omega}=p^{-1}( M)$ has the form $\Gamma \Omega_0$, where
$\Omega_0$ is a connected component of $\Omega$. Set
$\Gamma_0=Isot(\Omega_0,\Gamma)$. Notice that by Lemma \ref{l:puente} there is
a developing pair  $(\mathcal{D},\mathcal{H})$ such that
$\mathcal{D}(\widetilde M)=\Omega_0$  and $\mathcal{H}(\pi_1^ {Orb}
(M))=\Gamma_0$.\\

\begin{proaf}{\sc of Theorem \ref{t:main1}:}
By Corollary \ref{c:coneccion} one of the
following  cases must occur:

{\it Case} 1.- If $\Omega_0=\mathbb{H}^ 2_{\mathbb{C}}$ then  (\ref{i:co2}) of
Corollary \ref{l:comp} and Theorem \ref{t:pfuch} imply that
$\Gamma$ is complex hyperbolic.\vv

{\it  Case} 2.- If $\Omega_0$ is either $\mathbb{C}^ 2\setminus \{0\}$, $\mathbb{C}^*\times\mathbb{C}^*$, $\mathbb{C}^*\times\mathbb{C}$ or  $\mathbb{C}^2$, then it  is clear  that $\Gamma$ is virtually affine.\vv

{\it  Case} 3.- If $\Omega_0=\bigcup_{z\in
\mathbb{H}}\overleftrightarrow{e_1,z} \setminus \{e_1\}$ then by  
Theorem \ref{t: Inoue} the group is affine.\vv
 
{\it  Case} 4.- If $\Omega_0=\bigcup_{z\in D}\overleftrightarrow{z,e_1}\setminus (D\cup\{e_1\})$, where $\ell$ is a complex line  not  containing  $e_1$ and  $D$ is a hyperbolic domain in $\ell$, then by Theorem  \ref{t: elafsur} we know that $\Gamma$ is affine.
\end{proaf}


\vv
\begin{proaf}{\sc of Theorem \ref{t:main2}:}
 By Corollary \ref{c:coneccion} one of the
following  cases must occur:
\vv

\begin{enumerate}
\item If $\Omega_0=\mathbb{H}^ 2_{\mathbb{C}}$, then by (\ref{i:co2}) of
Corollary \ref{l:comp} and Theorem \ref{t:pfuch} we deduce that $\Omega_0=Eq(\Gamma)=\Omega_{Kul}(\Gamma)$ is the largest open set where $\Gamma$ acts properly discontinuously.\vv

\item  If $\Omega_0=\mathbb{C}^*\times\mathbb{C}^*$, then by Lemma \ref{c**} it holds that $\Omega_0=Eq(\Gamma)=\Omega_{Kul}(\Gamma)$ is the largest open set where $\Gamma$ acts properly discontinuously. \vv
  
\item  If $ \Omega_0=\mathbb{C}^*\times\mathbb{C}$, then by lemma
\ref{c*} it holds that $\Omega_0=Eq(\Gamma)=\Omega_{Kul}(\Gamma)$ is the largest open set where $\Gamma$ acts properly discontinuously. \vv

\item  If $ \Omega_0=\mathbb{C}^2$, then by Lemmas
\ref{c} and \ref{l:suwa}, we deduce  that $\Omega_0$ is the largest open set where $\Gamma$ acts properly discontinuously.  Now by lemma \ref{c} we conclude that $Eq(\Gamma)=\Omega_{Kul}(\Gamma)$ whenever $\Gamma$ is not a finite extension of the fundamental group of a primary Kodaira surface. \vv
  
\item If $\Omega_0=\bigcup_{z\in D}\overleftrightarrow{z,e_1}\setminus (D\cup\{e_1\})$, where $\ell$ is a complex line  not  containing  $e_1$ and  $D$ is a hyperbolic domain in $\ell$,  then  by Theorem  \ref{t: elafsur} we get   that $Eq(\Gamma)=\Omega_{Kul}(\Gamma)$ is the largest open set where $\Gamma$ acts properly discontinuously.  \vv

\item If $\Omega_0=\Bbb{C}\times \Bbb{H}$, then Theorem  \ref{t: Inoue}, yields  that $\Omega_{Kul}(\Gamma)$ is the largest open set where $\Gamma$ acts properly discontinuously.
 \end{enumerate}  
\end{proaf}


\vv

\begin{proaf}{\sc of Theorem \ref{t:main3}:}
By Corollary \ref{c:coneccion} we need to consider the
following  cases:\vv

\begin{enumerate}
\item If $\Omega_0=\mathbb{H}^ 2_{\mathbb{C}}$, then by (\ref{i:co2}) of
Corollary \ref{l:comp} and Theorem \ref{t:pfuch} we deduce that $\Omega_{Kul}(\Gamma)=\Omega_0$.\vv

\item  If $\Omega_0=\mathbb{C}^2\setminus \{0\}$, then $\Omega_{Kul}(\Gamma)$is either $\mathbb{C}^2\setminus \{0\}, \Bbb{C}^*\times \Bbb{C}$ or 
$ \Bbb{C}^*\times \Bbb{C}^*$. \vv

\item  If $\Omega_0=\mathbb{C}^*\times\mathbb{C}^*$, then by Lemma \ref{c**} it holds that $\Omega_{Kul}(\Gamma)=\Omega_0$. \vv
  
\item  If $ \Omega_0=\mathbb{C}^*\times\mathbb{C}$, then by lemma
\ref{c*} it holds that $\Omega_0=\Omega_{Kul}(\Gamma)$. \vv

\item  If $ \Omega_0=\mathbb{C}^2$, then by Lemmas
\ref{c} and \ref{l:suwa}, we deduce  that $\Omega_0=\Omega_{Kul}(\Gamma)$. \vv
  
\item If $\Omega_0=\bigcup_{z\in D}\overleftrightarrow{z,e_1}\setminus (D\cup\{e_1\})$, where $\ell$ is a complex line  not  containing  $e_1$ and  $D$ is a hyperbolic domain in $\ell$.  Then  by Theorem  \ref{t: elafsur} we get   that $\Omega_{Kul}(\Gamma)=\Bbb{C}^*\times W$ where $W$ is discontinuity set of some Kleinian group of M\"obius transformations.\vv
 
\item If $\Omega_0=\Bbb{C}\times \Bbb{H}$,  then Theorem  \ref{t: Inoue} yields  that $\Omega_{Kul}(\Gamma)$  consists of two copies  of $\Omega_0$.
 \end{enumerate}  
\end{proaf}


\vv

\begin{proaf}{\sc of Theorem \ref{t:main4}:}

\begin{enumerate}
\item If $\Omega_{\Kul}(\Gamma)=\Bbb{C}^2$, then by Theorem \ref{t:klingler} we have  that  $S_\Gamma$ is  biholomorphic to a complex torus $\s^1 \times \s^1  \times \s^1 \times \s^1$ or a
primary  Kodaira surface.\vv

\item If $\Omega_{\Kul}(\Gamma)= \Bbb{C}^2\setminus \{0\}$,  then $S_\Gamma$ is  biholomorphic to a complex torus or a primary Hopf surface  by Theorem \ref{t:klingler}. \vv

\item  If $\Omega_{\Kul}(\Gamma) =\Bbb{C}^*\times \Bbb{C}$,  then theorems \ref{t:klingler} and \ref{c**}  imply that $S_\Gamma$ is  biholomorphic  to a complex torus. \vv

\item  If $\Omega_{\Kul}(\Gamma) =\Bbb{C}^*\times \Bbb{C}^*$,  then $S_\Gamma$ is  biholomorphic  to a complex torus.
 \vv

\item  If $\Omega_{\Kul}(\Gamma) = \Bbb{C}^*\times (\Bbb{H}^+\cup \Bbb{H}^-)$,  then Theorems \ref{t:klingler} and \ref{t: Inoue} imply that  $S_\Gamma$ is  either $M$ or $M\sqcup M$ where $M$ is a  Inoue surface and the map $S_\Gamma \twoheadrightarrow \Omega_{\Kul}(\Gamma)/\Gamma$ is a covering with no ramification points.

 \vv

\item  If $\Omega_{\Kul}(\Gamma) =  D\times  \Bbb{C}^*$,  then the results of Section  \ref{s:eliptic} imply that  $S_\Gamma$ has countably many components with at least one of them being compact, and each connected component of $S_\Gamma$ is 
 an  elliptic surface with an affine structure.
 \vv

\item  Finally, if $\Omega_{\Kul}(\Gamma) = \Bbb{H}^2_{\Bbb{C}}$,  then  clearly $S_\Gamma$ is   a compact complex hyperbolic manifold.
\end{enumerate}
\end{proaf}


\vv
\begin{proaf}{\sc of Theorem \ref{t:main5}:}
Let $\Gamma\subset \PSL(3,\mathbb{C})$ be  a quasi-cocompact group. 
\begin{enumerate}
\item If $\Omega_{\Kul}(\Gamma)=\Bbb{C}^2$, then $\Gamma$ is affine and  by Theorem \ref{t:klingler}  it is a  finite extension of an   unipotent group. \vv

\item  If $\Omega_{\Kul}(\Gamma) =\Bbb{C}^*\times \Bbb{C}$, then,  by Theorems \ref{t:klingler} and \ref{c*},    $\Gamma$ has a normal subgroup  of finite index isomorphic to  $\Bbb{Z}\oplus\Bbb{Z}\oplus\Bbb{Z}$ which  is either  conjugated to a subgroup of $A_1$ or $A_2$.

\item  If $\Omega_{Kul}(\Gamma)=\mathbb{C}^*\times\mathbb{C}^*$, then  it follows from Theorems \ref{t:klingler} and \ref{c**} that $\Gamma$ has a subgroup $\check \Gamma$ with finite index,  isomorphic to $\Bbb{Z}\oplus \Bbb{Z}$ where each element has a  diagonal matrix as a lift.   \vv

\item  If $\Omega_{\Kul}(\Gamma) = \Bbb{C}^*\times (\Bbb{H}^+\cup \Bbb{H}^-)$,  then  $\Gamma$ is a finite extension of the fundamental group
 of an Inoue Surface by Theorem \ref{t:klingler}.  \vv

\item  If $\Omega_{\Kul}(\Gamma) =  D\times  \Bbb{C}^*$,  then    from  Section  \ref{s:eliptic}  follows that there is a quasi co-compact group $\Sigma\subset \PSL(2,\Bbb{C})$ such that $\Gamma$ is a controllable group with infinite kernel, control group $\Sigma$ and $\Omega(\Sigma)=D$.  \vv

\item  If $\Omega_{\Kul}(\Gamma) = \Bbb{H}^2_{\Bbb{C}}$,  then clearly $\Gamma$ is contained in $\PU(2,1)$  up to projective conjugation.  \vv
\end{enumerate}

\end{proaf}

\begin{proaf}{\sc of Theorem \ref{t:main6}:}
Let  $\Gamma$ be a virtually  cyclic quasi-cocompact
 group, and 
 let $\Omega$ be an open set in $\Bbb{P}^2_{\Bbb{C}}$ where $\Gamma$ acts properly discontinuously;
$p:\Omega\longrightarrow \Omega/\Gamma$ is the quotient map and
$M\subset \Omega/\Gamma$ is a compact connected component. Then
$\hat {\Omega}=p^{-1}( M)$ has the form $\Gamma \Omega_0$, where
$\Omega_0$ is a connected component of $\Omega$. Set
$\Gamma_0=Isot(\Omega_0,\Gamma)$. Notice that by Lemma \ref{l:puente} there is
a developing pair  $(\mathcal{D},\mathcal{H})$ such that
$\mathcal{D}(\widetilde M)=\Omega_0$  and $\mathcal{H}(\pi_1^ {Orb}
(M))=\Gamma_0$.
 Clearly $\Gamma_0$ is virtually cyclic.  By Theorem \ref{t:klingler} and the description  in (\ref{t:sol40}) of the Inoue surfaces, we have that $\Omega_0$ is either $\mathbb H^2_\Bbb{C}$, $\mathbb{C}^*\times\mathbb{C}^*$, $\C^2\setminus\{0\}, \mathbb{C}\times\mathbb{C}^*$, $\C^2$. 
 
 Since the group is virtually cyclic,   Theorem \ref{t:pfuch} yields that $\Omega_0\neq \mathbb H^2_\Bbb{C} $. 
On the other hand, observe that if  $\Omega_0$ is $ \mathbb{C}\times\mathbb{C}^*$ or $ \mathbb{C}\times\mathbb{C}$, then $\Omega/\Gamma_0$ is covered by a complex torus (trough a finite sheteed and possibly ramified covering map), so Lemma \ref{t:pi_1-orbifolds}    implies  that $\Gamma_0$ can not be virtually cyclic, ruling out these cases. Finally observe that the case $\Omega_0=\C^2$ is not possible by Lemma
\ref{l:suwa}. So we conclude that $\Omega_0=\C^2\setminus \{0\}$ and in consequence $\Omega_0/\Gamma_0$ is a covered by Hopf manifold trough a finite sheeted and possibly ramified covering map. Thus  we get $\Omega_0=\Omega$ and $\Gamma=\Gamma_0$.  The converse statement is clear, so we have proved the first part of Theorem \ref{t:main6}.

Let us show now that under the conditions of Theorem \ref{t:main6}, the equicontinuity set coincides with the Kulkarni region of discontinuity. If the group $\Gamma$  has one of these properties, then 
 Theorm  \ref{t:klingler} and Lemma \ref{l:puente} imply that $\Gamma$ has a finite index subgroup $\Gamma_0$, generated by a loxodromic element which has an invariant line and an invariant point. Since $\Gamma_0$ has finite index in $\Gamma$, the Kulkarni limit set of the latter equals that of the former; and also $Eq(\Omega_0)=Eq(\Omega)$.   Finally, 
 since $\Gamma_0$ is cyclic,  Theorem \ref{t:pciclic} implies  $Eq(\Gamma_0)=\Omega_{Kul}(\Gamma_0)$, and therefore 
 $Eq(\Gamma)=\Omega_{Kul}(\Gamma)$.
 
 Assume now that statement   (\ref{i:h3}) holds. Then  using Theorem \ref{t:klingler},  Lemma \ref{l:puente} and the above arguments, we get that  $\Omega_{Kul}(\Gamma)=\C^2\setminus \{0\}$.   Furthermore, since  $\Gamma$ has a finite index subgroup generated by a  loxodromic element $\gamma$, if   $\Gamma$ does not contain complex homotheties, then   $L_0(\Gamma)\cup L_1(\Gamma)$ is the union of the origin and a finite set of points contained in the line at infinity (thinking of $\PC^2$ as being $\C^2$ compactified by attaching to it the line at infinity). Let 
 $\ell\subset \PC^2$ be a complex projective line that does not meet any of these points. Then by the {\it Lambda Lemma} in \cite{pablo1} and 
 an easy computation, we get  that the set of cluster points of the orbit 
 $\{\gamma^n\ell: n\in \Bbb{N}\} $ consists of exactly two lines. Therefore  $\Lambda_{Kul}(\Omega)$ contains two projective lines, which is a contradiction. Thus 
  $\Gamma$ must contain a  complex   homothety and  statement   (\ref{i:h3}) implies statement  (\ref{i:h4}).

That  statement   (\ref{i:h4}) implies     (\ref{i:h2}) is clear. Let us prove now that  (\ref{i:h2}) implies  (\ref{i:h3}).
We know that with the hypothesis of (\ref{i:h2}), up to conjugation,  $\Gamma$ acts properly discontinuously on $\C^2\setminus \{0\}$.  Thus $\C^2\setminus \{0\}\subset \Omega_{Kul}(\Gamma)$.  Furthermore, since $\Gamma$ acts cocompactly on  $\C^2\setminus \{0\}$,  Corollary \ref{l:comp} implies  $\C^2\setminus \{0\}= \Omega_{Kul}(\Gamma)$ as stated.
\end{proaf}


\section{ An elementary quasi-cocompact group which is not affine.} \label{s:tor}

In this section we construct an example of a quasi-cocompact group which is virtually affine but  is not  affine nor complex hyperbolic.

Let  $a\in \Bbb{C}^*$ and $M_a, \,B \in \SL(3,\mathbb{C})$ be given by:
\[
M_a=\left(
\begin{array}{ccc}
a & 0 & 0 \\
0 & a & 0 \\
0 & 0 & a^{-2}
\end{array}
\right), \;
B=\left(
\begin{array}{ccc}
0 & 0 & 1 \\
1 & 0 & 0 \\
0 & 1 & 0 \\
\end{array}
\right) \,.
\]
Notice that $[[B]]$ has order three and cyclically permutes the three lines $\overleftrightarrow{e_1,e_2}$, $\overleftrightarrow{e_3,e_2}$ and $\overleftrightarrow{e_1,e_3}$:
\begin{equation} \label{e:per}
[[B]](\overleftrightarrow{e_1,e_2})=\overleftrightarrow{e_3,e_2}\,,\,
[[B]](\overleftrightarrow{e_3,e_2})=\overleftrightarrow{e_1,e_3}\,,\, \, \hbox{and} \, 
[[B]](\overleftrightarrow{e_1,e_3})=\overleftrightarrow{e_1,e_2}\,.
\end{equation}

\begin{lemma} \label{l:alg}  Define   $\Gamma_a$ as the  group
generated by $[[M_a]]$ and $[[B]]$. Then 

\begin{enumerate}

\item
\[ \label{e:r}
\Gamma_a=
\begin{scriptsize}
\left \{
\left
[
\left
[
B^k_2
\left(
\begin{array}{ccc}
a^{n_1-2n_2+n_3} & 0 & 0 \\
0 & a^{n_1+n_2-2n_3} & 0 \\
0 & 0 & a^{n_2+n_3-2n_1}
\end{array}
\right)
\right ]
\right ]
: n_i\in \mathbb{Z} \textrm{ and } k=0,\,1\,,2\right\}
\end{scriptsize}
\]
\item Let $ \widetilde {\Gamma}_a=
\langle [[B^2M_aB]],[[BM_aB^2]]\rangle $, then
$\widetilde{\Gamma}_a$ is a torsion free normal  subgroup of $\Gamma_a$ with index 3.
\end{enumerate}
\end{lemma}
\begin{proof} An easy computation  shows:
\begin{equation}\label{e:toro}
M_aB=B
\left(
\begin{array}{ccc}
a & 0 & 0 \\
0 & a^{-2} & 0 \\
0 & 0 & a
\end{array}
\right);\,
M_aB^2=B^2
\left(
\begin{array}{ccc}
a^{-2} & 0 & 0 \\
0 & a & 0 \\
0 & 0 & a
\end{array}
\right).
\end{equation}
By an inductive argument we deduce the result.
\end{proof}

\begin{proposition} The group $\Gamma_a$ has  the following properties:
\begin{enumerate}
 \item \label{e:o1} The Kulkarni limit set  is: $$\Lambda_{\Kul}(\Gamma_a)=\overleftrightarrow{e_1,e_2}\cup
\overleftrightarrow{e_1,e_3}\cup \overleftrightarrow{e_3,e_2}.$$
\item  \label{e:o5}  $\Gamma_a$  is a virtually affine group.
\item \label{e:o2} $\Gamma_a$  is not topologically conjugate  to a complex hyperbolic group.
\item \label{e:o3} $\Gamma_a$  is not topologically conjugate to an affine group.
\item \label{e:o4}The quotient  $\Omega_{\Kul}(\Gamma_a)/\Gamma_a$ is a compact orbifold with non-empty singular locus.
\end{enumerate}

\end{proposition}

\begin{proof}
Let us show (\ref{e:o1}). By equation \ref{e:toro}, it follows  that   $\widetilde{ \Gamma}_a$ acts properly 
 discontinuously and freely on $\Omega_0=\Bbb{P}^2_{\Bbb{C}}\setminus\overleftrightarrow{e_1,e_2}\cup
\overleftrightarrow{e_1,e_3}\cup \overleftrightarrow{e_3,e_2}$. Thus from  equation \ref{e:per}   follows that $ \Gamma_a$ acts properly 
 discontinuously  on $\Omega_0$. On the other hand, from equation  \ref{e:toro}  follows  that $L_0(\Gamma_a)=\overleftrightarrow{e_1,e_2}\cup
\{e_1,e_3\}$. Hence equation (\ref{e:per})  
yields $L_0(\Gamma_a)=\Bbb{P}^2_{\Bbb{C}}\setminus\Omega_0$. Thus we have shown that  for every compact set  $K\subset \Omega_0$ the set of cluster points of $\Gamma_a K$ is contained in $L_0(\Gamma_a)$. Thence from part (\ref{i:pk4}) of Proposition  \ref{p:pkg}  follows that
$\Lambda_{\Kul}(\Gamma_a)=\overleftrightarrow{e_1,e_2}\cup
\overleftrightarrow{e_3,e_2}\cup\overleftrightarrow{e_1,e_3}$.\vv

To show  (\ref{e:o5}) recall that $\Lambda_{\Kul}(\Gamma_a)$  is a $\Gamma_a$-invariant set, and 
 by  (\ref{e:o1}) it  consists of three lines. If $\ell$ is one of these lines and $\gamma \in \Gamma$ is an element with infinite order, then  for $n\in \Bbb{N}$ sufficiently large we must have  $\gamma^n\ell=\ell$. Hence the isotropy $Isot(\ell,\Gamma_a)$ is a non-trivial normal  affine subgroup of $\Gamma_a$. On the other hand $\Gamma_a$ leaves invariant the set of these three lines that form  $\Lambda_{\Kul}(\Gamma_a)$. Hence the action of $\Gamma_a$ restricted to  $\Lambda_{\Kul}(\Gamma_a)$ is given by $\Gamma_a/Isot(\ell,\Gamma_a)$, and it is a    subgroup of the permutations group $S_3$. Therefore $\Gamma_a$ is virtually affine.

\vv

Let us show  (\ref{e:o2}) notice that
 by Theorem \ref{t:pciclic}   we know that
 $\Lambda_{\Kul}([[M_a]])=\overleftrightarrow{e_1,e_2}\cup\{e_3\}$.  Then  by
 Theorem \ref{t:pfuch} we have   that  $[[M_a]]$ cannot be topologically conjugate to an
 element of $\PU(2,1)$.

\vv

To show (\ref{e:o3}) assume  on the contrary  that there exists  a homeomorphism
$\phi:\mathbb{P}^2_{\mathbb{C}}\longrightarrow
\mathbb{P}^2_{\mathbb{C}}$ such that
$\widehat{\Gamma}=\phi^{-1}\Gamma\phi $ is an affine
subgroup of $\PSL(3,\mathbb{C})$. Set  $\widehat {M}_a=
\phi^{-1} [[ M_a]] \phi$, thus $Fix( \widehat {M}_a)=\phi^{-1}(\overleftrightarrow{e_1,e_2})\cup
\phi^{-1}(e_1)$. Then  Theorem \ref{t:pciclic} yields that 
 $\phi^{-1}(\overleftrightarrow{e_1,e_2})$ is a
complex line. In consequence  $\phi^{-1}(e_1),\,
\phi^{-1}(e_2),\,\phi^{-1}(e_3)$ are non-collinear  points which are fixed
by $\widehat{M}_a$. After conjugating  with a projective transformation, if necessary,  we can assume that $\phi(e_1)=e_1,\,
\phi(e_2)=e_2,\,\phi(e_3)=e_3$. Hence   there exists $b\in \Bbb{C}^*$  such that 
\begin{equation} \label{e:m}
\widetilde{M}_a=
\left (
\begin{array}{lll}
b & 0 & 0\\
0 & b & 0\\
0 & 0 & b^{-2}\\
\end{array}
\right)
\end{equation}
is a lift of $\widehat{ M}_a$. Let  $\ell$ be the
invariant line under $\widehat{ \Gamma}_a$, then $\ell$ is invariant
under $ \widehat{M}_a$. By equation \ref{e:m} and Theorem \ref{l:invlin}
it follows that  $\ell\in\{\overleftrightarrow{e_3,p}\vert
p\in\overleftrightarrow{e_1,e_2}\}\cup \{
\overleftrightarrow{e_1,e_2}\}$. Since
 $\phi([[B]](\phi^{-1}(e_1)))=e_{2}$, $\phi([[B]](\phi^{-1}(e_2)))=e_{3}$, $\phi([[B]](\phi^{-1}(e_3)))=e_{1}$, we conclude that
$e_1,e_2,e_3\in \ell$, which is a contradiction.\\

T show (\ref{e:o4}) we observe
that $\Omega_{\Kul}(\Gamma_a)/\Gamma_a=((\mathbb{C}^*\times
\mathbb{C}^*)/\widetilde{\Gamma}_a)/(\Gamma_a/\widetilde{\Gamma}_a)$.  Equation \ref{e:toro} yields that
$(\mathbb{C}^*\times \mathbb{C}^*)/\widetilde{\Gamma}_a$ is a complex torus. Since $\Gamma_a/\widetilde{\Gamma}_a$ is finite, it  is enough to observe that   $[[B]]\in Isot([1:1:1],\, \Gamma_a)$ with   $[1:1:1]\in
\Omega_{\Kul}(\Gamma_a)$.
\end{proof}



\vskip.4cm

\section{ A Kissing Schottky Group.} \label{s:sch} 
Classical Schottky groups appear naturally
in the one  dimensional  case, for instance  from   planar covers  of
compact Riemann surfaces via K\"obe's Retrosection Theorem. In  higher  dimensions we say that  a group
$\Gamma\subset \PSL(n+1,\mathbb{C})$ is  a {\it Schottky group}  if there
exist a natural number $g\geq 2$,  a set of generators
$\gamma_1,\ldots,\gamma_g\in \Gamma$   and pairwise disjoint open
sets $R_1,\ldots,R_g,S_1,\ldots,S_g$, such that each of these open
sets is the interior of its closure, the closures of the $2g$ open
sets are pairwise disjoint and 
$\gamma_j(R_j)=\mathbb{P}^{n}_{\mathbb{C}}\setminus \overline{S_j}$ for $j = 1,..., g$.  These groups are discrete, free in $g$ generators and quasi-cocompact (by \cite{cano}). This type of groups  have  rich dynamics and yield to an interesting
 class of compact complex manifolds (see \cite {nori, sv3,  CNS}). Yet, as shown in \cite{cano}, these groups can only exist in odd dimensions.
 
 More generally, one may consider  {\it kissing-Schottky groups}, which spring from \cite{mumford}. These are defined as above but allowing the closures of the sets sets $R_1,\ldots,R_g,S_1,\ldots,S_g$  to touch each other tangentially, just demanding  
  that the group be discrete and 
$\overline{\bigcup_{j=1}^g R_j\cup S_j}\neq
\mathbb{P}^{n}_{\mathbb{C}}$.

Here  we construct a family of kissing-Schottky groups acting  on $\mathbb{P}^2_{\mathbb{C}}$, which  among other  interesting properties satisfy that they are  not  
 virtually affine  nor complex hyperbolic. 
 These are constructed by the suspension construction introduced in subsection \ref {ss:cont}.

\vv

Consider the following  M\"obius transformations:
\begin{equation} \label{e:gs1}
m_1(z)= \frac{(1+i)z  -i}{i z+ 1-i}\;  ; \; m_2(z)=\frac{(1-i)z -i  }{iz
+ 1+i }\; ;\; m_{3}(z)= \frac{3iz+ 10i }{i z+ 3i}\; .
\end{equation}
A direct computation shows that:
\begin{equation} \label{e:gs2}
\begin{array}{l}
m_1(B_1(1+i))=\mathbb{P}^1_{\mathbb{C}}\setminus  \overline{B_1(1-i)}\,;\\
m_2(B_1(-1+i))=\mathbb{P}^1_{\mathbb{C}}\setminus \overline{B_1(-1-i)}\,;\\
m_3(B_1(-3))=\mathbb{P}^1_{\mathbb{C}}\setminus \overline{B(3)}\,,
\end{array}
\end{equation}
where $B_r(x)$ denotes the Euclidean ball with center $x$ and radius $r$. It is easy to show that the group  $\Gamma_s$ generated by $m_1,\, m_2,\, m_3$ is a kissing-Schottky group. 

Now, for each 
$\varepsilon=(\varepsilon_1,\varepsilon_2,\varepsilon_3)\in \mathbb{C}^*\times
\mathbb{C}^2$ define:

\begin{scriptsize}
\[
M_1=\left(
\begin{array}{ccc}
-1-i & i  & 0\\
-i   & -1+i & 0\\
0   & 0   & 1\\
\end{array}
\right),\,
M_2=\left(
\begin{array}{ccc}
1-i & -i  & 0\\
i   & 1+i & 0\\
0   & 0   & 1\\
\end{array}
\right),\,
M_{\varepsilon}=
\left(
\begin{array}{ccc}
3i\varepsilon_1 & 10i\varepsilon_1 & 0\\
i\varepsilon_1  & 3i\varepsilon_1  & 0\\
\varepsilon_2   & \varepsilon_3    & \varepsilon_1^{-2}\\
\end{array}
\right).
\]
\end{scriptsize}

The proof of the following lemma is a simple exercise and is left to the reader:

\begin{lemma} \label{l:dinm2}
If $P_\varepsilon(\lambda)$ denotes the   characteristic polynomial of
 $M_\varepsilon$, then
\begin{enumerate}
 \item \label{i:di1}  $
P_\varepsilon(\lambda)=-\big(\lambda-\varepsilon_1^{-2}\big)\big(\lambda-i\varepsilon_1(3-\sqrt{10})\big)\big(\lambda-i\varepsilon_1(\sqrt{10}+3)\big).
$
\item \label{i:di2} If 
$\vert\varepsilon_1\vert>(\sqrt{10}+3)^{1/3}$, then  there exist $z_1,z_2\in \Bbb{C}^3\setminus \{0\}$ such that $\{z_1,z_2,e_3\}$ is an ordered basis of eigenvectors for  $M_\varepsilon$ with eigenvalues
$\{i\varepsilon_1(\sqrt{10}+3), i\varepsilon_1(-\sqrt{10}+3), \varepsilon_1^{-2}\}$. 

\item The Kulkarni limit set is   $\Lambda_{\Kul}(\langle [[M_\varepsilon]] \rangle)=\overleftrightarrow{z_2, e_3}\cup \overleftrightarrow{z_2,z_1}$.
\item \label{i:di3}
If we set  $\varepsilon_1=-(3+\sqrt{10})^{1/3}e^{-i\pi(1+4\vartheta)/6}$ and 
\[
k_\varepsilon^\pm=\frac{i(\pm\sqrt{10}
\varepsilon_2+\varepsilon_3)e^{i\pi(1+4\vartheta)/6}}{
(3+\sqrt{10})^{1/3} (3(1-e^{2i\pi\vartheta})-\sqrt{10}(\mp
1-e^{2i\pi\vartheta}))},
\]
then the vectors
$$\beta=\left \{
p_1(\varepsilon)=
\left (
\begin{array}{l}
-\sqrt{10}\\
1\\
k_\varepsilon^-\\
\end{array}
\right), \,
p_2(\varepsilon)=\left(
\begin{array}{l}
\sqrt{10}\\
1\\
k_\varepsilon^+\\
\end{array}
\right ) ,\,
\left (
\begin{array}{l}
0\\
0\\
1\\
\end{array}
\right ) \right \}$$ form an ordered basis of eigenvectors for $M_\varepsilon$ with eigenvalues  $\{\alpha_-,\,\alpha_+,\,e^{2\pi i\vartheta}\alpha_-\}$, where $$\alpha_\pm
=\frac{-i(3\pm\sqrt{10})(3+\sqrt{10})^{1/3}}{e^{i\pi(1+4\vartheta)/6}} \;.$$

\item \label{i:di4}   If we set $\varepsilon_1=-(3+\sqrt{10})^{1/3}e^{-i\pi(1+4\vartheta)/6}$
with $\vartheta\in \mathbb{R}\setminus  \mathbb{Q}$, then 
for every point $x$ in
$\mathbb{P}^2_{\mathbb{C}}\setminus 
(\overleftrightarrow{[p_1(\varepsilon)],[p_2(\varepsilon)]}\cup
\overleftrightarrow{[p_1(\varepsilon)],e_3}\cup
\overleftrightarrow{e_3,[p_2(\varepsilon)]})$  the set of cluster points of
$\{[M_\varepsilon^{-n}](x)\}$ is contained in
$\overleftrightarrow{[p_1(\varepsilon)],e_3}$ and is diffeomorphic to
$\mathbb{S}^1$.
\end{enumerate}
\end{lemma}

Now we have:

\begin{proposition} \label{p:kscho} Let $\Gamma_\varepsilon$ be the group generated by $[[M_1]],\,
[[M_2]],\,[[M_\varepsilon]]$, then:

\begin{enumerate}

\item \label{i:kscho1} The group  $\Gamma_{\varepsilon}$ is a complex
kissing-Schottky group with 3 generators. 
\item \label{i:kscho2}
The Kulkarni discontinuity region $\Omega_{\Kul}(\Gamma_\varepsilon)$  is the largest
open set on which  $\Gamma_{\varepsilon}$ acts properly discontinuously. 
Its complement,   the Kulkarni limit set,  is given by:
$$\Lambda_{\Kul}(\Gamma_\varepsilon)=\bigcup_{p\in\Lambda(\Gamma_s)}
\overleftrightarrow{p,e_3} .$$

\item \label{i:kscho3} The group $\Gamma_{\varepsilon}$ is
not topologically conjugate to a  complex hyperbolic
group.

\item  \label{i:kscho5} The group $\Gamma_{\varepsilon}$ is
not topologically conjugate to an elementary group nor to a virtually affine group.

 \item \label{i:kscho4} If $\vert K_\varepsilon^+\vert +
\vert K_\varepsilon^-\vert \neq 0$ and $\varepsilon_1=-(3+\sqrt{10})^{1/3}e^{-i\pi(1+4\vartheta)/6}$, where
$\vartheta \in\mathbb{R}\setminus  \mathbb{Q}$, then   $\Gamma_{\varepsilon}$ is
not topologically  conjugate  to an    affine   group.

\end{enumerate}
\end{proposition}

\begin{proof} Take $p=e_3,\,
l=\overleftrightarrow{e_2,e_1},\, \Pi=\Pi_{p,l},\, \pi=\pi_{p,l}
$. Then   $\Pi([[M_1]])=m_1, \, \Pi([[M_2]])=m_2,$ and $\Pi([[M_\varepsilon]])=m_3$, where $m_1,m_2,m_3$ are given by equation  \ref{e:gs1}  . \\

Let us prove (\ref{i:kscho1}). Consider the following disjoint family of open sets
\[
\begin{array}{ll}
 R_1=\pi^{-1}(\mathbb{D}+1+i), & S_1=\pi^{-1}(\mathbb{D}+1-i),\\
 R_2=\pi^{-1}(\mathbb{D}-1+i), & S_2=\pi^{-1}(\mathbb{D}-1-i),\\
 R_3=\pi^{-1}(\mathbb{D}-3),   & S_3=\pi^{-1}(\mathbb{D}+3).\\
\end{array}
\]
Clearly equation \ref{e:gs2} yields 
\[
\begin{array}{ll}
[[M_1]](R_1)=\mathbb{P}^2_{\mathbb{C}}\setminus  \overline{S_1}, &
[[M_2]](R_2)=\mathbb{P}^2_{\mathbb{C}}\setminus  \overline{S_2},\\

[[M_{\varepsilon}]](R_3)=\mathbb{P}^2_{\mathbb{C}}\setminus  \overline{S_3}, & \bigcup_{i=1}^3\overline{R_i\cup S_i}\neq \mathbb{P}^2_{\mathbb{C}}.\\\\
\end{array}
\]
Therefore $\Gamma_\varepsilon$ is a Kissing-Schottky group with 3 generators.\\

Now we prove (\ref{i:kscho2}). Since $tr^2(m_2)=4$ and $det(M_2+Id)=8$, it follows that   $M_2$ has the following  Jordan's normal form:
\[
\left (
\begin{array}{lll}
1 & 1 & 0\\
0 & 1 & 0\\
0 & 0 & 1\\
\end{array}
\right ).
\]
Then there exists a  complex line $\ell$  such that $e_3\in \ell=Fix
([[M_2]])$. Therefore $\pi(\ell\setminus \{e_3\})=Fix(m_2)$.  Now  equation
(\ref{e:idline}) yields 
\begin{equation}\label{e:sonia}
 \overline{\Gamma_\varepsilon
\ell}=\bigcup_{q\in\Lambda(\Gamma_s)}\overleftrightarrow{q,e_3}.
\end{equation}
Hence
$\bigcup_{q\in\Lambda(\Gamma_s)}\overleftrightarrow{q,e_3}\subset
L_0(\Gamma_\varepsilon)$. On the other hand, part  (\ref{i:kscho1}) of this lemma yields that  $Ker(\Pi)$ is the  trivial trivial group. By Theorem \ref{l:control2}, it follows that  $\Gamma_\varepsilon$
 acts properly  discontinuously
 on
$\mathbb{P}^2_\mathbb{C}\setminus 
\bigcup_{q\in\Lambda(\Gamma_s)}\overleftrightarrow{q,e_3}$. By part 
(\ref{i:pk4})     of Proposition \ref{p:pkg}, it follows that $\Lambda_{\Kul}(\Gamma_\varepsilon)=\bigcup_{q\in\Lambda(\Gamma_s)}
\overleftrightarrow{q,e_3}$.
 Which concludes the proof.\\

For  (\ref{i:kscho3}), notice that  by Lemma \ref{l:dinm2} and Theorem
\ref{t:pciclic} we have  that 
$\Lambda_{\Kul}([[M_\varepsilon]])=\overleftrightarrow{[p_1(\varepsilon)],e_3}\cup
\{[p_2(\varepsilon)]\}$, where $p_1(\varepsilon),p_2(\varepsilon)$ are defined as in Lemma \ref{l:dinm2}.  Thence   Theorem \ref{t:pfuch}  says that  
$[[M_\varepsilon]]$ cannot be topologically conjugate to an element
in $\PU(2,1)$.\\

Let us show  (\ref{i:kscho4}).  Assume on the contrary that there exists a homeomorphism
$\phi:\mathbb{P}^2_{\mathbb{C}}\longrightarrow \mathbb{P}^2_{\mathbb{C}}$ such
that $\phi^{-1}\Gamma_\varepsilon\phi$ is a subgroup of $\PSL(3,\mathbb{C})$ leaving invariant a complex line $\ell$. \\

{\it Claim} 1.- Let $p_1(\varepsilon),p_2(\varepsilon)$ be  defined as in Lemma \ref{l:dinm2}, then  $\phi( \ell)\subset
\overleftrightarrow{[p_1(\varepsilon)],[p_2(\varepsilon)]}\cup\overleftrightarrow{[p_2(\varepsilon)],e_3}
\cup\overleftrightarrow{[p_1(\varepsilon)],e_3}$.  Let us assume on the contrary, that there exists a point $q\in \phi(\ell) \cap \Bbb{P}^2_{\Bbb{C}}\setminus ( \overleftrightarrow{[p_1(\varepsilon)],[p_2(\varepsilon)]}\cup
\overleftrightarrow{[p_2(\varepsilon)],e_3}\cup
\overleftrightarrow{[p_1(\varepsilon)],e_3})$. By part (\ref{i:di2}) of Lemma
\ref{l:dinm2} the set of cluster points  of
$\{[[M_\varepsilon^{-n}]]( q)\}_{n \in\mathbb{N}}$ is contained in
$\overleftrightarrow{[p_1(\varepsilon)],e_3}$ and is diffeomorphic to
$\mathbb{S}^1$. 
Hence 
$\phi^{-1}(\overleftrightarrow{[p_1(\varepsilon)],e_3})\cap \ell$ contains a set homeomorphic to a circle. On the other hand $\phi^{-1}(\overleftrightarrow{[p_1(\varepsilon)],e_3})\subset Fix (\phi^{-1}[[M_\varepsilon]]\phi)$. Thus by   Theorem \ref{t:pciclic}, it follows  that $\phi^{-1}(\overleftrightarrow{[p_1(\varepsilon)],e_3})$ is a complex line. Since $\phi^{-1}(\overleftrightarrow{[p_1(\varepsilon)],e_3})\cap \ell$ is a circle, it follows that
$\phi( \ell)=\overleftrightarrow{[p_1(\varepsilon)],e_3}$, which is a contradiction.\\

{\it Claim} 2. We have $\phi(\ell)= \overleftrightarrow{[p_1(\varepsilon)],[p_2(\varepsilon)]}$. Since    $\phi(\ell)\setminus 
\{[p_1(\varepsilon)],[p_2(\varepsilon)],e_3\}$  is a connected set, the previous claim yields that
either  $\phi(\ell)$ is $\overleftrightarrow{[p_1(\varepsilon)],[p_2(\varepsilon)]}$ or $
\overleftrightarrow{[p_2(\varepsilon)],e_3}$ or $
\overleftrightarrow{[p_1(\varepsilon)],e_3}$. From the equalities
$$\overline{\Gamma_\varepsilon
\overleftrightarrow{[p_1(\varepsilon)],e_3}}=\overline {\Gamma_\varepsilon
\overleftrightarrow{[p_2(\varepsilon)],e_3}}=\bigcup_{q\in \Lambda(\Gamma_s)}
\overleftrightarrow{q,e_3},$$  we get  $\phi(\ell)= \overleftrightarrow{[p_1(\varepsilon)],[p_2(\varepsilon)]}$.\\

{\it Claim} 3. $\phi(\ell)=\overleftrightarrow{e_1,e_2}$. 
One can  check that  $[1:1:0]$ and $[0:0:1]$
are the unique fixed points of $[[M_1]]$ and its matrix representation   with respect  the ordered
basis $\{(1,1,0),\,(0,1,0),\,(0,0,1)\}$  is:
\begin{equation} \label{e:msho}
\left (
\begin{array}{ccl}
-1 & 1  & 0\\
0  & -1 & 0\\
0  & 0  &1\\
\end{array}
\right ).
\end{equation}
Thus Theorem \ref{l:invlin} yields that
$\overleftrightarrow{e_1,e_2}$
and $\overleftrightarrow{[1:1:0],e_3}$ are the unique invariant
complex lines under  $[[M_1]]$. By the previous claims either $\phi(\ell)=\overleftrightarrow{[1:1:0],e_3}$  or $\phi(\ell)=\overleftrightarrow{e_1,e_2}$. Since
$
\Lambda_{\Kul}(\Gamma_\varepsilon) \subset 
\overline
{
\bigcup_{\gamma\in\Gamma_s}
\overleftrightarrow{\gamma([0:1:0]),e_3}
}=
\overline{\Gamma\overleftrightarrow{[1:1:0],e_3}},
$
we conclude
$\phi(\ell)=\overleftrightarrow{e_1,e_2}$.\\

Now notice that  the previous claims yield 
$k^+_\varepsilon=k^-_\varepsilon=0$, which is a contradiction.\\

Finally, let us prove (\ref{i:kscho5}). Let $\phi:\Bbb{P}^2_{\Bbb{C}}\rightarrow \Bbb{P}^2_{\Bbb{C}}$ be a homeomorphism such that $\phi^{-1}\Gamma_\varepsilon\phi$ is a subgroup of $\PSL(3,\Bbb{C})$. As we pointed in the proof of part (\ref{i:kscho2}), there exists a complex line $\ell$  such that $\ell\subset Fix([[M_2]])$. Then   for each $\varphi \in \Gamma_\varepsilon$ we have that $\phi((\varphi(\ell)))$ is a complex line	.  By Proposition \ref{p:pkg}, it follows that $\bigcup_{\varphi \in \Gamma_\varepsilon}\phi((\varphi(\ell)))$ is contained in the complement of any discontinuity region of $\phi^{-1}\Gamma_\varepsilon\phi$. To conclude the proof observe that by equation \ref{e:sonia}  the set $\{\varphi(\ell):\varphi \in\Gamma_\varepsilon\}$ is infinite, so the group is not elementary.

Now suppose that $\Gamma_\varepsilon$ is virtually  affine.  Then there  exists  $n_0\in \Bbb{N}$ such that the group $G$ generated by  $M_1^{n_0}, M_2^{n_0}, M_\varepsilon^{n_0}$  is affine. Now the same arguments used before to show that
 $\Gamma_\varepsilon$ is not  topologically conjugate to an affine group  show that   $G$  cannot be topologically conjugate to an affine group, which is a contradiction.

\end{proof}

\vskip.4cm

\section{APPENDIX: The Equicontinuity Region for Subgroups of
$\PSL(2,\C)$}\label{a:equicontinuo}

We use the notation and definitions of Section \ref{s:equicontinuo}. We need the following lemmas:

\begin{lemma} \label{l:pfis}
Let $\gamma \in \SO(3)\setminus Rot_\infty$ and  let $\{p_+,\, p_-\}$ be its fixed points, then:
$$p_+\overline{p_-}=-1.$$ 
\end{lemma}
\begin{proof}
Just notice that if $z$ is a fixed point of $\gamma$, then $\gamma(-\frac{1}{\bar z}) = -\frac{1}{\bar z}$.
\end{proof}

\begin{lemma} \label{l:2}
Let $\gamma\in \PSL(2,\Bbb{C})$ be an elliptic element  such that  $Fix(\gamma)= \{1,p\}$ with $|p| \le 1$  and $$\gamma(z)=\frac{az+ b}{cz+ d, }$$ where  
$ad-bc=1$.     
Then:
\begin{enumerate}
\item \label{p:simetria} $a=\bar{d}$ if and only if $p\in
\mathbb{R}$. \item \label{p:pabolicos} If $p\in \mathbb{R}$ then $\vert
a\vert =1$ if and only if $p=0$. \item \label{p:loxos} If $p\in
\mathbb{R}$ then  $\vert a\vert <1$ if and only if $p<0$.
\end{enumerate}
\end{lemma}

\begin{proof}  Since $\gamma$  is elliptic, a simple calculation shows that:
\begin{equation} \label{e:forma}
\gamma(z)=
 \frac{
\frac{p\bar{\lambda}-\lambda}{p-1}z +\frac{p(\lambda-\bar{\lambda})}{p-1}
}
{
\frac{\bar{\lambda}-\lambda}{p-1}z + \frac{p\lambda-\bar{\lambda}}{p-1}
} \;,
\end{equation}
for some $\lambda=e^{\pi i\vartheta}$.\\

Let us show  (\ref{p:simetria}). By equation (\ref{e:forma}) it follows that:
\begin{equation}\label{e:coef}
a=\pm\frac{p\overline{\lambda}-\lambda}{p-1} \;; \, \;
d=\pm\frac{p\lambda-\overline{\lambda}}{p-1}\; \hbox{ and } \,
Re(a-d)=\mp\frac{4Im(\lambda)Im(p)}{\vert p-1 \vert^2}.
\end{equation}
Now the assertion follows.\\

To show (\ref{p:pabolicos}) and (\ref{p:loxos}) notice that by equation \ref{e:coef} the inequality  $\vert a\vert \leq 1$ is equivalent to the inequality $\vert
p\overline{\lambda}-\lambda\vert\leq \vert p-1 \vert  $, which is
equivalent  to $4pIm(\lambda)^2\leq 0$. This  proves those statements.
\end{proof}

Now we have: 

\vv
\noindent
{\bf Proof of Proposition (\ref {p:cr1}):}


We  show  statement (\ref{p:min1}) first, {\it i.e.}, that for each $z\in \mathbb{P}^1_{\mathbb{C}}$ one has
that $Cr(-1)z=\mathbb{P}^1_{\mathbb{C}}$. For this,  let  $z\in \mathbb{P}^1_{\mathbb{C}}$, then  
$i\vert z\vert\in Rot_{\infty}z$. We get:
$$
\tau_{-1}^{-1}Rot_{\infty}\tau_{-1}(i\vert z\vert )= 
i\mathbb{R}\cup {\infty} \quad  \hbox{and} \quad
Rot_\infty(i\mathbb{R}\cup {\infty})=\mathbb{P}^1_{\mathbb{C}},
$$
proving statement  (\ref{p:min1}).

Now we prove the first part of statement  (\ref{p:pgro}) in (\ref {p:cr1}), {\it i.e.},  that 
$Cr(-1)= \SO(3)$.
 Let  $\theta,\vartheta\in \mathbb{R}\setminus \mathbb{Q} $,  $\gamma_1(z)=e^{2\pi i \theta}z$ and
\[
\gamma_2(z)=
\frac{
cos (\pi \vartheta )z + i sin(\pi \vartheta)
}
{
i sin (\pi \vartheta)z + cos (\pi \vartheta)
}.
\]
A straightforward calculation shows that $\overline{\langle \gamma_1,\gamma_2\rangle }=Cr(-1)$. Since $\gamma_1$ and $\gamma_2$ are elements in $\SO(3)$ we have  $\langle\gamma_1,\gamma_2\rangle\subset \SO(3)$ and therefore $\overline{\langle \gamma_1,\gamma_2\rangle }
\subseteq  \SO(3)$ because $ \SO(3)$ is closed.

 Now let $\tau\in \SO(3)\setminus Rot_\infty$ and  let $\{p_+,\, p_-\}$ be its fixed points. By part  (\ref{p:min1}) proved above, there exists an element $\gamma_0\in
Cr(-1)$ such that $\gamma_0(\infty)=p_+$. Since $Fix(\gamma_0
\tau\gamma_0^{-1})=\{\gamma_0(0),\gamma_0(\infty)\}$, it
follows from Lemma \ref{l:pfis} that $\gamma_0(0)=p_-$. Therefore $\gamma_0^{-1}
\tau\gamma_0\in Rot_\infty$. Hence $\SO(3)\subseteq Cr(-1)$.\vv


Let us prove statement (\ref{p:esim1}). We start by constructing a certain path  $\phi$ in $Cr(p)$.  Each point   in $\phi$ is an elliptic M\"obius transformation. To construct this path we use the following auxiliary functions: 
\[
a(x)=\frac{x(p-1)-i(p+1)\sqrt{1-x^2}}{p-1} \quad \hbox{and}\quad
c(x)=\frac{-2i\sqrt{1-x^2}}{p-1} \;.
\]
Now define $\phi:(0,1)\rightarrow Cr(p) $ by:
\[
\phi(x)(z)=
\frac{
\vert a(x)\vert z+\frac{ \vert a(x)\vert p \overline{c(x)}}{a(x)}
}
{
\frac{c(x)a(x)}{\vert a(x)\vert}z + \vert a(x) \vert
} \;.
\]
It is clear that this path is in 
 $Cr(p)$. Notice that if we consider the function
$ \eta:(0,1)\rightarrow \Bbb{C} \,$ 
defined by:
\[ \eta(x)=\frac{i\vert a(x)\vert \sqrt{1-\vert a(x)\vert^2}}{c(x)a(x)} \;, \]  
then a straightforward computations shows 
$\phi(\pm \eta(x)) = \pm \eta(x)$.
This yields to  the   statement in (\ref{p:esim1}) about the fixed points of $\gamma_p$. 

To finish the proof, it is now sufficient to show that  in the path
 $\phi$ there is an element with infinite order.
 For this we observe that if $T$ is  a certain M\"obius transformation, to prove that 
  $T$ has infinite order  it is enough to show that at one of its fixed points,  the derivative of its real part is not a root of unity. 
 We have:
\[
Re(\frac{d\phi(x)}{dz}(\eta(x)))=Re(\frac{d\phi(x)}{dz}(-\eta_x))
=\frac{-8px^2+p^2+6p+1}{(p-1)^2} \,.
\]
For simplicity we let
$ f:(0,1)\rightarrow (-1,1) \,$ be defined by:

\[
f(x)=\frac{-8px^2+p^2+6p+1}{(p-1)^2} \,,
\] 
and consider the set 
$$Ur=\{x\in (0,1):x\pm i\sqrt{1-x^2} \; \textrm{ is a root of the
unity} \}.$$
It is clear that the  set
$Ur$ is countable. Since $f$ is not a constant function, the Mean Value Theorem implies that there exists a point 
$r_0\in (0,1)\setminus 
f^{-1}(Ur)$. This completes the proof of statement (\ref{p:esim1}).

\medskip

Finally, let us complete the proof of statement (\ref{p:pgro}) for $p<0$. By part (\ref{p:esim1}) above,   there exists
$\gamma_p\in Cr(p)$ an element with infinite order and its fixed points satisfy
$\{z_p,-z_p\}$, $z_p\in \C$. Define $\kappa(z)=w_p z$,
then $\kappa^{-1}Rot_\infty \kappa =Rot_\infty$ and also 
$\overline{\langle\kappa^{-1}\gamma_p\kappa\rangle}=\tau_{-1}Rot_\infty \tau_{-1}$.
Then $$\overline{\kappa^{-1} \langle 
Rot_\infty,\kappa^{-1}\gamma_p\kappa\rangle \kappa}= S0(3)\,.$$ 
Since $Cr(p)$ is purelly parabolic, Theorem \ref{t:greenberg} yields that $\overline{\kappa^{-1} \langle 
Rot_\infty,\kappa^{-1}\gamma_p\kappa\rangle \kappa}=Cr(p)$.
\qed

\vv

The next lemma is used for proving Theorem \ref{c:5}. Recall that the cross-ratio of a 4-tuple of distinct points
$z_1$, $z_2$, $z_3$, $z_4$ 
 in the complex line   is given by
$$S(z_1,z_2,z_3,z_4) = \frac{(z_1-z_3)}{(z_1-z_4)} \frac{(z_2-z_3)}{(z_2-z_4)} \,.$$
For the point $(0,\infty,1,p)$ we have $S(0,\infty,1,p) = p$, so we can map $(z_1,z_2,z_3,z_4)$ to $(0,\infty,1,p)$ by an element of $\PSL(2,\C)$.

\begin{lemma} \label{c:3}
Let $\gamma_1,\gamma_2\in \PSL(2,\C)$ be elliptic elements
with infinite order  whose fixed point  sets are given by  $Fix(\gamma_1)=\{z_1,z_2\},\,
Fix(\gamma_2)=\{w_1,w_2\}$. If the cross ratio satisfies
$[z_1:z_2:w_1:w_2]\in \mathbb{C}\setminus (\mathbb{R}_-\cup \{0\})$, then
$\langle\gamma_1,\gamma_2\rangle$ contains a loxodromic element.
\end{lemma}

\begin{proof} Since M\"obius transformations
 preserve the cross ratio we can assume that
 $Fix(\gamma_1)=\{0,\infty \}$ and $Fix(\gamma_2)=\{1,p\}$
 with $p\in\mathbb{C}\setminus (\mathbb{R}_-\cup \{0\})$.  Hence there exist
 $\theta\in \mathbb{R}\setminus \mathbb{Q}$, and  $a,b,c,d\in \mathbb{C}$ such that $ad-bc=1$ and:
$$
\gamma_1(z)=\lambda^2z;
$$
\[
\gamma_2(z)=
\frac{az + b }
{cz + d},
\]
where $\lambda=e^{\pi i \theta}$. Without loss of generality we
can  assume that  $a\neq 0$. Now let  $(n_m)_{m\in
\mathbb{N}}\subset (m)_{m\in\mathbb{N}}$ be a subsequence such that
$\lambda^{n_m}\xymatrix{ \ar[r]_{m \rightarrow  \infty}&}
\frac{\vert a \vert}{a}$, then: $$\gamma_1^{n_m}\circ\gamma_2\xymatrix{
\ar[r]_{m \rightarrow  \infty}&}f(z)=\frac{\vert a\vert z+\frac{\vert
a\vert b}{a}}{\frac{acz}{\vert a \vert}+\frac{ad}{\vert a\vert }},
\textrm{ uniformly on } \mathbb{P}^1_{\mathbb{C}}.$$ If $f$ is
loxodromic the result follows easily, so we assume that $f$
is not loxodromic. By   Lemma \ref{l:2}  and the equality:
$$Tr^2(f)=\vert a\vert^2 \left (1+\frac{ad}{\vert a\vert^ 2}\right)^2\,,$$
we conclude that $r=d\bar{a}^{-1}\in \mathbb{R}-\{1\}$. On the other hand, since   $0\leq Tr^2(\gamma_1)=(a+d)^2$ we conclude $a+d\in
\mathbb{R}$. Hence   $0=Im(1-r)Im(a)$, therefore    $a,d\in \mathbb{R}$. To finish the proof, define $H:\mathbb{S}^1\rightarrow \mathbb{R}$ by
   $H(z)=(a^2-d^2)Im (z)$ and observe that
   $H(\lambda^n)=Im(Tr^2(\gamma_1^n\gamma_2))$, $H(i)=a-d$ and  $\overline{\{e^{i\pi \theta n}:n\in \mathbb{Z}\}}=\mathbb{S}^1$.
\end{proof}

Now we recall the statement of Theorem  \ref{c:5}:

\vv
\noindent
{\bf Theorem \ref{c:5}.}
 {\it  Let $\Gamma\subset \PSL(2,\C)$
be a subgroup, then:

\vv \noindent
(1) $\Eq(\Gamma)=\mathbb{P}^1_{\mathbb{C}}$ if and only
if $\Gamma$ is  either  finite or   conjugate  to a subgroup of   $\SO(3)$ or $Dih_\infty$, where $Dih_\infty$ and $\SO(3)$ are as in examples
\ref{e:dih} and \ref{e:cr} respectively. 

\vv \noindent
(2) 
 $\Eq(\Gamma)$ is $\mathbb{C}$, up to a projective transformation,  if
and only if $\Gamma$ is conjugate to a subgroup $ \Gamma_*$ of
${\rm Epa}(\mathbb{C})$ such that $\overline{\Gamma_*}$ contains a
parabolic element,  where ${\rm Epa}(\mathbb{C})$ is as in example
\ref{e:punto}.

\vv \noindent
(3)  $\Eq(\Gamma) $ is  $\mathbb{C}^*$, up to a
projective transformation,  if and only if $\Gamma$ is
conjugate to a subgroup $ \Gamma_*$ of $M\ddot{o}b(\mathbb{C}^*)$
such that
 $\Gamma_*$ contains a loxodromic element.}

\medskip

 Let us prove (\ref{i:el1}). Let $\Gamma$ be an infinite group. By 
Remark \ref{r:in} we  only need to  consider the following cases:\vv

{\it Case} 1. $o(\gamma)<\infty$ for all $\gamma\in\Gamma\setminus  \{Id\}$.
In this case Selberg's lemma and the classification of
the finite groups of $\PSL(2,\C)$, see \cite{maskit}, imply that $\mathcal{B}=\{\langle A \rangle: A  \textrm{ is a non-empty finite
subset of } \Gamma \}$ is an infinite set where each element is
either a cyclic or a dihedral group. Therefore   $\Gamma$ is conjugate  to a subgroup
of $Dih_\infty$.\vv

 {\it Case} 2. $\Gamma$ contains an element $\gamma_1$ with
infinite order. Assume first that $Fix(\gamma)=Fix(\gamma_1)$ for each element
$\gamma\in \Gamma$ with $o(\gamma_1)=\infty$. Then by lemma \ref{c:3}  we conclude that $\Gamma$ is conjugate  to a
subgroup of $Dih_\infty$.\vv

Now we assume that  there exists an
element $\gamma_2$ with infinite order  and such that
$Fix(\gamma_1)\neq Fix(\gamma_2)$. Then from Remark \ref{r:in},
lemmas \ref{l:2}, \ref{c:3} and  part  (\ref{p:pgro}) of Proposition    \ref{p:cr1}  we see that  up to conjugation, 
 $\overline{\langle\gamma_1,\gamma_2\rangle}=\SO(3)$.
Finally, if $\gamma_3\in \Gamma$ is another element, then  by Remark
\ref{r:in}, the proof of  part  (\ref{p:pgro}) of Proposition    \ref{p:cr1}  and
  (\ref{p:min1})  of   that same proposition, we deduce that there exists  $z\in \mathbb{C} $ such that
$Fix(\gamma_3)=\{z,-\overline{z}^{-1}\}$. And again from the  arguments used in the proof of
  (\ref{p:pgro}) of Proposition \ref{p:cr1}  we
deduce that $\gamma_3\in \SO(3)$.  Hence  $\overline{\Gamma}$
is conjugate to   $\SO(3)$.

\vv

The proofs of (\ref{i:el2}) and (\ref{i:el3})   follow easily from Remark
\ref{r:in},  completing the proof of Theorem \ref{c:5}.
\qed

\vv
\vv

The following corollary is now immediate. This is the second  statement in Theorem \ref{t: eq2}:

\begin{corollary} \label{c:6}
Let $\Gamma\subset \PSL(2,\C)$ be  an infinite closed group, then
$\Gamma$ is purely elliptic if and only if
$\Eq(\Gamma)=\mathbb{P}^1_{\mathbb{C}}$.
\end{corollary}

The next corollary is  the third statement in Theorem \ref{t: eq2}:

\begin{corollary} \label{c:7} Let $\Gamma\subset \PSL(2,\C)$ be a group and
${H}\subset \Gamma$ an infinite normal subgroup, such that $Card(
\Lambda_{\Gr}({H}))= 2,0 $. If  $H$ is not conjugate to a subgroup
of $\SO(3)$, then $\Gamma$ is elementary.
\end{corollary}
\begin{proof}
By the proof of Theorem  \ref{c:5} there exists an
${H}$-invariant set  $\mathcal{P}$ with $Card( \mathcal{P})=2$ and
with the following property: If $\mathcal{R}$ is another finite
${H}$-invariant set, then $\mathcal{R}\subset\mathcal{P}$. Since  $
{H}$ is a normal  subgroup we have  that
${H}g(\mathcal{P})=g(\mathcal{P})$ for all $g\in \Gamma$. This
implies  $\mathcal{P}=g(\mathcal{P})$.
\end{proof}

\begin{lemma} \label{l:8}
Let $\gamma_1,\gamma_2\in \PSL(2,\C)$ be parabolic elements
such that $Fix(\gamma_1)\cap Fix(\gamma_2)=\emptyset$, then
$\langle\gamma_1,\gamma_2\rangle$ contains a loxodromic element.
\end{lemma}
\begin{proof}
 After conjugating with a M\"obius transformation we can assume that $\gamma_1(z)=z+\alpha$ and $$\gamma_2(z)=\frac{z}{\beta z+1},$$ for
some $\alpha,\beta\in\mathbb{C}^*$. Then  
$Tr^2(\gamma_2^m\gamma_1)=(m\alpha\beta+2)^2\xymatrix{ \ar[r]_{m
\rightarrow  \infty}&}\infty$. Hence  $\gamma_2^m\gamma_1$  is
loxodromic for $m$ large.
\end{proof}

Corollaries \ref{l:10} and  \ref{c:tran} below form together the statement (2) in Theorem \ref{t:eq3}:

\begin{corollary} \label{l:10}
Let $\Gamma\subset \PSL(2,\C)$ be a non-elementary subgroup. Then
$\Gamma$ contains a loxodromic element.
\end{corollary}
\begin{proof}
Assume that $\Gamma$ does not contain loxodromic elements, then by
Corollary \ref{c:6} and Remark  \ref{r:in}  we deduce that
$\overline{\Gamma}$ contains a parabolic element $\gamma_0$.
 After conjugating with a M\"obius transformation, we can assume that   $\infty$ is the unique fixed point of $\gamma_0$. By
Lemma \ref{l:8} we conclude that $\overline{\Gamma}\infty=\infty$.
Therefore every element in $\Gamma$ has the form $az+b$ with
$\vert a\vert=1$. That is $\Gamma\subset {\rm Epa}(\mathbb{C})$. This is a
contradiction by (\ref{i:el2}) of Theorem  \ref{c:5}.
\end{proof}

From  Corollary \ref{l:10} and using standard arguments as in  \cite{maskit} we
get  the following corollaries. Notice that \ref{c:11} is statement (3) in Theorem \ref{t:eq3}, while \ref{c:12} and \ref{c:12-2} complete Theorem   \ref{t: eq2}.

\begin{corollary} \label{c:tran}
If $\Gamma\subset \PSL(2,\C)$ is  non-elementary, then
 $\mathbb{P}^1_{\mathbb{C}}\setminus \Eq(\Gamma)$ is the closure of the loxodromic
fixed points.

\end{corollary}

\begin{definition}\label{d:exceptional}
Let $\Gamma\subset \PSL(2,\C)$ be non-elementary.
Define the set of {\it exceptional points} of $\Gamma$ as $$Ex(\Gamma)=\{z\in
\mathbb{P}^1_{\mathbb{C}}\setminus \Eq(\Gamma):\overline{\Gamma z}\neq
\mathbb{P}^1_{\mathbb{C}}\setminus \Eq(\Gamma)\}.$$ 
\end{definition}

\begin{corollary} \label{c:11}
Let $\Gamma\subset \PSL(2,\C)$ be non-elementary. Then $Card( Ex(\Gamma)) < 2$.
\end{corollary}

\begin{proof}
Let  $p$ be an exceptional point for the group. From Corollary \ref{l:10} we know that 
 $\Gamma $ has a loxodromic element. 

We claim that every $p$ is a fixed point of every loxodromic element in  $\Gamma $. To show this, assume on the contrary that there exists a  loxodromic element $\gamma$ in  $\Gamma $ such that $\gamma(p) \ne p$. Then the set 
 $A=\{\gamma^n p:n\in \Bbb{N} \}$  has infinite cardinality.  Let 
  $z$ be a point in $\PC^1\setminus Eq(\Gamma) $ and  $W$ a neighborhood of $z$. Then by Corollary \ref{c:tran} we have that there exists  a  loxodromic element $\gamma_1\in \Gamma$ and  $p_2\in Fix(\gamma_1)\cap W$. On the other hand, since  $A$ is an infinite set we have that there exists   $n_0\in \Bbb{N}$ such that  $\gamma^{n_0}p$ is not fixed for  $\gamma_1$. Since $\gamma_1$ is loxodromic and  $p_1$ is one of its fixed points, we see that there exists $n_1\in \Bbb{N}$ such that  $\gamma_1^{n_1}(\gamma^{n_0}p)\in W$. Therefore  $z$ is a cluster point for the orbit of $p$.  Hence  $\PC^1\setminus Eq(\Gamma)\subset\overline{\Gamma p} $ which is not possible $p$ is exceptional. Thus $p$ is a fixed point of every loxodromic element in  $\Gamma $.
  
  Now we claim that $p$ actually is   a fixed point of every  element in  $\Gamma $. Assume on the contrary that there exists an element $\gamma$ in  $\Gamma $ such that $\gamma(p) \ne p$. Since 
  $\gamma (p)\in Ex(\Gamma)$ we have from the previous claim that  $\gamma (p)$ is a fixed point of all  loxodromic elements. Since $\gamma (p)\neq p $ and the  loxodromic  elements have exactly 2 fixed points, from  Corollary \ref{c:tran} we get that 
  $\Gamma$ is elementary,  which is a contradiction. Thus  $p$ is  $\Gamma$-invariant.

Thus we have that every point in 
 $Ex(\gamma)$ is  $\Gamma$-invariant. Since  M\"obius transformations have at most 2 fixed points, we conclude that  $Ex(\Gamma)$ has at most  2 points. On the other hand, if $Ex(\gamma)$ had 2 points, this would imply that 
  $\Gamma$ is  conjugate to  a  subgroup of  $M\ddot{o}b(C^* )$, so $\Gamma$ is elementary, which is a contradiction, and we arrive to (\ref{c:11}).
\end{proof}


\begin{corollary}\label{c:12} Let $\Gamma\subset \PSL(2,\C)$ be a subgroup
and $\mathcal{C}\neq Ex(\Gamma)$ a closed $\Gamma$-invariant set.
Then $\Lambda_{\Gr}(\Gamma)\subset \mathcal{C}$.
\end{corollary}

\begin{corollary}\label{c:12-2} 
Let $\Gamma\subset \PSL(2,\C)$ be a  group, then
$\Eq(\Gamma)=\mathbb{P}^1_{\mathbb{C}}\setminus \Lambda_{\Gr}(\Gamma)$.
\end{corollary}

The following two lemmas are used in the proof of the proposition below: 

\begin{lemma} \label{l:13}
Let $\Gamma\subset \PSL(2,\C)$ be a purely parabolic closed Lie
group with $dim_{\mathbb{R}}(\Gamma)=1$ and   $\gamma_1,\,\gamma_2\in
\PSL(2,\C)$ be loxodromic elements such that
$\Gamma\infty=\gamma_1(\infty)=\gamma_2(\infty)=\infty$ and
$Fix(\gamma_2)\not\subset \Gamma(Fix(\gamma_1))$. Then
$\Gamma_0=\{\gamma\in \langle\Gamma,\gamma_1,\gamma_2\rangle:Tr^2(\gamma)=4\}$
is a lie group  with $dim_{\mathbb{R}}(\Gamma_0)=2$.
\end{lemma}

\begin{proof} Notice that the set of parabolic elements in $\PSL(2,\C)$ that fix $\infty$ is 2-dimensional. Hence the dimension of $\Gamma_0$ is 1 or  2. 
Assume   that $dim_{\mathbb{R}}(\Gamma_0)=1$.  After conjugating with a M\"obius, we   we can assume that: $\gamma_1(z)=t^2z$ , $\gamma_2(z)= a^2z + ab$ and 
$\Gamma=\{z+r:r\in \mathbb{R}\}.$
Hence $\gamma_2\Gamma\gamma_2^{-1}=\{z+t^2r:r\in \mathbb{R}\},\,
\gamma_1\Gamma\gamma_1^{-1}=\{z+a^{2}r:r \in\mathbb{R}\}$ and in consequence  $a^{2},t^2\in \mathbb{R}$. On the other hand, for all $n\in \mathbb{Z}$ one has:
$\gamma_1^n\gamma_2\gamma_1^{-n}\gamma_2^{-1}(z)= z
+ ab(t^{2n}-1).$ Therefore $\{z+abr:r\in \mathbb{R}\}\subset
\langle\Gamma,\gamma_1,\gamma_2\rangle$ and in consequence $ab\in \mathbb{R}$.
This is a contradiction since $ab(1-a^2)^{-1}\in Fix(\gamma_2)$.
\end{proof}

\begin{lemma} \label{l:14}
Let $\Gamma\subset \PSL(2,{\mathbb{C}})$ be a connected Lie group, $g\in
\Gamma$ a  loxodromic element and $U$ be a neighborhood of $g$
such that $hgh^{-1}g^{-1}=Id$ for each $h\in U$. Then
$hgh^{-1}g^{-1}=Id$ for each $h\in \Gamma$.
\end{lemma}


\begin{proof}
Define $W=\{h\in \Gamma:hgh^{-1}g^{-1}=Id\}$. We claim that $W$ is open and closed in $\Gamma$. For this,
let $\tau \in W$, then $\tau g^{-1}U$ is an open set that contains $\tau$ and it is contained in $W$. Hence $W$ is open.
Now, given  $\tau \in \overline{W}$,
notice that there exists a sequence $(\tau_n)\subset W$ such that
$\tau_n \xymatrix{\ar[r]_{n \rightarrow \infty}&}\tau$. We observe that
$$Id=\tau_ng\tau_n^{-1}g^{-1}
\xymatrix{\ar[r]_{n \rightarrow \infty}&} \tau g\tau^{-1}g^{-1}\,.$$ 
Then $\tau g\tau^{-1}g^{-1}=id\,.$
 Therefore $W$ is closed.  Since $\Gamma$ is connected, this implies  $W=\Gamma$.
\end{proof}


Now we have the following proposition, which completes Theorem \ref{t:eq3}.

\begin{proposition} \label{p:15}
Let $\Gamma\subset \PSL(2,\C)$ be  a non-discrete,  non-elementary
group  with no empty equicontinuity set. Then  $\Lambda_{\Gr}(\Gamma)$
is a circle in $\mathbb{P}^1_\mathbb{C}$.
\end{proposition}
\begin{proof}
Since $\PSL(2,{\mathbb {C}})$ is a Lie group we deduce that
$\overline{\Gamma}$ is a Lie group (see \cite{varadajan}). Let
${H}$ be  the connected component of the identity in
$\overline{\Gamma}$, so  ${H}$ is a connected and normal
subgroup of $\overline{\Gamma}$. By Theorem \ref{t:greenberg},
 (\ref{p:pgro}) of Proposition \ref{p:cr1}, and corollaries \ref{c:7}
and \ref{c:tran},  to finish  the proof  of (\ref{p:15}) we only need to consider  the following
cases:\vv

{\it Case} 1.- There exists exactly one point $p\in\mathbb{P}^1_{\mathbb{C}}$ such
that ${H}p=p$.- In this case, since ${H}$ is a normal subgroup, it follows that $\Gamma p=p$. Then  there exists a
purely elliptic Lie group $\mathcal{K}\subset {H}$ with
$dim_{\mathbb{R}}(\mathcal{K})=1$. Let $\gamma_0\in \Gamma$ be a
loxodromic element;  by Lemma \ref{l:13} we deduce that for
each loxodromic element $\gamma\in \Gamma$ one has  that
$Fix(\gamma)\subset \mathcal{K}(Fix(\gamma_0))$. Moreover by
Corollary \ref{c:12} we conclude that
$\mathcal{K}(Fix(\gamma_0))=\Lambda_{\Gr}(\Gamma)$. The result now
follows because $\mathcal{K}(Fix(\gamma_0))$ is a circle in
$\mathbb{P}^1_{\mathbb{C}}$.\vv

{\it Case} 2.- If the set $Card(\Lambda_{\Gr}(H))$ has at last two points, then   there exists an $H$-invariant  circle $\mathcal{C}$-  After conjugating with a M\"obius transformation, we can assume that
$\mathcal{C}=\overline{\mathbb{R}}$, thus there exists a loxodromic element  $\gamma_0\in {H}$  such
that $Fix(\gamma_0)\subset \mathbb{R}$. Let $p_1,p_2$ be the fixed points of
$\gamma_0$, then there exist  a neighborhood $W\subset
\mathbb{R}^{dim_{\mathbb{R}}({H})}$ of $0$ and   real analytic maps
$a,b,c,d:W\rightarrow \mathbb{C}$ such that $\phi:W\rightarrow {H}$
defined by: $$\phi(w)(z)=\frac{a(w)z+b(w)}{c(w)+d(w)}$$ is a chart of
$\gamma_0$ (that is $\phi(0)=\gamma_0$). Set   $F:W\times\mathbb{P}^{1}_{\mathbb{C}} \rightarrow \mathbb{P}^1_{\mathbb{C}}$ defined  by $F(w,z)=\phi(w)(z)-z$. Then
$\partial_{z}F(0,p_i)=g^\prime(p_i)-1 \neq 0$. By the Implicit
Function Theorem, there exist a neighborhood $W_0\subset W$ of $0$ and
continuous functions  $\tau_i:W_0\rightarrow \mathbb{C}$, $i=1,2$, such that $F(w,\tau_i(w))=0$ and $\tau_1(w)\neq \tau_2(w)$ for each $w\in W_0$. Hence
$\{\tau_1(w),\tau_2(w)\}=Fix(\phi(w))$ for all $w\in W_0 $. By Lemma \ref{l:14} we can
assume that  $\tau_1$ is non constant and   $\phi (W_0)$
contains only loxodromic elements. Thus $\tau_1(W_0)\subset
\Lambda_{\Gr}$ and contains  an open interval, which clearly
implies that $\Lambda_{\Gr}(\Gamma)=\overline{\mathbb{R}}$.
\end{proof}


\section*{Acknowledgments}
The authors would like to thank  IMPA (in Rio de Janeiro, Brazil) and its people, for their support and hospitality while working on  this paper; our special thanks to Paulo Sad and Jorge Vitorio Pereira, for very fruitful conversations. We are grateful to Juan Pablo Navarrete  for his valuable   comments  and 
to J. Francisco Estrada  for helpful discussions. And we are most indebted to the referee for his/her infinite patience and very  useful comments and suggestions.


\bibliographystyle{amsplain}

\vskip.5cm

\noindent
{Angel Cano and Jos\'e Seade: \\ Instituto de Matem\'aticas, Unidad Cuernavaca, \\ Universidad Nacional Aut\'onoma de M\'exico, \\ 
 Av. Universidad s/n. Col. Lomas de Chamilpa, \\
C.P. 62210, Cuernavaca, Morelos, M\'exico.}

\vskip.1cm

\noindent
email: angel@matcuer.unam.mx, $\,$ jseade@matcuer.unam.mx
\end{document}